%% file: SBM_SH.tex
\documentclass[3p,times,english]{elsarticle}
\usepackage{graphicx}
\usepackage{amssymb,amsmath,amsthm,amsbsy}
\usepackage{calligra,calrsfs,mathrsfs}
\usepackage{multirow}
\usepackage{cases}
\usepackage{pgfplots}
\usepackage{tikz}
\usepackage{tkz-euclide}
\usetikzlibrary[patterns]
\usepackage{hyperref}
\usepackage{algorithm}
\usepackage{algpseudocode}
\usepackage{dsfont}
\usepackage{caption, comment}
\usepackage{subcaption}
\usepackage{epstopdf}
\usepackage{epsfig}
\usepackage{array}
\newcolumntype{P}[1]{>{\centering\arraybackslash}p{#1}}
\usepackage{float}
\usepackage[T1]{fontenc}
\usepackage[latin9]{inputenc}
\usepackage{geometry}
\usepackage{color}
\usepackage{stmaryrd}
\usepackage{setspace}
\usepackage{amsthm}
\usepackage{textcomp}
\usepackage{hyperref}
\usepackage{soul,xcolor}
\usepgfplotslibrary{fillbetween}
\usepackage{ulem}
\usepackage{textcomp}
\setstcolor{red}
\usepackage{subcaption}
\allowdisplaybreaks
\makeatletter

\newcolumntype{M}[1]{>{\centering\arraybackslash}m{#1}}

\definecolor{darkgreen}{rgb}{0.0, 0.5, 0.0}

\input{format_definitions}

\theoremstyle{plain}

\newtheorem{prop}{Proposition}

\theoremstyle{remark}
\newtheorem{remark}{Remark}

\biboptions{sort&compress}
\allowdisplaybreaks

\makeatother

\usepackage{babel}

\journal{Journal of Computational Physics}

\begin{document}

\setlength{\unitlength}{1cm}

\begin{frontmatter}

\author[llnl]{Nabil M. Atallah\corref{ca}}
\ead{atallah1@llnl.gov}
\author[llnl]{Vladimir Z. Tomov}
\ead{tomov2@llnl.gov}
\author[duke]{Guglielmo Scovazzi}
\ead{guglielmo.scovazzi@duke.edu}
\address[llnl]{Lawrence Livermore National Laboratory, Livermore (CA), USA}
\address[duke]{Department of Civil and Environmental Engineering,
               Duke University, Durham, North Carolina 27708, USA}
\cortext[ca]{Corresponding author: Nabil M. Atallah}

\title{Weak Boundary Conditions for Lagrangian Shock Hydrodynamics: \\
A High-Order Finite Element Implementation on Curved Boundaries}

\begin{abstract}
We propose a new Nitsche-type approach for weak enforcement of normal
velocity boundary conditions for a Lagrangian
discretization of the compressible shock-hydrodynamics equations using
high-order finite elements on curved boundaries.
Specifically, the variational formulation is appropriately modified
to enforce free-slip wall boundary conditions, without perturbing the structure of the  function spaces used to represent the solution, 
with a considerable simplification with respect to traditional approaches. Total energy is conserved and the resulting
mass matrices are constant in time.
The robustness and accuracy of the proposed method are validated
with an extensive set of tests involving nontrivial curved boundaries.
\end{abstract}

\begin{keyword}
Lagrangian hydrodynamics; wall boundary conditions; curved boundaries;
high-order finite elements.
\end{keyword}

\end{frontmatter}

\section{Introduction}
\label{sec_intro}

This work is motivated by the need to perform Lagrangian hydrodynamics
simulations in domains with nontrivial and possibly curved boundaries.
The most common boundary condition (BC) in simulations of Lagrangian shock
hydrodynamics is the free-slip wall BC, $\bs{v} \cdot \bs{n} = 0$, where the fluid
particles are allowed to slip tangentially along the wall's surface.
Correct enforcement of these conditions is essential for accurately modeling
and simulating fluid behavior in confined geometries, such as containers or
channels, or when the flow must go around internal obstacles.
A few approaches have been proposed in the past for incompressible flows
\cite{engelman1982implementation, behr2004application},
but, to the best of our knowledge, there is no method that enforces wall BC robustly
for high-order FE Lagrangian simulations of compressible flows in
general curved geometries.

The starting point of this work is the high-order Finite Element
(FE) method of Dobrev {\it et. al.}~\cite{Dobrev2012}
(open-source version available at \cite{Laghos2019}).
This method has proven itself over the years,
demonstrating high-order accuracy, robust behavior for various problems,
good symmetry preservation, accurate capturing of the flow geometry;
the method has been used as a backbone of a next-generation multiphysics
simulation code \cite{Rieben2020}.
However, a major restriction of the original formulation is that it is only
applicable to straight boundaries, i.e., when the boundary normals $\bs{n}$ are
parallel to one of the coordinate axes.
This limitation is not just specific to the algorithms described in
\cite{Dobrev2012, Laghos2019, Rieben2020},
but to vast majority of the high-order Lagrangian and
Arbitrary Lagrangian-Eulerian (ALE) shock-hydrodynamic codes in the literature
\cite{Despres2019, Morgan2019, Abgrall2020, Gaburro2020, Waltz2021}.

Enforcing wall BC on general boundaries with FE can be done by strong
enforcement, i.e., by posing constraints on the linear system level, or by
weak enforcement, i.e., by adding certain penalty force integrals in
the variational formulation.
While strong enforcement is generally more accurate, it becomes
complex when the different velocity components need to couple, and requires
manipulations in the linear algebra operators.
Weak enforcement, on the other hand, provides more flexibility as it allows
unified treatment of different cases.
Furthermore, enforcing the wall BC weakly through penalty integrals allows
the use of numerical techniques like partial assembly and
matrix-free computations, which enable high performance on the latest
computer architectures \cite{Vargas2022, Kolev2021}.
For these reasons, the focus of this work is weak enforcement.

We propose a new Nitsche-type approach for weak enforcement of
free-slip wall boundary conditions.
Nitsche's method~\cite{nitscheweak} has been traditionally applied for the imposition of boundary conditions to elliptic and parabolic Partial Differential Equations (PDEs), and only more recently has been considered in the context of hyperbolic systems of PDEs~\cite{song2015nitsche,scovazzi2017velocity,song2018shifted}, associated with acoustics, waves in solids, and shallow water flows. 
Our method is inspired by the developments in~\cite{song2015nitsche,song2018shifted}, but aims at tackling the complex challenges associated with the strong nonlinearities of shock hydrodynamics, and introduces some important new ideas.
The variational form of the momentum equation is enhanced by two penalty terms,
affecting the mass matrix and the right-hand side, respectively.
These terms are chosen in a way that the mass matrix stays constant in time.
The variational form of the specific internal energy is adjusted in a similar
manner, incorporating terms to ensure the preservation of total energy.

This article is organized as follows:
Section~\ref{sec_prelim} introduces the equations of Lagrangian shock hydrodynamics;
Section~\ref{sec:sbm_equations} describes the principles and implementation of  weak slip wall boundary conditions;
Section~\ref{sec_results} present a suite of numerical experiments;
and Section~\ref{sec_concl} summarizes the conclusions and future work.

\input{1_Governing_Equations}

\input{2_SBM}

\input{3_NumericalResults}

\section{Conclusions}
\label{sec_concl}

We developed a framework for the weak imposition of slip wall boundary conditions in Lagrangian hydrodynamics simulations.
The advantage of the proposed approach was assessed in a number of computations involving domains with nontrivial and curved boundaries.
Through these simulations, we demonstrated the flexibility, robustness, and accuracy of the proposed approach, and consequently its overall superiority with respect to the strong imposition of similar boundary conditions.
Future work will be directed to extend these developments in the context of ALE methods, and to obtain an integrated approach to multi-material ALE/Lagrangian hydrodynamics \cite{Tomov2018}.
We will also leverage the new capability in order to develop novel methods for
shifted slip wall boundary conditions, a technology that will enable
hydrodynamics simulations in complex domains that won't have to be meshed
exactly.


\section*{Acknowledgments}
The authors of Duke University are gratefully thanking the generous support of Lawrence Livermore National Laboratories, through a Laboratory Directed Research \& Development (LDRD) Agreement.
This work performed under the auspices of the U.S. Department of Energy
by Lawrence Livermore National Laboratory under
Contract DE-AC52-07NA27344, LLNL-JRNL-853773.
Guglielmo Scovazzi has also been partially supported by the National Science Foundation, Division of Mathematical Sciences (DMS), under Grant 2207164.

\bibliographystyle{plain}
\bibliography{./SBM_SH}

\end{document}

%% file: format_definitions.tex
\newcommand{\ifcomment}{\iffalse}

\newcommand{\RN}[1]{\textup{\uppercase\expandafter{\romannumeral#1}}}

\newcommand{\bs}[1]{\boldsymbol{#1}}

\newcommand{\G}{\Gamma}

\newcommand{\Om}{\Omega}

\newcommand{\cT}{{\mathcal T}}

\DeclareRobustCommand{\tspsb}[2]{{%
		\m@th\ensuremath{%
			^{#1}%
			_{#2}%
		}%
}}
\DeclareRobustCommand{\tsp}[1]{{%
		\m@th\ensuremath{%
			^{#1}%
		}%
}}

\DeclareRobustCommand{\tsb}[1]{{%
		\m@th\ensuremath{%
			_{#1}%
		}%
}}

\newdefinition{rem}{Remark}

%% file: 1_Governing_Equations.tex
\section{General equations of Lagrangian shock hydrodynamics}
\label{sec_prelim}

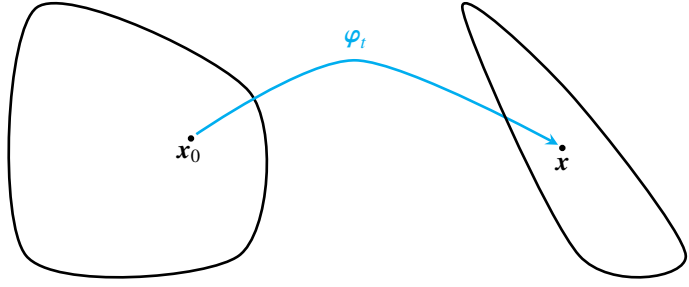
\begin{figure}
    \centering
    \begin{tikzpicture}[line width=1.0, scale = 1.4]
		\draw [black] plot [smooth cycle] coordinates{(-2.5, -3.7) (-.5, -4.5) (-0.6, -6.05)  (-2.6, -6.05)};
		\draw [cyan,-stealth] plot [smooth] coordinates{(-1.0, -4.9) (0.5, -4.2) (2.4, -5.0)};
		\draw [black] plot [smooth cycle] coordinates{(1.5, -3.7) (2.5, -4.5) (3.6, -6.05)  (2.6, -6.05)};
		\filldraw [black] (-1.05, -4.94) circle (0.5pt);
		\draw [black] (-1.07, -5.09) node{$\bs{x}_{0}$};
		\draw [cyan] (0.5, -4.0) node{$\bs{\varphi}_{t}$};
		\filldraw [black] (2.44, -5.03) circle (0.5pt);
		\draw [black] (2.44, -5.18) node{$\bs{x}$};
	\end{tikzpicture}
	\captionsetup{font=normalsize}
   \caption{ Sketch of the Lagrangian map $\bs{\varphi}_{t}$. \label{lagmap}}
\end{figure}

The classical equations of Lagrangian shock hydrodynamics govern the rate of change in position, momentum and energy of a compressible body of fluid, as it deforms. Let  $\Om_{0}$ and $\Om$ be open sets in $\mathbb{R}^{d}$ (where $d$ is the number of spatial dimensions.) The {\it motion}
\begin{subequations}
\label{eq:motion}
\begin{align}
\bs{\varphi}_{t}: \Om_{0}  &\rightarrow  \Om=\bs{\varphi}_{t}(\Om_{0}) \; , \\
\bs{x}_0 &\mapsto \bs{x}=\bs{\varphi}_{t}(\bs{x}_0) \;  , \quad \forall \bs{x}_0 \in \Om_{0}, \ t \geq 0 \; ,
\end{align}
\end{subequations}
maps the material coordinate $\bs{x}_{0}$, representing the initial position of an infinitesimal material particle of the body, to $\bs{x}$, the position of that particle in the current configuration (see Fig. \ref{lagmap}).
Here $\Om_{0}$ is the domain occupied by the body in its initial configuration, with boundary $\partial \Om_{0} = \G_{0}$ and outward-pointing boundary normal $\bs{n}_0$.
The transformation $\bs{\varphi}_{t}$ maps $\Om_{0}$ to $\Om$,
the domain occupied by the body in its current configuration, with
boundary $\partial \Om = \G$ and outward-pointing boundary normal $\bs{n}$.
Usually $\bs{\varphi}_{t}$ is a {\it smooth}, invertible map, and the {\it deformation gradient} $\bs{F} = \nabla_{\bs{x}_{0}} \bs{\varphi}_{t}$ and {\it deformation Jacobian determinant} $J=J(\bs{x}_{0},t) := \det(\bs{F})$ can be defined by means of the original configuration gradient $\nabla_{\bs{x}_{0}}$. 

In the domain $\Om$, we utilize the non-conservative form of the
Euler system \cite{benson92} and evolve the material position, density,
velocity, and specific internal energy:
\begin{subequations}
\label{eq:SH_Strong_updated_Lagrangian}
\begin{align}
\dot{\bs{x}} &= \bs{v} \; , \label{displ} \\
\rho \,  J           &= \rho_0 \; , \label{mass} \\
\rho \, \dot{\bs{v}} &= \rho \ \bs{b}  + \nabla_{\bs{x}} \cdot \bs{\sigma} \; , \label{mom} \\
\rho \, \dot{{e}}  &= \rho \ r + \nabla_{\bs{x}} \bs{v} \!:\! \bs{\sigma} - \nabla_{\bs{x}} \!\cdot \bs{q} \; . \label{ene} 
\end{align}
\end{subequations}
Here $\nabla_{\bs{x}}$ and $\nabla_{\bs{x}} \cdot$ are the current configuration gradient and divergence operators, respectively, and $\dot{(\cdot)}$ indicates the material, or Lagrangian, time derivative.
Furthermore $\rho_0$ is the reference (initial) density with
$\rho_{\max} := \sup_{\bs{x}_{0} \in \Om_{0}} \rho_{0}(\bs{x}_{0}) < \infty$, $\rho$ is the (current) density,  $\bs{v}$ is the velocity,
$\bs{b}$ is the body force (e.g., gravity),
$\bs{\sigma}$ is the symmetric  Cauchy stress tensor,
$r$ is the energy source term, and $\bs{q}$ is the heat flux.
Using index notation,
$\bs{\sigma}^T \!:\! \nabla_{\bs{x}}
 \bs{v} = \sigma_{ji}\  \partial_{x_i} \! v_j$, and
$\nabla_{\bs{x}} \bs{v} \!:\! \bs{\sigma} = \bs{\sigma} \!:\!
 \nabla_{\bs{x}} \bs{v} = \bs{\sigma}^T \!:\! \nabla_{\bs{x}} \bs{v}$,
since $\bs{\sigma}$ is symmetric.
We also denote by $E={e}+ \bs{v} \cdot \bs{v}/2$ the total energy per unit mass, 
the sum of the specific internal energy ${e}$ and
the kinetic energy $\bs{v} \cdot \bs{v} / 2$.
Obviously, $E$, ${e}$, $\bs{b}$, $r$ are measured per unit mass.

The system of equations~\eqref{eq:SH_Strong_updated_Lagrangian} are most commonly adopted in shock-hydrodynamics algorithms~\cite{benson92} and make use of the quasi-linear rather than the conservative form of the internal energy equation. The sum of the internal energy equation~\eqref{ene} and the kinetic energy equation (the product of \eqref{mom} by the velocity vector $\bs{v}$) yields the equation for the conservation of total energy. 

The system of equations~\eqref{eq:SH_Strong_updated_Lagrangian} completely defines the evolution of the system, once constitutive relationships for the stress $\bs{\sigma}$ and heat flux $\bs{q}$ are specified, together with appropriate initial and boundary conditions.

\subsection{Constitutive laws}
For a compressible inviscid fluid, the Cauchy stress $\bs{\sigma}$ reduces to an isotropic tensor, dependent only on the thermodynamic pressure, namely
\begin{equation}
\bs{\sigma} \ = \ - p \bs{I}_{d \times d} + \bs{\sigma}' \; , \label{stressV}
\end{equation}
where an equation of state of the type 
\begin{equation}
p = \mathcal{P}(\rho,{e}) \; , \label{eos}
\end{equation}
is assumed.
The deviatoric stress $\bs{\sigma}'$ is calculated from either a hypoelastic or hyperelastic constitutive model.
Mie-Gr\"uneisen equations of state are of the form \eqref{eos} with $\mathcal{P}(\rho,{e}) = \ f_1(\rho) + f_2(\rho) {e}$, and apply to materials such as compressible ideal gases, co-volume gases, high explosives, etc. (see \cite{menikoffplohr} for more details).
Ideal gases satisfy a Mie-Gr\"uneisen equation of state with $f_1=0$ and $f_2=(\gamma-1)  \rho$, namely
\begin{equation}
\mathcal{P}(\rho,{e}) =  (\gamma-1)  \rho  {e} \; ,
\end{equation}
where $\gamma$ is the exponent of the isentropic transformation of the gas. 
Consequently, the sound speed $c_{s}$ takes the form
\begin{align}
c_{s} = \sqrt{\gamma \, (\gamma-1) \, e} \, .
\end{align}


\subsection{Boundary conditions}
\label{bcs}
We assume that {\it slip boundary conditions}
are enforced on the entire boundary $\G$ given as
\begin{subequations}
\begin{align}
\label{wallbc} 
\bs{v} \cdot \bs{n}_{| \G} &=  0 \; ,
\\
\label{strescond}
\bs{\tau}_{i} \cdot (\bs{\sigma} \, \bs{n}_{| \G} ) &=  0 \; ,
\end{align}
\end{subequations}
where $\bs{n}_{| \G}$ is the outward-pointing normal to $\G$.
In particular~\eqref{wallbc} and~\eqref{strescond} enforce that the normal component of the velocity and the tangential component of the distributed traction force must vanish at the wall.


\subsection{General notation: Inner products, boundary functionals and norms}
\label{sec:gen}

Throughout the paper, we will denote by
\begin{align}
(v \, , w)_{\omega} = \int_{\omega} v \, w \, \text{d} \omega 
\qquad \text{and} \qquad
(\bs{v} \, , \bs{w})_{\omega} = \int_{\omega} \bs{v} \cdot \bs{w} \, \text{d} \omega 
\end{align} 
the $L^{2}(\omega)$ and the $(L^{2}(\omega))^{d}$ inner products on the interior of the domain $\omega \subset \Om$, and by
\begin{align}
\langle{v \, , w\rangle}_{\gamma} = \int_{\gamma} v \, w \, \text{d} \gamma
\qquad \text{and} \qquad
\langle\bs{v} \, , \bs{w}\rangle_{\gamma} = \int_{\gamma} \bs{v} \cdot \bs{w} \, \text{d} \gamma 
\end{align} 
a boundary functional on $\gamma \subset \Gamma$.

Let $L^{2}(\Om)$ be the space of square integrable functions on $\Om$. We will use the Sobolev spaces $H^m(\Om)=W^{m,2}(\Om)$ of index of regularity $m \geq 0$ (where $H^0(\Om)=L^{2}(\Om)$), equipped with the (scaled) norm
\begin{equation}
\|v \|_{H^{m}(\Om)} 
= \left( \| \, v \, \|^2_{L^2(\Om)} + \sum_{k = 1}^{m} \| \, l(\Om)^k  \bs{D}^k v \, \|^2_{L^2(\Om)} \right)^{1/2} \; ,
\end{equation}
where $\bs{D}^{k}$ is the $k$th-order spatial derivative operator and $l(A)=\mathrm{meas}_{d}(A)^{1/d}$ is a characteristic length of the domain $A$ with
$\text{meas}(A)$ the Lebesgue measure of the set $A$. As usual, we use a simplified notation for norms and semi-norms, i.e., we set $\| \, v  \, \|_{m,\Om}=\|\, v \, \|_{H^m(\Om)}$ and $| \, v \, |_{k,\Om}= 
\| \, \bs{D}^k v \,\|_{0,\Om}= \| \, \bs{D}^k v \, \|_{L^2(\Om)}$.

%% file: 2_SBM.tex
\section{A Nitsche approach to boundary conditions}
\label{sec:sbm_equations}
Let $\mathcal{T}^h$ be a family of admissible and shape-regular triangulations of $\Om$. 
We will denote by $h_T$ the circumscribed diameter of an element $T \in {\cT}^h$ and by $h$ the piecewise constant function in $\Om$ such that $h_{|T}=h_T$ for all $T \in {\cT}^h$.

%

\subsection{Discrete approximation spaces for the kinematic and thermodynamic variables}
\label{sec:functionspaces}
We rely on a semi-discrete formulation to derive the Nitsche weak variational formulation of the Euler equations in the Lagrangian reference frame for \eqref{displ}- \eqref{ene}. The discretization is determined by two finite dimensional functional spaces on the initial domain $\Om_{0}$:  
\begin{subequations}
	\begin{itemize}
	\item	$\mathcal{V}(\Om_{0} ) \subset  (H^{1}(\Om_{0}))^{d} $: the discrete space for the kinematic variables, with basis  $\{w_{a}\}_{a=1}^{n_{\mathcal{V}}} $ \,  .
	\item $\mathcal{E}(\Om_{0}) \subset  L_2(\Om_{0})$: the discrete space for the thermodynamic variables,  with basis  $\{\phi_{l}\}_{l=1}^{n_{\mathcal{E}}} $\,  .
	\end{itemize}
\end{subequations}
Note that we can define Lagrangian (moving) extensions of the kinematic and thermodynamic basis functions on $\Om(t)$ through the formulas $w_{a}(\bs{x},t) := w_{a} \circ \bs{\varphi}_{t}^{-1}(\bs{x})= w_{a}(\bs{x}_{0})$ and $\phi_{l}(\bs{x},t) := \phi_{l} \circ \bs{\varphi}_{t}^{-1}(\bs{x})= \phi_{l} (\bs{x}_{0})$, for $1 \leq a \leq n_{\mathcal{V}}$ and $1 \leq l \leq n_{\mathcal{E}}$, respectively.
These moving bases are constant along particle trajectories and therefore have zero material derivatives, that is,
\begin{align}
\dot{w}_{a} = 0 \, \text{ and } \, \dot{\phi}_{l} = 0 \, .
\end{align}
The spaces associated with the deformed domain $\Om(t)$ will be denoted by $\mathcal{V}(\Om(t))$ and $\mathcal{E}(\Om(t))$, respectively.
A mild restriction on the space $\mathcal{V}(\Om(t))$ is the requirement that 
\begin{align}
\bs{x}(\bs{x}_{0},t_{0})
 = \bs{x}_{0} \, , \qquad \forall \bs{x}_{0} \in \Om_{0} \, ,
\end{align}
expressing that we can represent exactly the initial geometry.
We discretize the position $\bs{x}(t)$ of the particle $\bs{x}_{0}$ at time $t$ using the expansion
\begin{align}
\label{eq:discrete_motion}
\bs{x}(\bs{x}_{0},t)
:=
\bs{\varphi}_{t}(\bs{x}_{0})
=
\mathrm{x}_{i;a}(t) \, w_{a}(\bs{x}_{0}) \, \bs{\zeta}_i \, ,
\end{align}
where $\mathrm{x}_{i;a}(t)$ is the time-dependent, $i$th coordinate of the position unknown of index $a$ and $\bs{\zeta}_i$ is the unit vector in the direction $i$, for $1 \leq i \leq d$.
The discrete velocity field corresponding to the motion~\eqref{eq:discrete_motion} is given by 
\begin{align}
\label{eq:velocity_motion}
\bs{v}(\bs{x}_{0},t)
:=
\dot{\bs{\varphi}}_{t}(\bs{x}_{0})
=
\mathrm{v}_{i;a}(t) \, w_{a}(\bs{x}_{0}) \, \bs{\zeta}_i \, ,
\end{align}
where $\mathrm{v}_{i;a}(t)=\dot{\mathrm{x}}_{i;a}(t)$. 
Note that we can also think of the velocity as a function on $\Om(t)$ with the expansion 
$$
\bs{v}(\bs{x}, t) 
=
\mathrm{v}_{i;a}(t) \, w_{a}(\bs{x},t) \, \bs{\zeta}_i \, .
$$
Using the same coordinates, but in the moving kinematic basis.
The thermodynamic discretization starts with the expansion of the internal energy
in the basis $\{\phi_{l}\}_{l=1}^{n_{\mathcal{E}}} $: 
\begin{align}
e(\bs{x}_{0},t) = \mathrm{e}_{l}(t)\,\phi_{l}(\bs{x}_{0}) \, .
\end{align}
Also the internal energy can also be expressed in the moving thermodynamic basis: $e(\bs{x},t) = \mathrm{e}_{l}(t)\,\phi_{l}(\bs{x},t) $.


\subsection{Semi-discrete mass conservation law}
Given an initial density field $\rho_{0}(\bs{x}_{0}) = \rho(\bs{x}_{0} , t_{0})$, we use the strong mass conservation principle \eqref{mass} to define the density pointwise at any time $t$,
\begin{equation}
\rho(\bs{x}_{0}, t) = \frac{\rho_{0}(\bs{x}_{0})}{J(\bs{x}_{0},t)} \, ,
\end{equation}
which implies that the mass in every Lagrangian volume is preserved exactly.

\subsection{Semi-discrete momentum conservation law}
We formulate the discrete momentum conservation equation by applying a Galerkin variational formulation to the continuous equation \eqref{mom}.
At any given time $t$, we multiply \eqref{mom} by a moving test function basis $\{ \bs{w}_{i;a} \}_{a=1}^{n_{\mathcal{V}}}$ constructed as $\bs{w}_{i;a} = w_a \, \bs{\zeta}_i$, where $w_a \in \mathcal{V}(\Om(t))$.

Integrating by parts over $\Om(t)$ and enforcing continuity of the stress $\bs{\sigma}$ across internal faces and~\eqref{strescond}  on $\G(t)$, we obtain
\begin{align}
\label{eq:mom_weak_form}
( \rho \, \dot{\bs{v}} , \bs{w}_{i;a})_{\Om(t)} 
+ ( \bs{\sigma} ,  \nabla_{\bs{x}}  \bs{w}_{i;a})_{\Om(t)} 
- \langle{ \bs{n} \cdot (\bs{\sigma}\, \bs{n}) \, \bs{n}, \bs{w}_{i;a} \rangle}_{\G(t)}
- ( \rho \, \bs{b},  \bs{w}_{i;a})_{\Om(t)} 
 =
0 \, .
\end{align}
Unlike \cite{Dobrev2012}, the boundary integral term will not vanish since \eqref{wallbc} is not embedded in the function space $\mathcal{V}(\Om(t))$. Expanding the velocity in terms of the moving velocity basis, observing that $\bs{\zeta}_i \cdot \bs{\zeta}_j = \delta_{ij}$ (the Kronecker delta tensor) and $b_i := \bs{b} \cdot \bs{\zeta}_i$, and using index notation give us
\begin{align}
\label{eq:mom_weak_form_index}
( \rho \, \dot{\mathrm{v}}_{i;b}(t) \, w_b , w_a  )_{\Om(t)} 
+ ( \sigma_{ik},  \nabla_{x_k} w_a )_{\Om(t)} 
- \langle{ (n_{k}\sigma_{kj} n_j) \, n_i, w_a  \rangle}_{\G(t)}
- ( \rho \, b_i,  w_a )_{\Om(t)} 
 =
0 \, .
\end{align}
We proceed now to the weak enforcement of~\eqref{wallbc},
by adding two penalty terms to~\eqref{eq:mom_weak_form}
(or~\eqref{eq:mom_weak_form_index}).
The first term is
\begin{align}
\label{eq:boundPenMass}
\langle{ \beta \, \rho \, c_{s} \, \bs{v} \cdot \bs{n} \, , \, w_{a} \, \bs{n} \cdot \bs{\zeta}_i \rangle}_{\G(t)} 
=
\langle{ \beta \, \rho \, c_{s} \, v_k n_k \, , \, w_{a} \, n_i \rangle}_{\G(t)} 
\, ,
\end{align}
which specifically enforces the slip boundary condition~\eqref{wallbc}.
Here $\beta = 20 \, C_{I} $ is a non-dimensional constant, while the choice of
$\rho \, c_{s}$ is standard to ensure the term has the correct units.
Following \cite{WarburtonHesthaven2003}, the penalty scale is chosen to
increase with the increase of the polynomial degree:
\begin{equation}
	C_{I} = \left\{
	\begin{array}{ll}
	(k+1)\,(k+d)\,d^{-1} & \mbox{for simplices in $d$ dimensions} \\
	(k+1)^2 & \mbox{for quadrilaterals and hexahedra,}
	\end{array}
	\right.
\end{equation}
where $k$ is the order of the spatial polynomial discretization for the velocity.
The second penalty term is
\begin{align}
\label{eq:boundPenMass_dert}
\langle{ \alpha_{0} \, \rho_{\max} \, L \,  \dot{\bs{v}} \cdot \bs{n}_{0} \, , \, w_{a} \, \bs{n}_{0} \cdot \bs{\zeta}_i \rangle}_{\G_{0}}
=
\dot{\mathrm{v}}_{j;b} \,
\langle{ \alpha_{0} \, \rho_{\max} \, L \,  w_b \, n_{0;j} \, , \, w_{a} \, n_{0;i} \rangle}_{\G_{0}}
 \, ,
\end{align}
which enforces the condition $\dot{\bs{v}}(\bs{x}_0,t) \cdot \bs{n}_{0}(\bs{x}_0)=0$, that is, that the acceleration in the initial configuration frame is orthogonal to the boundary normal in the initial configuration.
This condition is equivalent to $\mathrm{d}/\mathrm{d}t(\bs{v} \cdot \bs{n})=0$ for boundary surfaces that do not move in the normal direction, like the ones considered in this work.
In this case, $\bs{n}$ stays constant over time and in particular $\bs{n}_0=\bs{n}$, so that $\mathrm{d}/\mathrm{d}t(\bs{v} \cdot \bs{n})= \dot{\bs{v}} \cdot \bs{n}_{0} $.
In \eqref{eq:boundPenMass_dert}, $L$ is the perimeter
(unit of length) of the bounding box of $\Om_{0}$,
$\alpha_{0} = \beta \, L / J_{\scriptscriptstyle\square}^{1/d}$ is a non-dimensional constant, and $J_{\scriptscriptstyle\square} = \text{det} (\bs{F}_{\scriptscriptstyle\square})$ with $\bs{F}_{\scriptscriptstyle\square} = \nabla_{\bs{\xi}} \bs{X}$ is the deformation gradient of the mapping from the parent domain to the original configuration.
In choosing the scaling in \eqref{eq:boundPenMass_dert}, the two key points are dimensional consistency and time-independence. Hence, a natural choice would be $\beta \, \rho_{\max} \, L$. However, multiplying the latter scaling by the dimensionless ratio $  L / J_{\scriptscriptstyle\square}^{1/d}$ is critical to ensure that~\eqref{eq:boundPenMass_dert} does not decrease with mesh refinement and, at the same time, allows $\rho_{\max}$ to mimic, to a certain degree, the time-dependent pointwise density $\rho$ present in \eqref{eq:boundPenMass} (inversely proportional to $J$). 
Our tests indicate that this choice of $\alpha_{0}$ is critical to obtain accurate tangential motion at the boundary as the mesh refinement.

\noindent
Let $\mathbf{v} = \mathrm{v}_{j;b}$ be an $[ n_{\mathcal{V}} \times d   ]$-vector,
\begin{subequations}
	\begin{align}
	\label{eq:massMat}
\mathbf{M}_{\mathcal{V}}	= \mathrm{M}_{a_i b_j} 
	= 
	( \rho \, w_{a} , w_b \, \delta_{ij})_{\Om(t)} 
	+ \langle{ \alpha_{0} \, \rho_{\max} \, L \, w_{a} \, n_i \, , \, w_b \, n_j  \rangle}_{\G_{0}}
	\end{align}
	\text{be a $[ (n_{\mathcal{V}} \times d ) \times (n_{\mathcal{V}} \times d )  ]$-square matrix with $a, b \in [1 \, , \, n_{\mathcal{V}}] $ and $i,j \in [1 \, , \, d]$, }
\begin{align}
	\label{eq:Frhs}
\mathbf{F} = \mathrm{F}_{a_i \, l} 
=
( \sigma_{ik}, \phi_{l}\,  \nabla_{x_k} w_a )_{\Om(t)}  
- \langle{ (n_{k}\sigma_{kj} n_j) \, n_i, w_a  \, \phi_{l} \rangle}_{\G(t)}
+ 	\langle{ \beta \, \rho \, c_{s} \, \bs{v} \cdot \bs{n}, w_a  \, n_{i} \, \phi_{l} \rangle}_{\G(t)}
\end{align}
	\text{be a $[ (n_{\mathcal{V}} \times d ) \times n_{\mathcal{E}}   ]$-matrix with $1 \leq l \leq n_{\mathcal{E}}$, and}
\begin{align}
		\label{eq:Brhs}
	\mathbf{B}	= \mathrm{B}_{a_i} = ( \rho \, b_{i},  w_a)_{\Om(t)}
	\end{align}
	\text{be a $[ n_{\mathcal{V}} \times d   ]$-vector. }
\end{subequations}
Then the semidiscrete Lagrangian momentum equation can be  written in matrix vector form as
\begin{align}
\label{eq:final_mom_weak_form_mat}
\mathbf{M}_{\mathcal{V}} \frac{d\mathbf{v}}{dt} 
= 
-\mathbf{F} \, \bs{1}
+ \mathbf{B}
\, ,
\end{align}
where $\bs{1}$ is a $[ n_{\mathcal{E}} ]$-vector whose entries are all equal to one.

\subsection{Semi-discrete energy conservation law}
We formulate the discrete energy conservation equation by multiplying \eqref{ene} with a test function $\phi_l \in \mathcal{E}(\Om(t))$ and integrating over the domain $\Om(t)$
\begin{align}
\label{eq:ene_weak_form}
(\rho \, \dot{{e}} , \phi_{l})_{\Om(t)}  
- (  \bs{\sigma} : \nabla_{\bs{x}} \bs{v} , \phi_{l})_{\Om(t)} 
+ (\nabla_{\bs{x}} \cdot \bs{q} -  \rho \ r, \phi_{l})_{\Om(t)} 
= 0\, .
\end{align}
To ensure conservation of total energy, we add 
	\begin{align}
	\label{eq:boundPenItg}
 \langle{\bs{n} \cdot (\bs{\sigma}\, \bs{n}) \, \bs{v} \cdot \bs{n},  \, \phi_{l} \rangle}_{\G(t)} 
 	-\langle{ \beta \, \rho \, c_{s} \, (\bs{v} \cdot \bs{n})^{2}, \phi_{l} \rangle}_{\G(t)} \, .
	\end{align}
to the left hand side of~\eqref{eq:ene_weak_form}. Observe that these are residual terms, since they weakly enforce the boundary condition $\bs{v} \cdot \bs{n}=0$. 
Expanding the internal energy in terms of the moving thermodynamic basis and using tensor notation give us
\begin{align}
\label{eq:temp1_e}
\frac{d\mathrm{e}_{m}}{dt} ( \rho \, \phi_{m} , \phi_{l} )_{\Om(t)}
= 
(  \bs{\sigma} : \nabla_{\bs{x}} \cdot \bs{v} , \phi_{l})_{\Om(t)} 
+ (\rho \ r  - \nabla_{\bs{x}} \cdot \bs{q}   , \phi_{l})_{\Om(t)}  
- \langle{\bs{n} \cdot (\bs{\sigma}\, \bs{n})  \, \bs{v} \cdot \bs{n}, \phi_{l} \rangle}_{\G(t)} 
+\langle{ \beta \, \rho \, c_{s} \, (\bs{v} \cdot \bs{n})^{2}, \phi_{l} \rangle}_{\G(t)}
 \, .
\end{align}
Defining
\begin{subequations}
	\label{eq:discreteE}
	\begin{align}
	\label{eq:massMatE}
	\mathbf{M}_{\mathcal{E}}	= \mathrm{M}_{ml} 
	= 
	( \rho \, \phi_{m} , \phi_{l} )_{\Om(t)} 
	\end{align}
	\text{be a $[ n_{\mathcal{E}}  \times n_{\mathcal{E}}   ]$-square matrix with $1 \leq m, l \leq \, n_{\mathcal{E}}$ and } \,
	\begin{align}
	\label{eq:Rrhs}
	\mathbf{R} = \mathrm{R}_{l} 
	=
(\rho \ r - \nabla_{\bs{x}} \cdot \bs{q}, \phi_{l}) 
	\, ,
	\end{align}
	\text{be a $[  n_{\mathcal{E}}   ]$-vector,}
\end{subequations}
the semi-discrete energy conservation can  be written in matrix-vector form as
\begin{align}
\label{eq:final_ene_weak_form_mat}
\mathbf{M}_{\mathcal{E}} \frac{d\mathbf{e}}{dt} 
= 
\mathbf{F}^{\mathrm{T}} \mathbf{v}
+ \mathbf{R}
\, ,
\end{align}

\begin{rem}
	\label{rem:fixedMassMat}
	The mass matrices $\mathbf{M}_{\mathcal{V}}$ and $\mathbf{M}_{\mathcal{E}} $ are independent of time due to~\eqref{mass} and the fact that all the shape functions are independent of time.
	Namely,
	\begin{align*}
	\frac{d \mathbf{M}_{\mathcal{V}}}{dt} 
	= 
	\frac{d}{dt}	
	( \rho \, w_{a} , w_b \, \delta_{ij})_{\Om(t)} 
	+ 
	\frac{d}{dt} 
	\langle{ \alpha_{0} \, \rho_{\max} \, L \, w_{a} \, , \, w_b \, \delta_{ij} \rangle}_{\G_{0}}
	= 0 \, .
	\end{align*}
		\begin{align*}
	\frac{d \mathbf{M}_{\mathcal{E}}}{dt} 
	= 
	\frac{d}{dt}	
( \rho \, \phi_{m} , \phi_{l} )_{\Om(t)}
	= 0 \, .
	\end{align*}
\end{rem}


\subsection{Artificial viscosity operator}
\label{sec_visc}

To ensure a comprehensive presentation of the proposed method, this section details
the formulas associated with the artificial viscosity operator
$\bs{\sigma}_{\mathrm{art}}$, which align with those in \cite{Dobrev2012}.

The artificial viscosity tensor $\bs{\sigma}_{\mathrm{art}}$
is added to the semidiscrete equations
\eqref{eq:final_mom_weak_form_mat} and \eqref{eq:final_ene_weak_form_mat}
to regularize shock wave propagation.

This technique was originally introduced by Von Neumann and Richtmyer~\cite{VonNeumann1950}, whereby the discrete Euler equations are augmented with a diffusion term scaled by a special mesh dependent nonlinear coefficient $\mu$. In particular, we add 
$$
 ( \bs{\sigma}_{\mathrm{art}},  \nabla_{\bs{x}}\bs{w}_{i;a})_{\Om(t)} 
$$
in the momentum equation and
$$
- (\bs{\sigma}_{\mathrm{art} :\nabla_{\bs{x}} \bs{v} }, \phi_{l})_{\Om(t)} 
$$
in the energy equation.
From \cite{Dobrev2012} and \cite{StabShockHydroP12007}, we choose
\begin{align}
\bs{\sigma}_{\mathrm{art}} =  \mu_{s_{1}} \epsilon(\bs{v})
\end{align}
 with $\epsilon(\bs{v}) = \frac{1}{2} \left(\nabla \bs{v} + (\nabla \bs{v})^{T}\right)$ is the symmetrized velocity gradient and 
 \begin{align}
\mu_{s_{1}} = \rho \left(q_{2} \, l^{2}_{\bs{s}} \, | \Delta_{\bs{s}} \bs{v} | + q_{1} \, \psi_{0} \, \psi_{1} \, l_{\bs{s}} \, c_{s}\right) \, ,
 \end{align}
 where $q_{1}$ and $q_{2}$ are linear and quadratic scaling coefficients chosen as $1/2$ and $2$ respectively, $ \Delta_{\bs{s}} \bs{v}  $ is the directional measure of compression defined in section 6.1 of~\cite{Dobrev2012}, $l_{\bs{s}} = l_{\bs{s}}(x)$ is a directional length scale defined in the direction of $\bs{s}$ evaluated at a point $x$ defined in section 6.3 in~\cite{Dobrev2012}, $\psi_{1}$ is a compression switch which forces the linear term to vanish at points in expansion, and
$\psi_{0}$ is a vorticity switch that suppresses
the linear term at points where vorticity dominates the flow:
\begin{equation}
\psi_{1} =
  \begin{cases}
  1 & \Delta_{\bs{s}} \bs{v} < 0 \, , \\
  0 & \Delta_{\bs{s}} \bs{v} \geq 0 \, , 
  \end{cases}
  \quad \quad
\psi_{0} = \frac{|\nabla \cdot \bs{v}| }{\| \nabla \bs{v} \|} \, .
\end{equation}


\subsection{Conservation of total linear momentum and total energy}
\label{sec_conservation}

Linear momentum is conserved by the proposed numerical approach, up to
$O(h^{m+1})$, where $m$ is the regularity index from Section \ref{sec:gen}.
Premultiplying~\eqref{eq:final_mom_weak_form_mat}  by $\bs{c}_{i;a}$, where $\bs{c}_{i;a} = c_a \, \bs{\zeta}_i$ is a constant $[ n_{\mathcal{V}} \times d   ]$-vector and $c_a w_{a}= 1$, we obtain:
\begin{align}
\bs{c}_{i;a} \cdot \left(\mathbf{M}_{\mathcal{V}} \frac{d\mathbf{v}}{dt} \right)
&=
- \bs{c}_{i;a} \cdot \left(\mathbf{F} \, \bs{1} \right)
+ \bs{c}_{i;a} \cdot \mathbf{B}
\nonumber \\
\frac{d}{dt}\int_{\Om_{0}}  \rho_{0} \, \bs{v} \cdot \bs{\zeta}_i
&=
\int_{\Om(t)} \rho \, \bs{b} \cdot \bs{\zeta}_{i}
-\int_{\G(t)}(n_{k}\sigma_{kj} n_j) \, (\bs{n} \cdot \bs{\zeta}_{i}) 
-
\int_{\G_{0}}  \rho_{0} \, \underbrace{\dot{\bs{v}} \cdot \bs{n}_{0}}_{O(h^{m+1})} (\bs{n}_{0} \cdot \bs{\zeta}_{i})
+
\int_{\G(t)} \beta \, \rho \, c_{s} \, \underbrace{\bs{v} \cdot \bs{n}}_{O(h^{(m+1)})} \, (\bs{n} \cdot \bs{\zeta}_{i})
\end{align}
The first three terms represent the statement of balance of accelerations, internal forces and boundary forces, respectively,
which is the standard statement of conservation of momentum.
The last two terms represent the error due to weakly enforcing
$\bs{v} \cdot \bs{n} = 0$, which scales as $O(h^{m+1})$.

Total energy is conserved exactly.
Premultiplying~\eqref{eq:final_mom_weak_form_mat} by $\mathbf{v}^T$ and~\eqref{eq:final_ene_weak_form_mat} by $\mathbf{1}^T$, we obtain:
$$
\frac{d\mathbf{E}}{dt}
:=
\frac{d}{dt}
\left( 
\frac{1}{2} \mathbf{v}^T \mathbf{M}_{\mathcal{V}} \mathbf{v}
+
\mathbf{1}^T \mathbf{M}_{\mathcal{E}} \mathbf{e}
\right)
=
\mathbf{v}^T \mathbf{B} + \mathbf{R}
\, ,
$$
were we applied Remark~\ref{rem:fixedMassMat}.

\begin{remark}
Observe that the total energy defined above converges to the exact total energy as the numerical solution converges to the exact solution, similarly to the case of strong imposition of boundary conditions.
In particular, while the definition of the numerical internal energy does not seem to pose any particular problem, the definition of the numerical kinetic energy must be considered with care.
In fact, the numerical kinetic energy reads
\begin{align}
\frac{1}{2} \mathbf{v}^T \mathbf{M}_{\mathcal{V}} \mathbf{v}
\;=\;
\frac{1}{2}
( \rho \, \mathrm{v}_{i;b} \, w_b , w_a \mathrm{v}_{i;a} )_{\Om(t)} 
+
\frac{1}{2}
\langle{ \alpha_{0} \, \rho_{\max} \, L \, , \, (\bs{n}_{0} \cdot {\bs{v}})^2 \rangle}_{\G_{0}} \, .
\end{align}
The second term on the right hand side goes to zero as the grid is refined or the polynomial order is increased, since the boundary condition on $\G_{0}$ will be more and more accurately satisfied.
Then the proposed definition of the numerical kinetic energy is consistent with the infinite dimensional limit.
\end{remark}


\subsection{Time integration and fully discrete approximation}
\label{sec_time}

So far we have focused exclusively on the spatial discretization.
Now we discuss the discretization of the time derivatives in the nonlinear system of ODEs~\eqref{displ},~\eqref{eq:final_mom_weak_form_mat} and~\eqref{eq:final_ene_weak_form_mat}, obtained from the spatial discretization of the Euler equations.
In this section we consider a general high-order temporal discretization method and demonstrate its impact on the semidiscrete conservations laws. 

We adopt the same {\it{modified}} midpoint Runge-Kutta second-order scheme proposed in~\cite{caramana1998construction,barlow2008compatible,scovazzi2008multi,Dobrev2012}.  This choice of time integrator guarantees the proposed formulation conserves total energy without resorting to any staggered approach in time. 
Let $t \in \{t_{n} \}_{n=0}^{N_{t}}$ and associate with each moment in time, $t_{n}$, the computational domain $\Om^{n} \equiv \Om(t_{n})$.
Let $\mathrm{Y}=(\mathbf{v}; \mathbf{e}; \mathrm{x})$
be the hydrodynamic state vector.
We identify the quantities of interest defined on $\Om^{n}$ with a superscript $n$, denote by $\Delta t$ the time-step. The fully discrete numerical algorithm then reads
\begin{align*}
\mathbf{v}^{\mathrm{n+\frac{1}{2}}} 
&= 
\mathbf{v}^{\mathrm{n}} - \frac{\Delta t}{2} \mathbf{M}_{\mathcal{V}}^{-1} 
\left(\mathbf{F}^{n} \mathbf{1} - \mathbf{B}^{n} \right) \, ,
& 
\mathbf{v}^{\mathrm{n+1}} 
&= 
\mathbf{v}^{\mathrm{n}} - 
\Delta t \, \mathbf{M}_{\mathcal{V}}^{-1} 
\left(\mathbf{F}^{n+\frac{1}{2}} \mathbf{1} - \mathbf{B}^{n+\frac{1}{2}} \right) \, ,
\\
\mathbf{e}^{\mathrm{n+\frac{1}{2}}} 
&= 
\mathbf{e}^{\mathrm{n}} + \frac{\Delta t}{2} \,  \mathbf{M}_{\mathcal{E}}^{-1} 
\left( (\mathbf{F}^{n})^{\mathrm{T}} \mathbf{v}^{n+\frac{1}{2}} + \mathbf{R}^{n} \right) \, ,
  & 
\mathbf{e}^{\mathrm{n+1}} 
&= 
\mathbf{e}^{\mathrm{n}} + \Delta t \, \mathbf{M}_{\mathcal{E}}^{-1} 
\left( (\mathbf{F}^{n+\frac{1}{2}})^{\mathrm{T}}  \bar{\mathbf{v}}^{n+\frac{1}{2}} + \mathbf{R}^{n+\frac{1}{2}} \right) \, ,
\\
\mathrm{x}^{\mathrm{n+\frac{1}{2}}} 
&= 
\mathrm{x}^{\mathrm{n}}
+ \frac{\Delta t}{2} \, \mathbf{v}^{n+\frac{1}{2}} \, ,
&    
\mathrm{x}^{\mathrm{n+1}} 
 &= 
\mathrm{x}^{\mathrm{n}}
+ \Delta t \, \bar{\mathbf{v}}^{n+\frac{1}{2}} \, ,
\end{align*}
where $\mathbf{F}^{\mathrm{k}} = \mathbf{F}(\mathrm{Y}^{\mathrm{k}})$ and $\bar{\mathbf{v}}^{n+\frac{1}{2}} = (\mathbf{v}^{\mathrm{n+1}}+\mathbf{v}^{\mathrm{n}})/2$ \, . \\

\begin{prop}
	\label{prop:discreteTotalMomCons}
	The RK2-average scheme described above conserves the discrete total linear momentum up to $O(h^{m+1})$.
	\begin{proof}
	Let $\bs{c}_{i;a} = c_a \, \bs{\zeta}_i$ is a constant $[ n_{\mathcal{V}} \times d   ]$-vector and $c_a w_{a}= 1$, then the change in linear momentum (LM) can be expressed as
			\begin{subequations}
				\begin{align*}
				LM^{n+1} - LM^{n} = \bs{c}_{i;a} \cdot  \left( \mathbf{M}_{\mathcal{V}} (\mathbf{v}^{\mathrm{n+1}}  - \mathbf{v}^{\mathrm{n}}  \right)
				= - \Delta t \, \bs{c}_{i;a} \cdot 
				\left(\mathbf{F}^{n+\frac{1}{2}} \mathbf{1} - \mathbf{B}^{n+\frac{1}{2}} \right) 
				 & = 
				 \Delta t \, \int_{\Om(t)^{n+\frac{1}{2}}} (\rho \, \bs{b})^{n+\frac{1}{2}} \cdot \bs{\zeta}_{i}
				 \nonumber \\
				 & \phantom{=}
				- \Delta t \, \int_{\G(t)^{n+\frac{1}{2}}}(n_{k}\sigma_{kj} n_j)^{n+\frac{1}{2}} \, (\bs{n}^{n+\frac{1}{2}} \cdot \bs{\zeta}_{i}) 
				\nonumber \\
				&  \phantom{=}
				+\Delta t \,\int_{\G(t)^{n+\frac{1}{2}}} \beta \left (\rho \, c_{s} \right)^{n+\frac{1}{2}}  {\underbrace{(\bs{v} \cdot \bs{n} )}_{O(h^{m+1})}}^{n+\frac{1}{2}}  (\bs{n}^{n+\frac{1}{2}} \cdot \bs{\zeta}_{i})
				\end{align*}
			\end{subequations}
	 The first two terms on the right hand side above are the same as the ones for a strong boundary condition enforcement and the last term is the $O(h^{m+1})$ error from weakly enforcing $\bs{v} \cdot \bs{n} = 0$.
	\end{proof}
\end{prop}

\begin{prop}
	\label{prop:discreteTotalEneCons}
The RK2-average scheme described above conserves the discrete total energy exactly.
\begin{proof}
The change in kinetic energy (KE) and internal energy (IE) can be expressed as
	\begin{subequations}
		\begin{align*}
		KE^{n+1} - KE^{n} = (\mathbf{v}^{\mathrm{n+1}}  - \mathbf{v}^{\mathrm{n}} ) \cdot  \left( \mathbf{M}_{\mathcal{V}} \, \bar{\mathbf{v}}^{n+\frac{1}{2}}  \right)
		= - \Delta t \left(\mathbf{F}^{n+\frac{1}{2}} \mathbf{1}  - \mathbf{B}^{n+\frac{1}{2}}  \right) \cdot \bar{\mathbf{v}}^{n+\frac{1}{2}} 
		\end{align*}
			\begin{align*}
		IE^{n+1} - IE^{n} 
		= \mathbf{1} \cdot  \left( \mathbf{M}_{\mathcal{E}} \, (\mathbf{e}^{\mathrm{n+1}} - \mathbf{e}^{\mathrm{n}})  \right)
		= \Delta t \, \mathbf{1} \cdot
		\left( (\mathbf{F}^{n+\frac{1}{2}})^{\mathrm{T}}  \bar{\mathbf{v}}^{n+\frac{1}{2}}+ \mathbf{R}^{n+\frac{1}{2}}  \right)
		\end{align*}
	\end{subequations}
Thus, the change in discrete total energy $KE^{n+1} + IE^{n+1}  - KE^{n} - IE^{n} = \Delta t 
\left( \mathbf{B}^{n+\frac{1}{2}} \cdot \bar{\mathbf{v}}^{n+\frac{1}{2}} + \mathbf{1} \cdot \mathbf{R}^{n+\frac{1}{2}}  \right)$ is due to the presence of body force, source and heat flux terms. If $r=0$ and $\bs{b} = \bs{q} = 0$, we see that the discrete total energy is conserved exactly, namely:
$KE^{n+1} + IE^{n+1}  = KE^{n} + IE^{n} $.
\end{proof}
\end{prop}

\begin{remark}
The RK2-average scheme described above has similar properties to the one described in~\cite{caramana1998construction,barlow2008compatible,scovazzi2008multi,love2009angular,Dobrev2012}.
In particular, it conserves exactly also angular momentum in the limit of a large number of corrector passes~\cite{love2009angular}.
The RK2-average scheme can be extended to higher orders of time integration
through the work presented in \cite{Sandu2021}.
\end{remark}

%% file: 3_NumericalResults.tex
\section{Numerical Results}
\label{sec_results}

We consider the standard shock hydrodynamic benchmark of the Sedov explosion~\cite{Sedov} in various 2D and 3D domains. In all test cases, the Sedov problem consists of an ideal gas ($\gamma$ = 1.4) with a delta function source
of internal energy deposited at the origin such that the total energy $E_{total}$ = 1. The sudden release of the energy creates an expanding shock wave, converting the initial
internal energy into kinetic energy. The delta function energy source is approximated
by setting the internal energy $e$ to zero in all degrees of freedom except at the origin
where the value is chosen so that the total internal energy is 1. In all of the tests, we enforce $\bs{v} \cdot \bs{n} = 0$ on all boundaries. 
Note that the density plots are in logarithmic scale.

All simulations are performed in a customized version of the open-source
Laghos proxy application~\cite{Laghos2019}, which is based on the
MFEM finite element library~\cite{MFEM2021}. 

\subsection{Two-dimensional Sedov explosion in a square}
\label{planar_sedov_2d}
\begin{figure}[tb]
	\centering
	\begin{tabular}{ccc}
		\multicolumn{3}{c}{Velocity}  \\
		\includegraphics[width=0.30\linewidth]{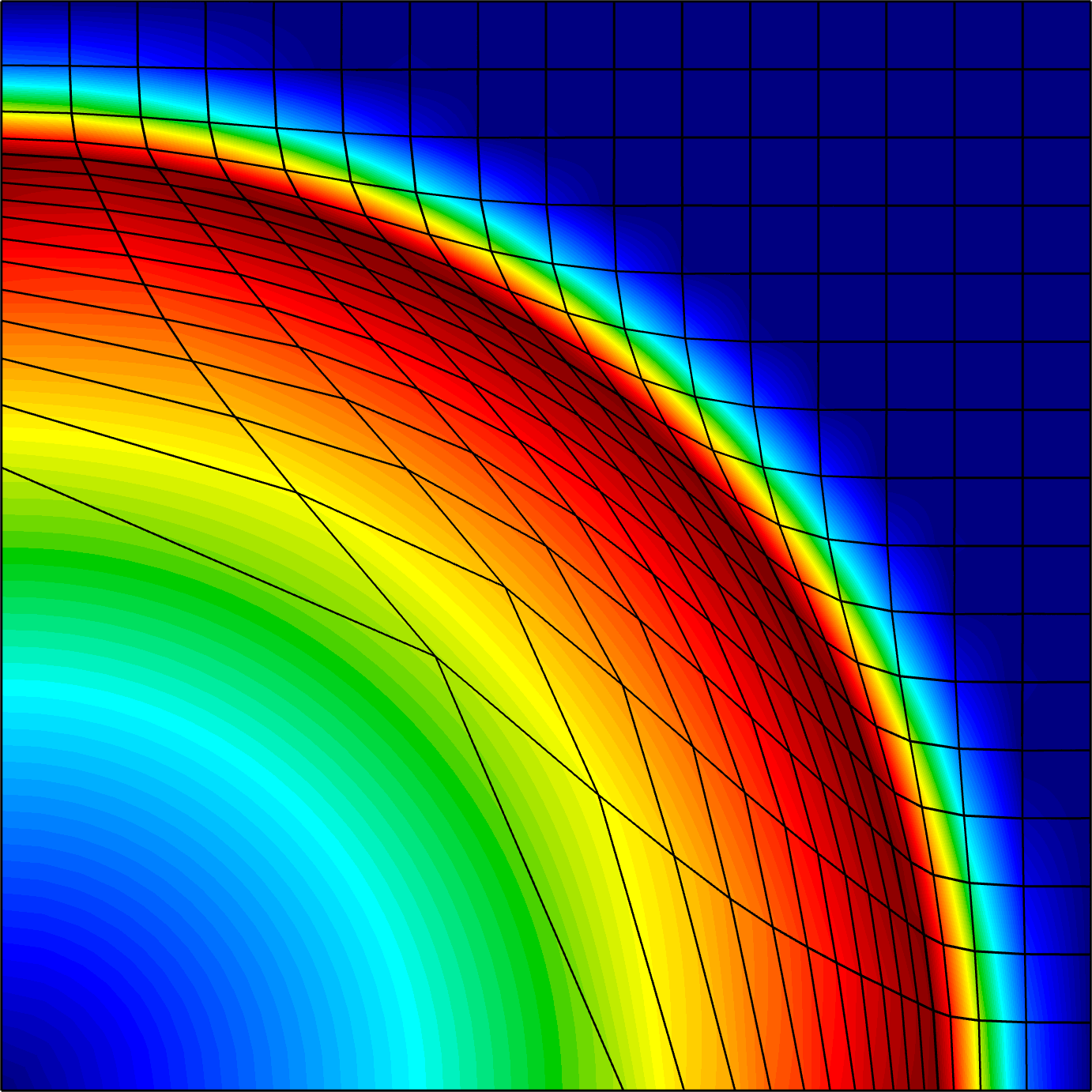} 
		\qquad & \qquad  
		\includegraphics[width=0.30\linewidth]{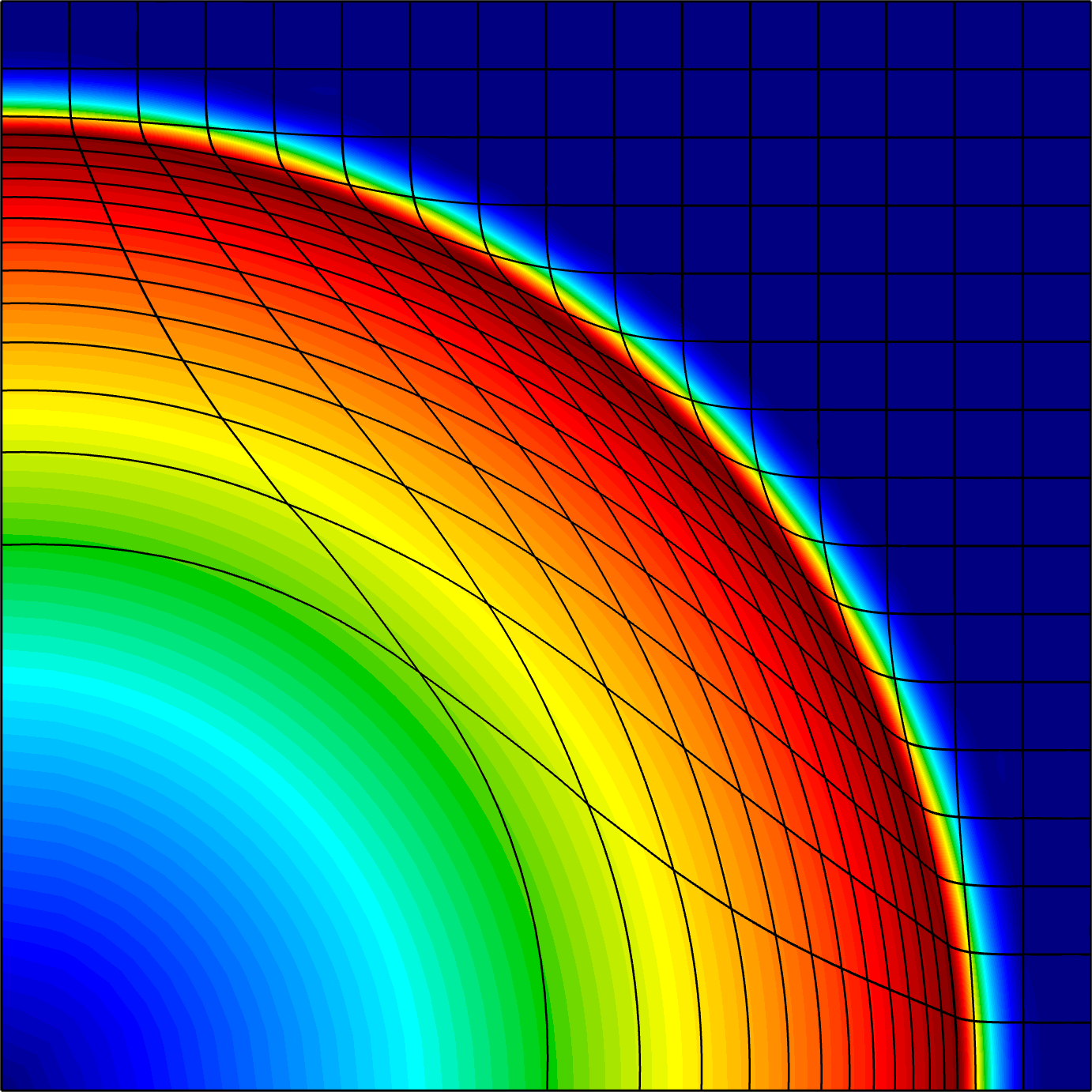}
		\qquad & \qquad  
		\includegraphics[width=0.30\linewidth]{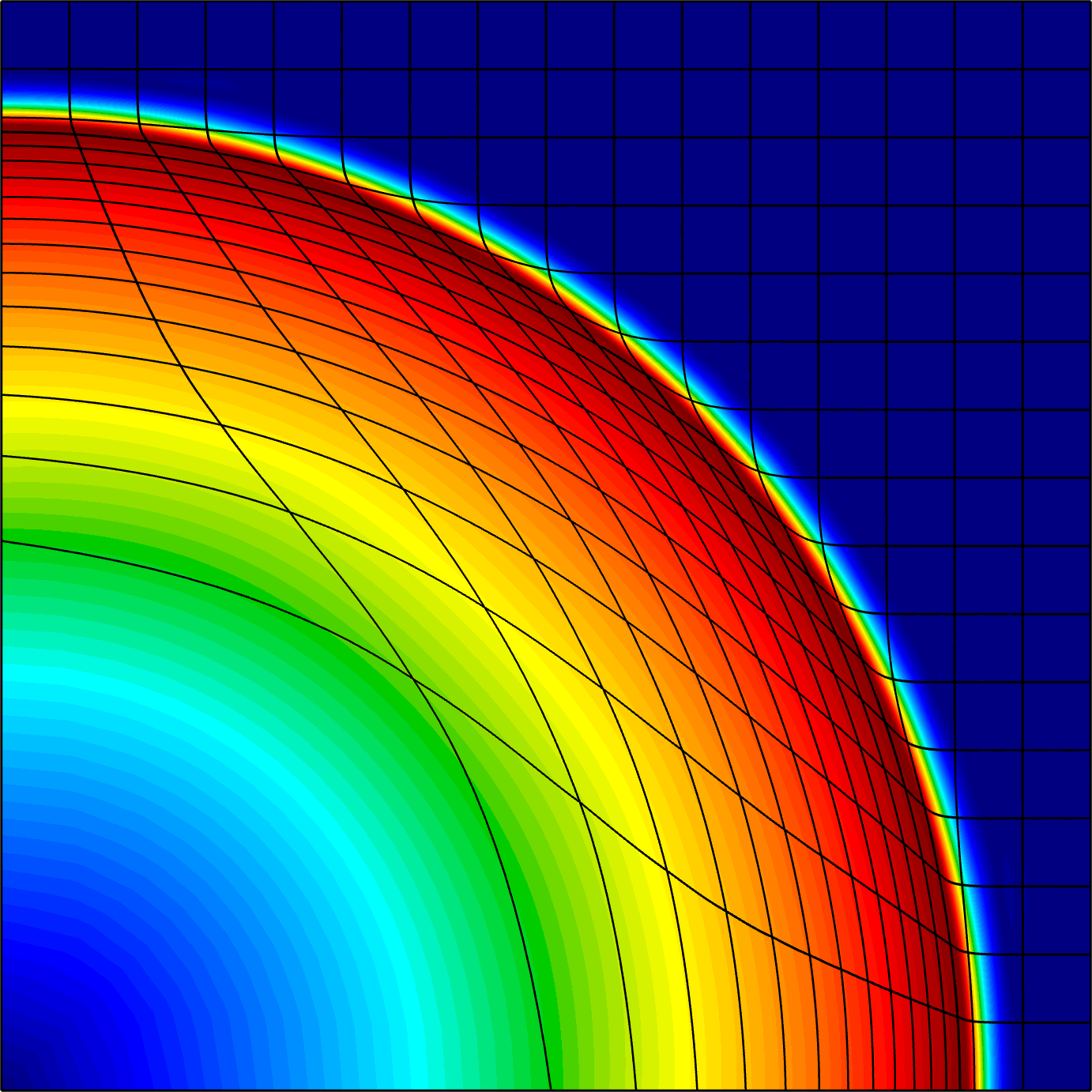} \\[.4cm]
		\multicolumn{3}{c}{Density}  \\
		\includegraphics[width=0.30\linewidth]{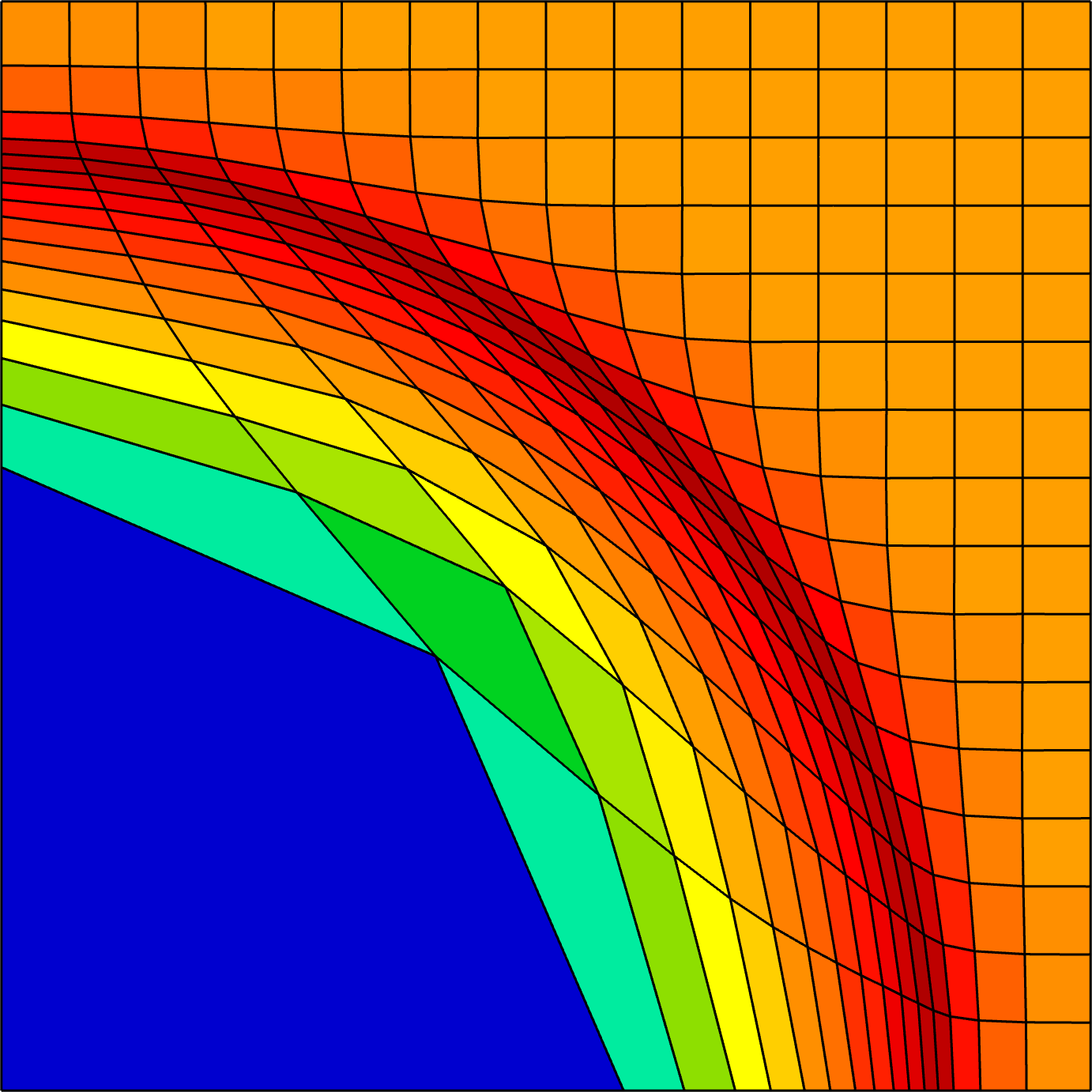} 
		\qquad & \qquad  
		\includegraphics[width=0.30\linewidth]{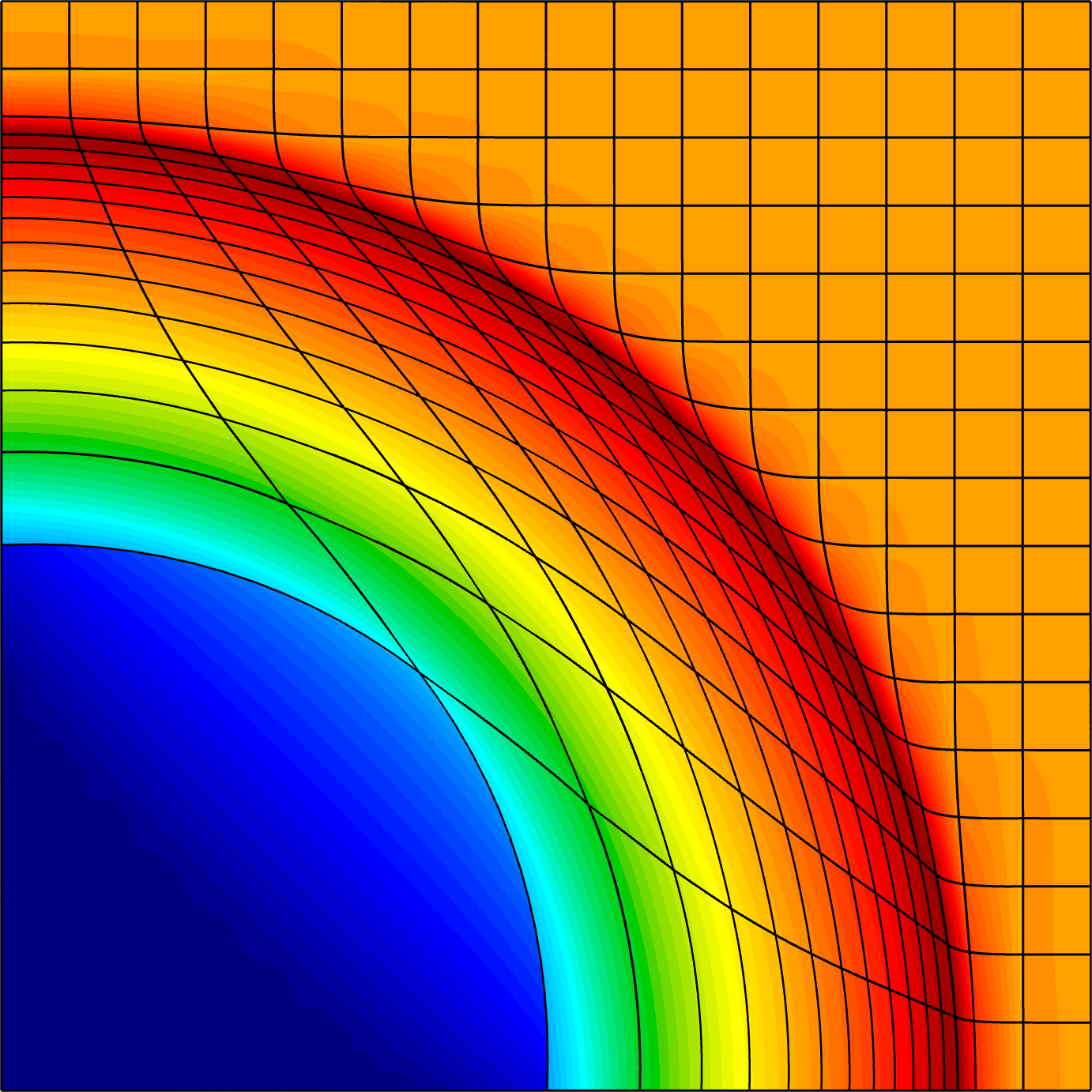}
		\qquad & \qquad  
		\includegraphics[width=0.30\linewidth]{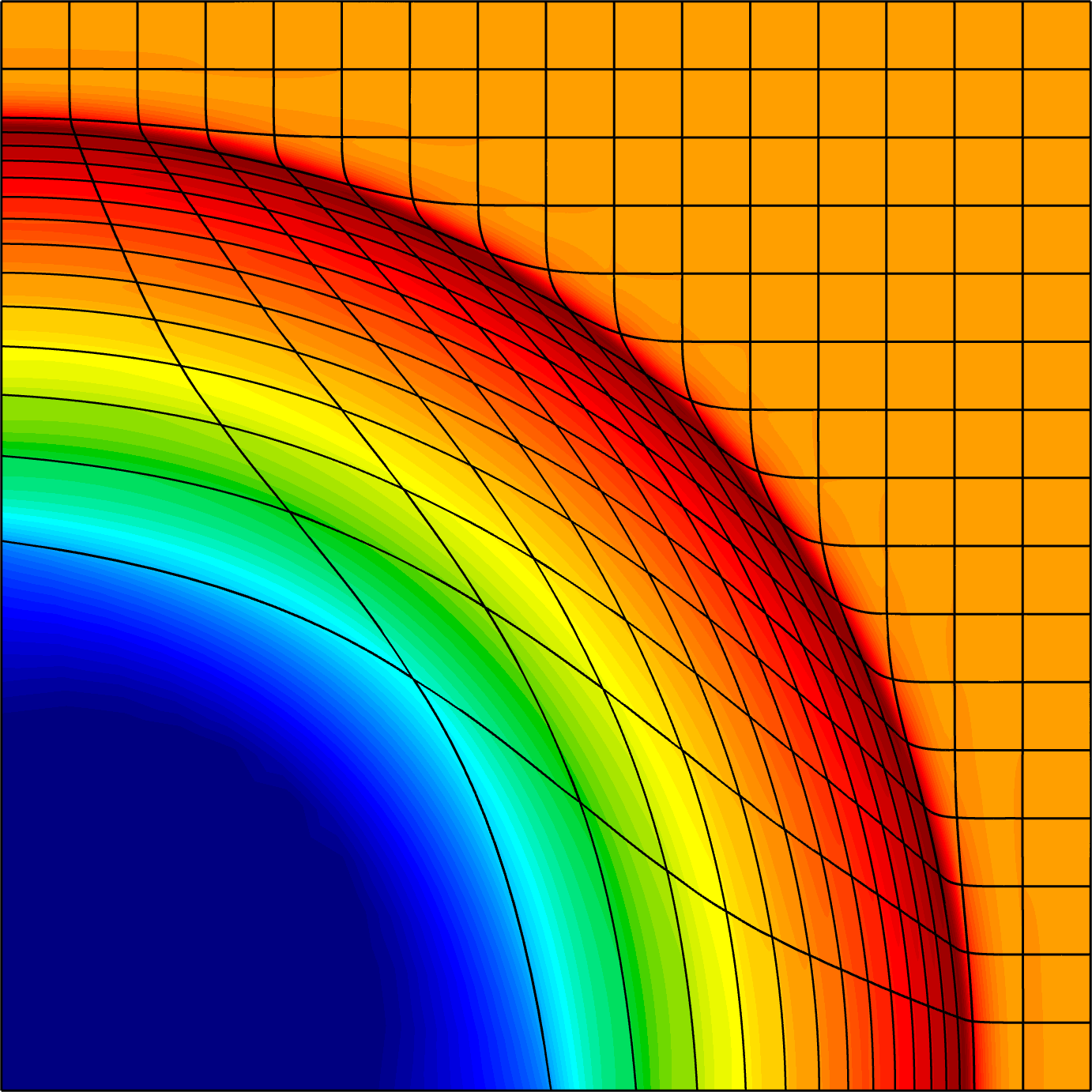}  \\[.4cm]
		\multicolumn{3}{c}{Mesh Deformation}  \\
		\includegraphics[width=0.30\linewidth]{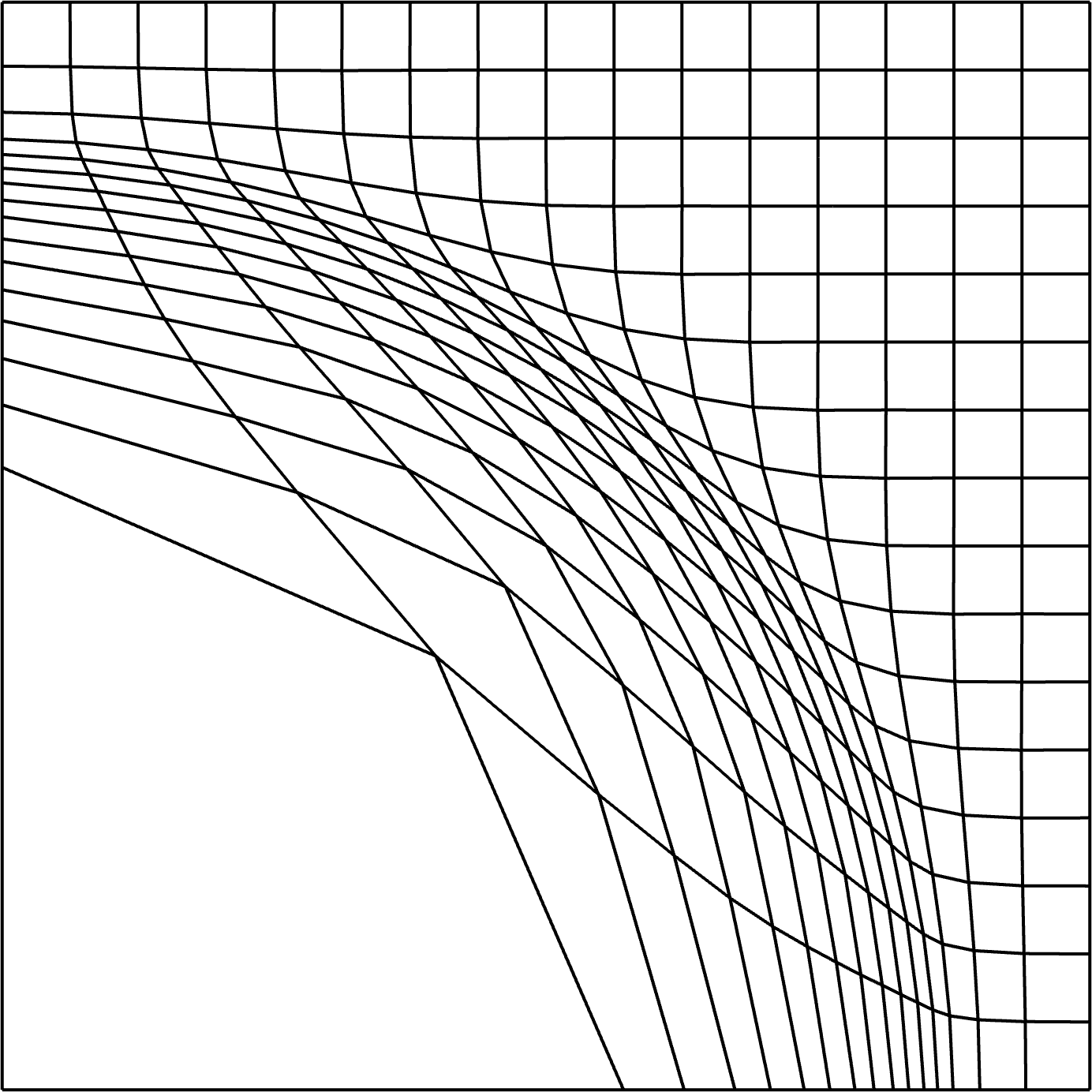} 
		\qquad & \qquad  
		\includegraphics[width=0.30\linewidth]{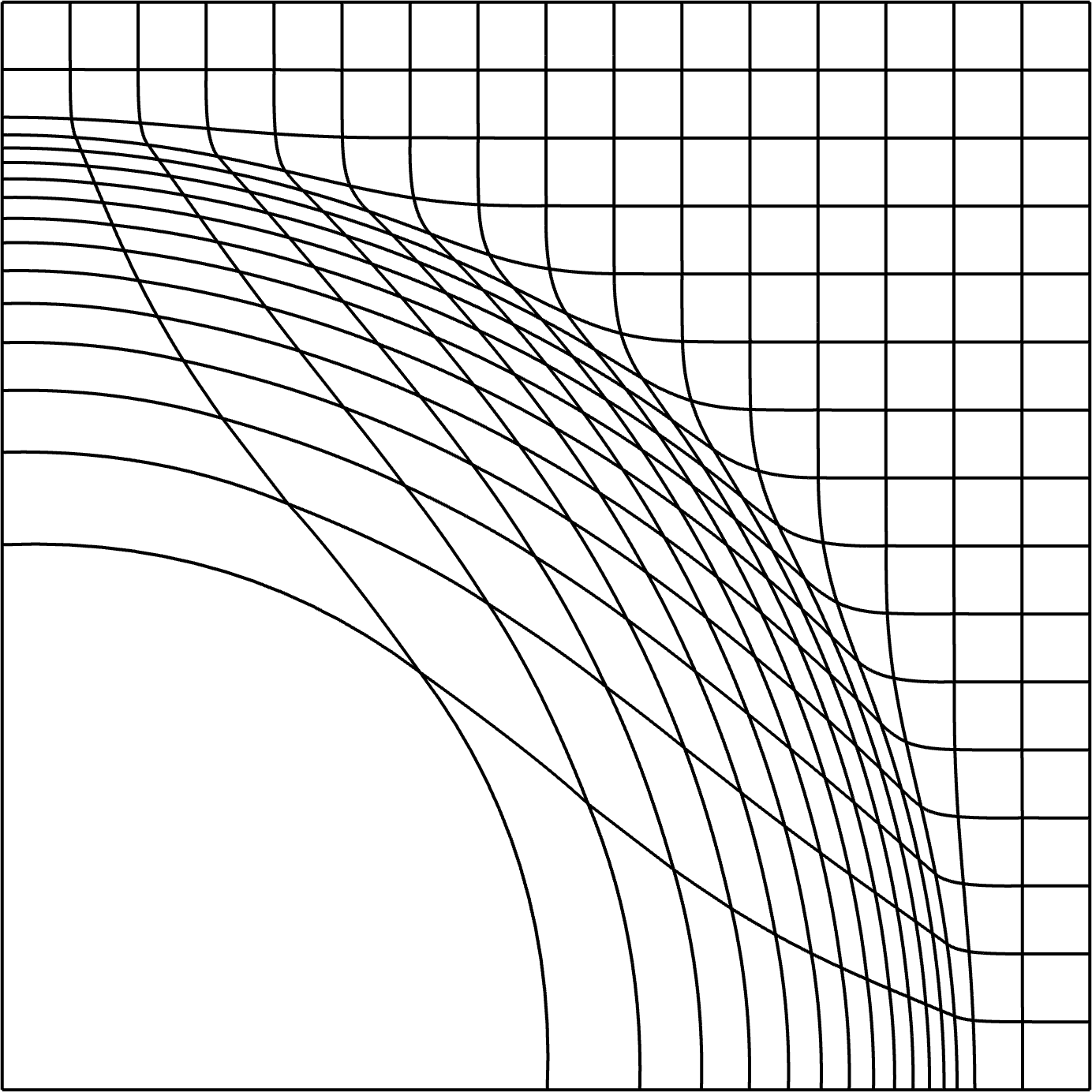}
		\qquad & \qquad  
		\includegraphics[width=0.30\linewidth]{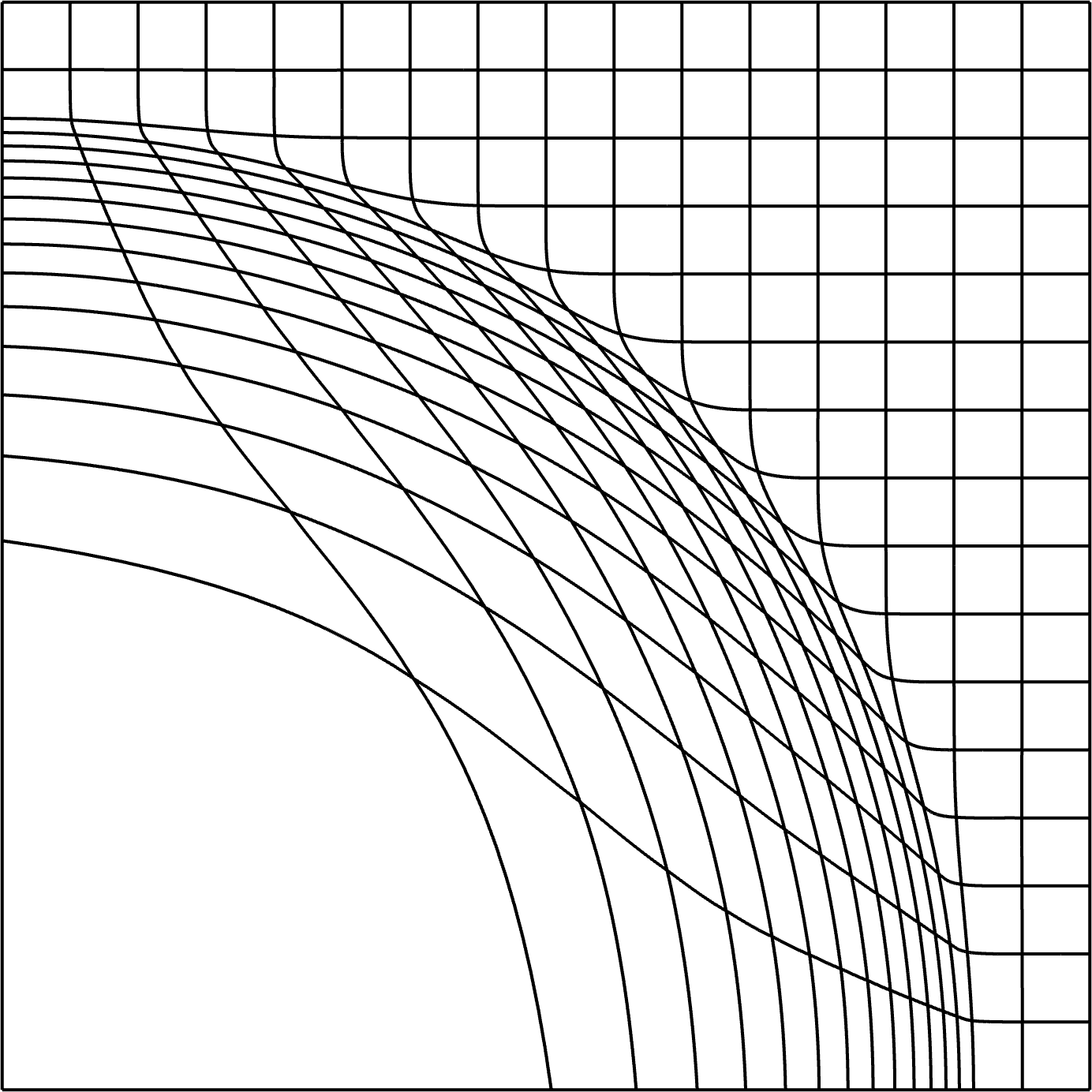}  \\[.0cm]
		$Q_{1}-Q_{0}$ 
		&
		$Q_{2}-Q_{1}$
		&
		$Q_{3}-Q_{2}$
	\end{tabular}
	\caption{Plots of the velocity and density fields in addition to the mesh deformation for the planar Sedov test using $Q_{1}-Q_{0}$, $Q_{2}-Q_{1}$, $Q_{3}-Q_{2}$ velocity-energy pairs.}
	\label{fig:PlanarSedovresults}
\end{figure}

\begin{figure}[tb]
	\centering
	\begin{subfigure}[t]{.4\linewidth}
	\centering	\includegraphics[width=1.0\linewidth]{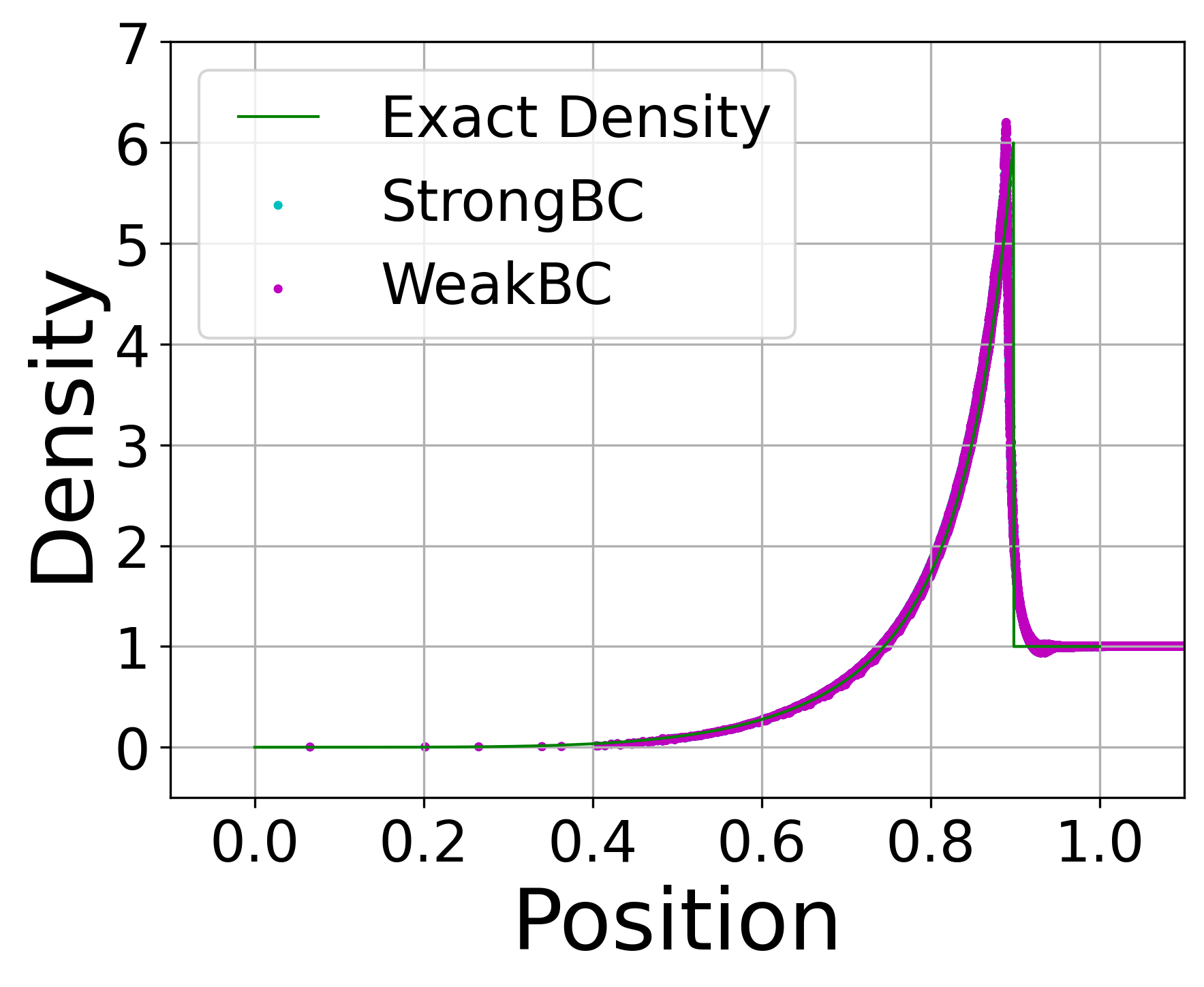}
		\caption{Comparison between the shock front position between strong and weak boundary enforcement with the exact solution.}
					\label{fig:SedovShockFrontP1}
		\end{subfigure}
		\qquad \qquad
		\begin{subfigure}[t]{.4\linewidth}
		\centering \includegraphics[width=1.0\linewidth]{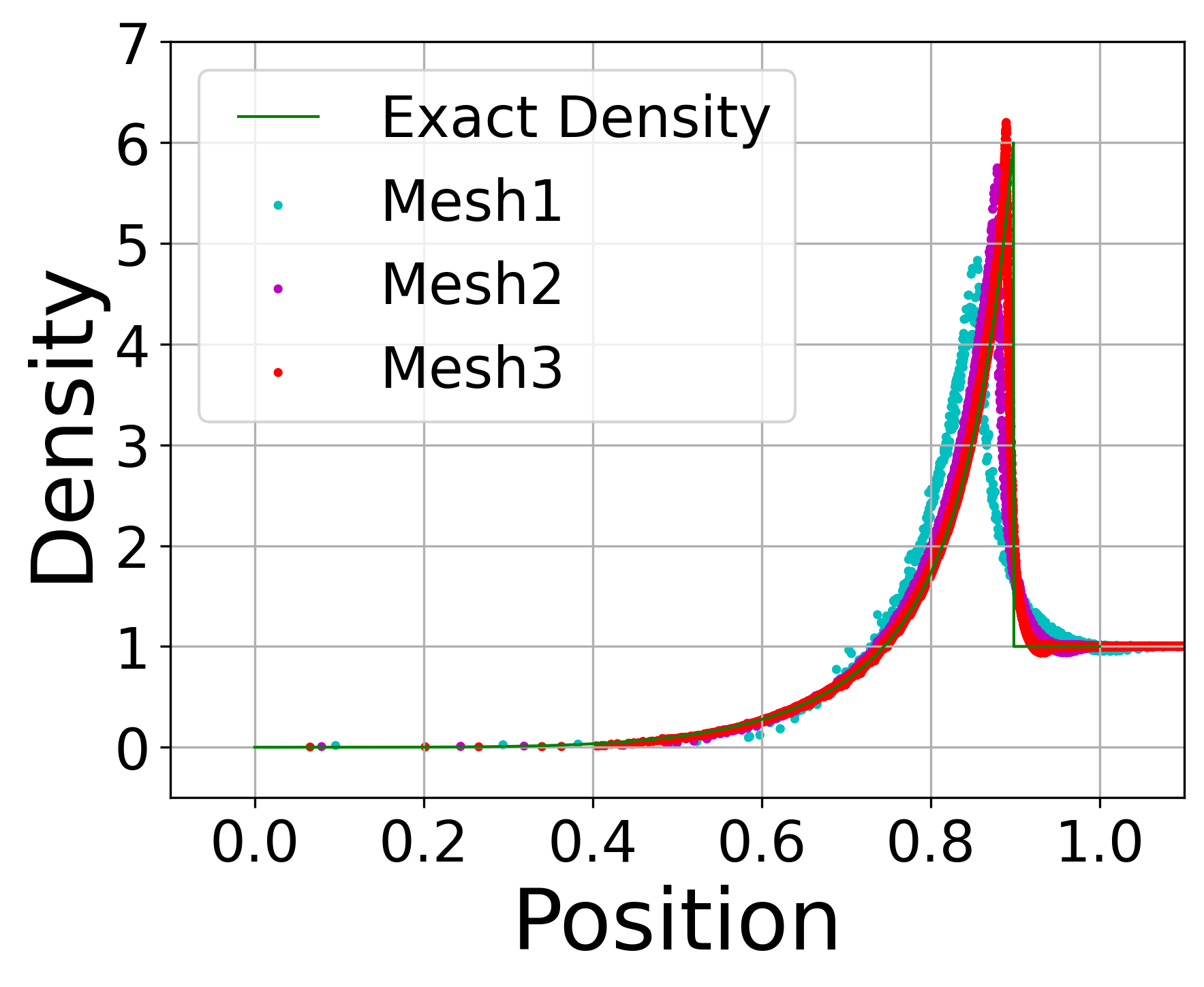}
		\caption{Convergence study for the shock front with weak boundary condition enforcement.}
						\label{fig:SedovShockFrontP2}
\end{subfigure}
	\caption{Plots comparing the shock front between strong and weak boundary condition enforcement (left) and a convergence study for weak boundary condition enforcement (right) for a Q2-Q1 velocity-energy pair.}
	\label{fig:SedovShockFront}
\end{figure}

We consider a $[0, 1] \times [0, 1]$ domain and a final time $t=0.8$. In Figure~\ref{fig:PlanarSedovresults}, we show plots of the velocity and density fields in addition to the mesh deformation at the final time of $t = 0.8$ for the $Q_{1}-Q_{0}$, $Q_{2}-Q_{1}$, $Q_{3}-Q_{2}$ velocity-energy pairs. 
As it is apparent, the weak wall boundary produce solutions that are indistinguishable from those obtained with strong enforcement which are smooth and without any unphysical oscillations. A plot of the shock front locations using a weak and strong boundary enforcement for a $Q_{2}-Q_{1}$ discretization in Figure~\ref{fig:SedovShockFrontP1} shows that they are indistinguishable from one another. Finally, in Figure~\ref{fig:SedovShockFrontP2} we see that shock front location attained from weak boundary condition enforcement converges to the exact location with mesh refinement.


\subsection{Two-dimensional Sedov explosion in a trapezoid}
\label{trapezoid_sedov_2d}
\begin{figure}[tb]
	\centering
	\begin{tabular}{ccc}
		\multicolumn{3}{c}{Velocity}  \\
		\includegraphics[width=0.30\linewidth]{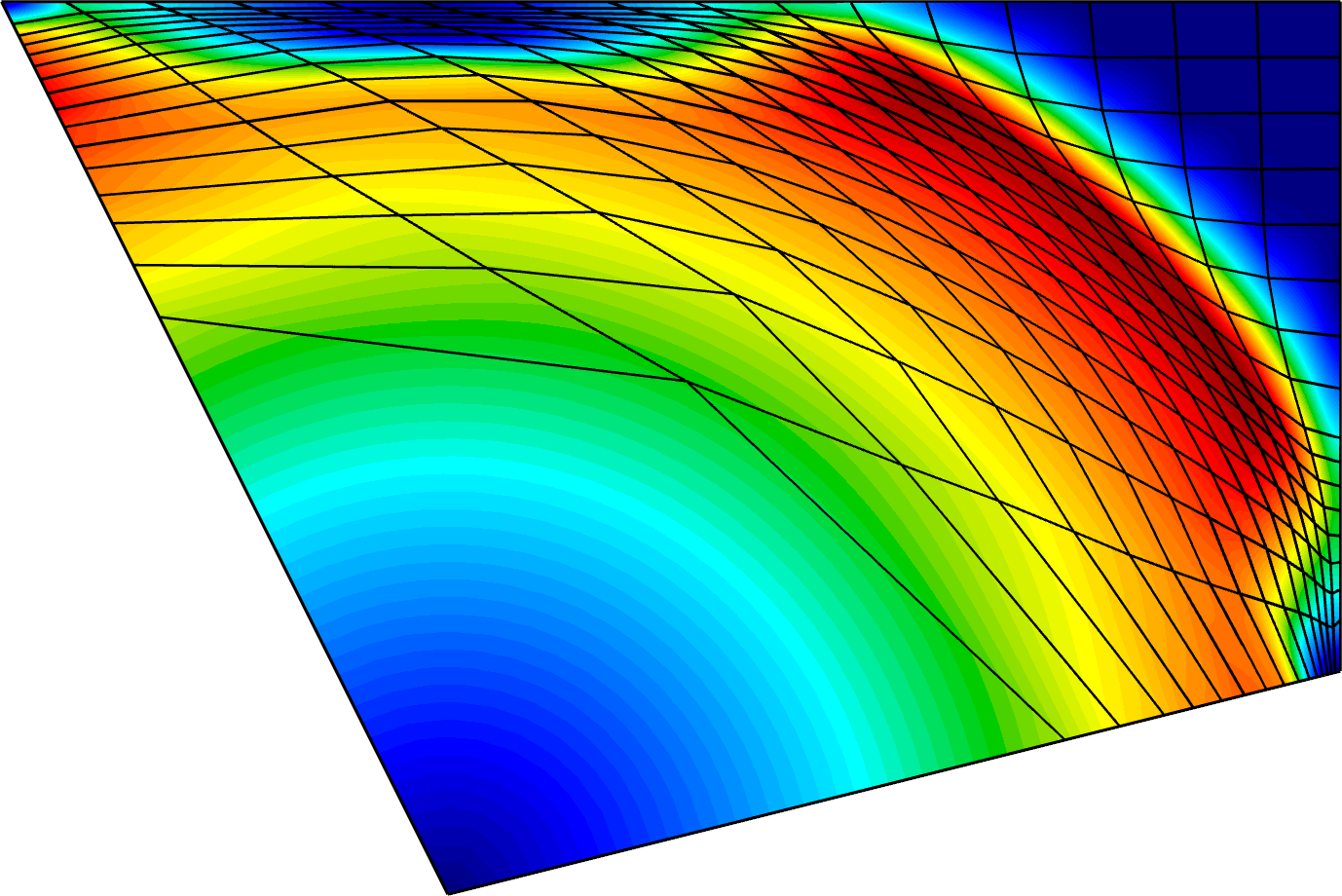} 
		\qquad & \qquad  
		\includegraphics[width=0.30\linewidth]{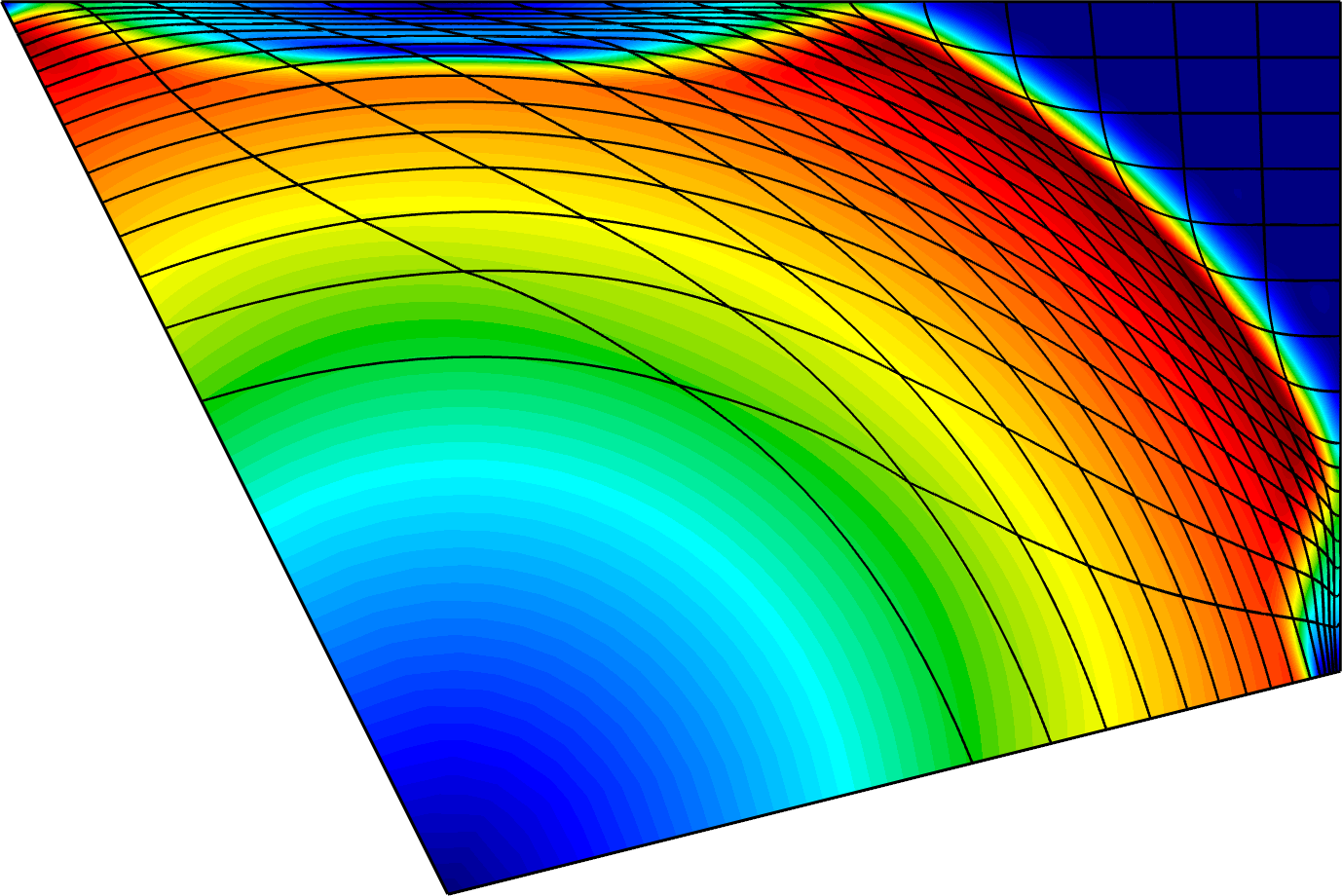}
		\qquad & \qquad  
		\includegraphics[width=0.30\linewidth]{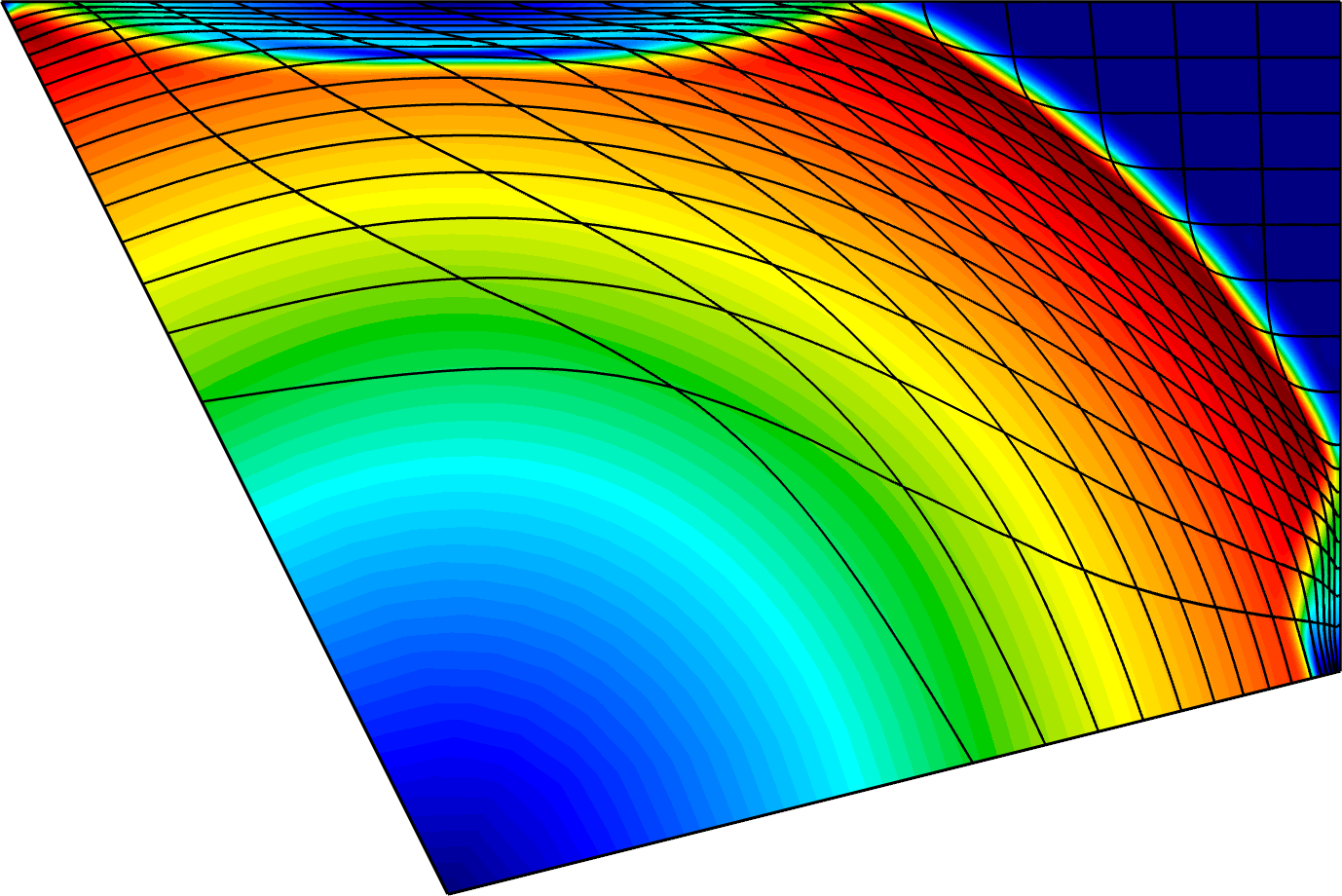} \\[.4cm]
		\multicolumn{3}{c}{Density}  \\
		\includegraphics[width=0.30\linewidth]{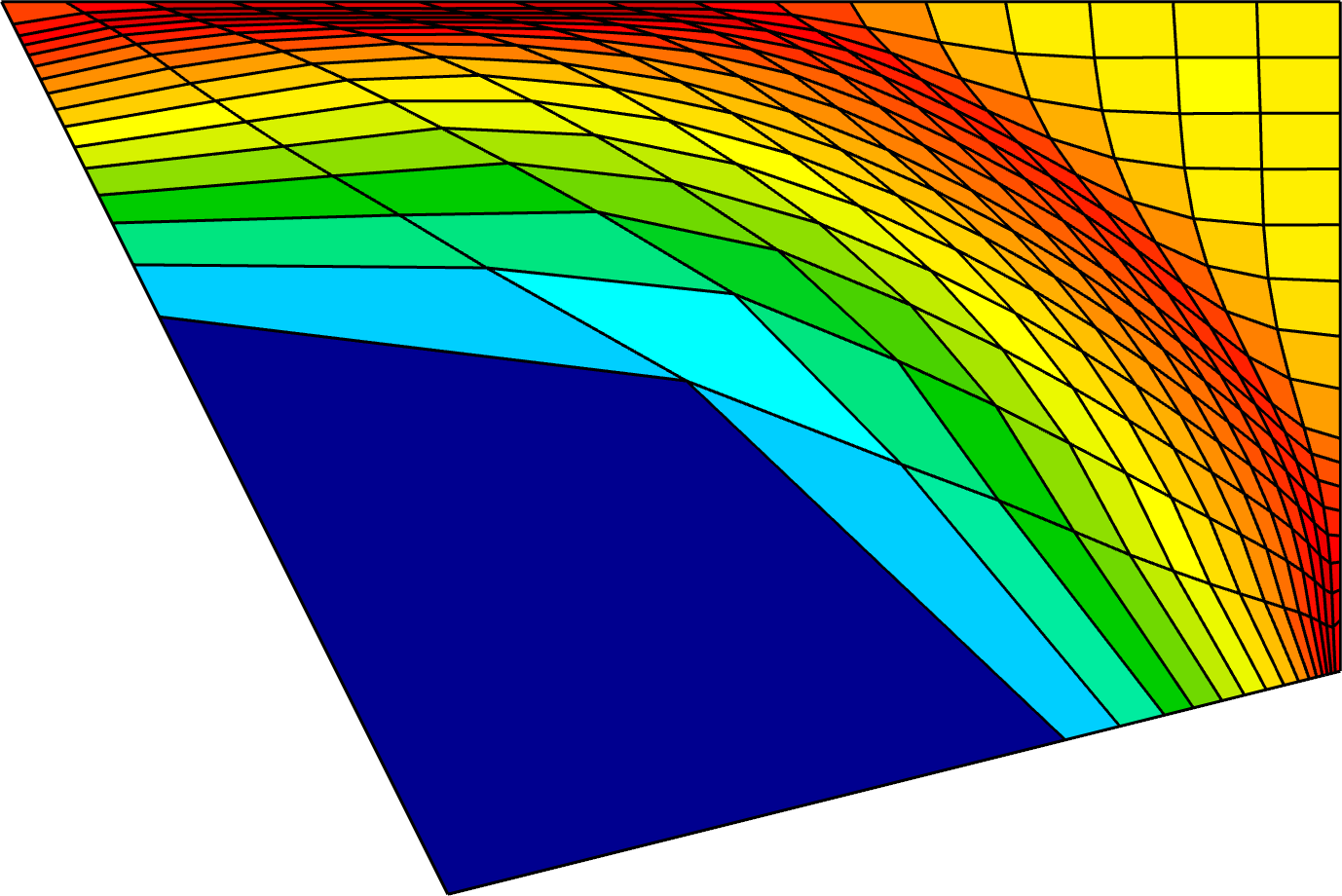} 
		\qquad & \qquad  
		\includegraphics[width=0.30\linewidth]{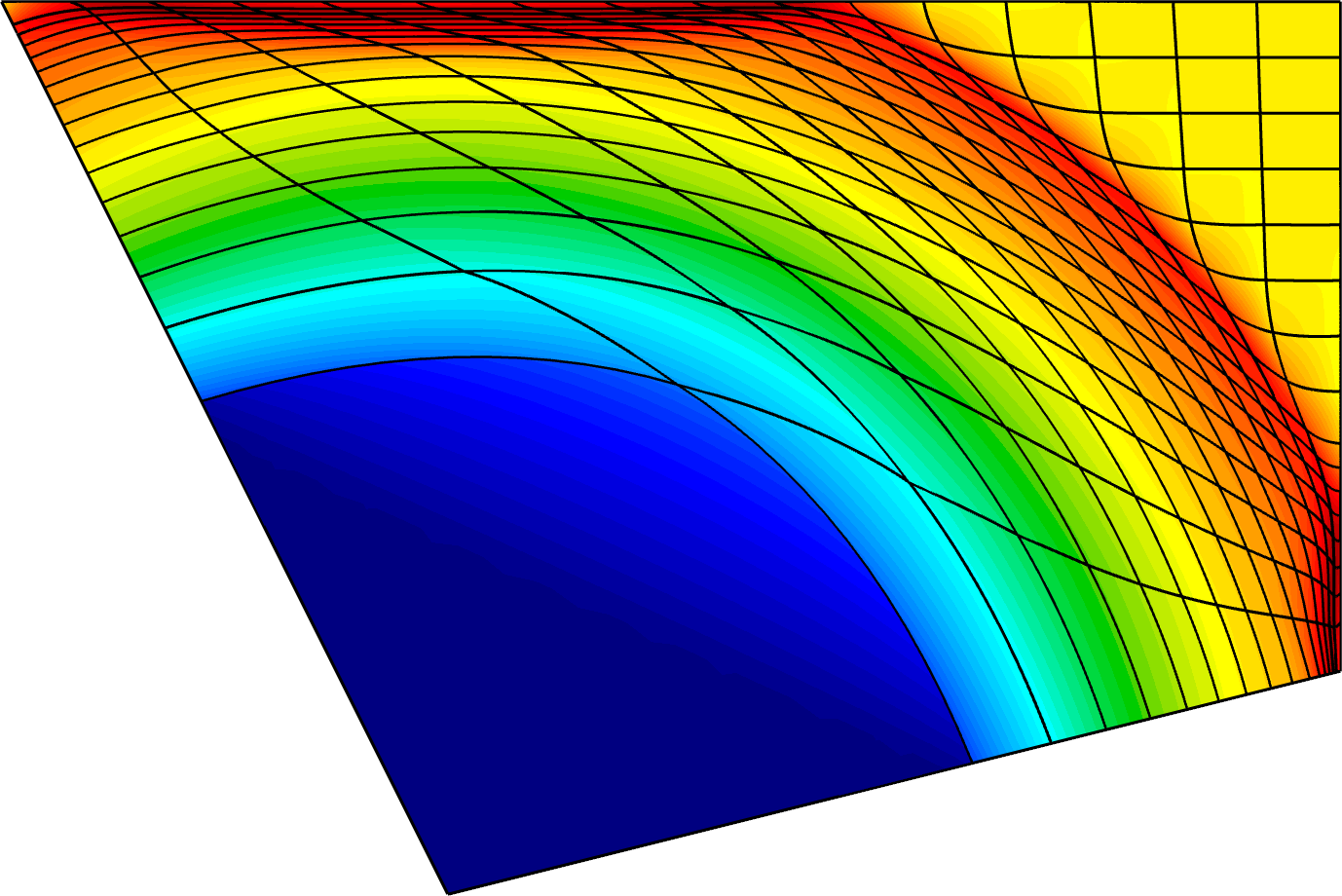}
		\qquad & \qquad  
		\includegraphics[width=0.30\linewidth]{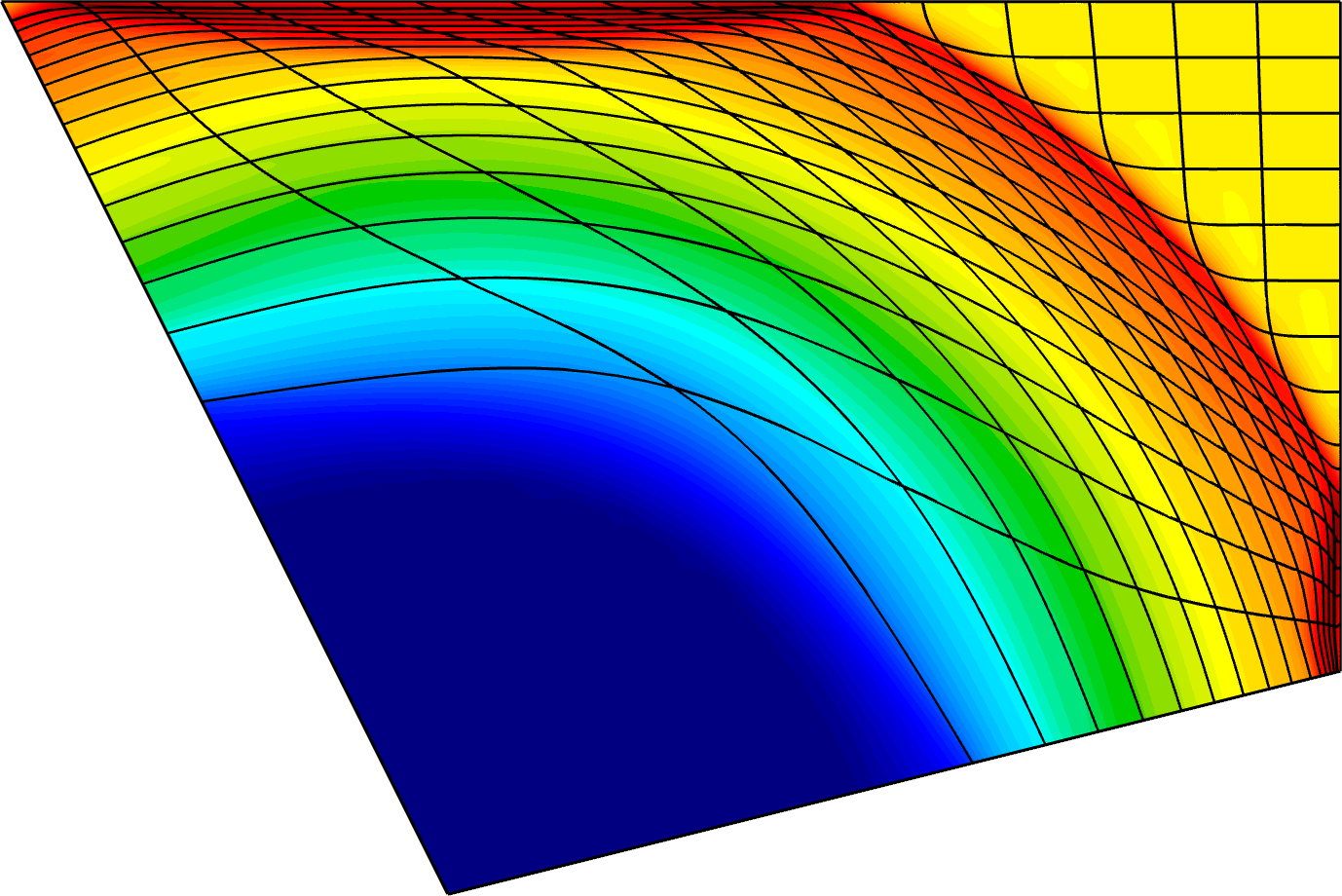}  \\[.4cm]
		\multicolumn{3}{c}{Mesh Deformation}  \\
		\includegraphics[width=0.30\linewidth]{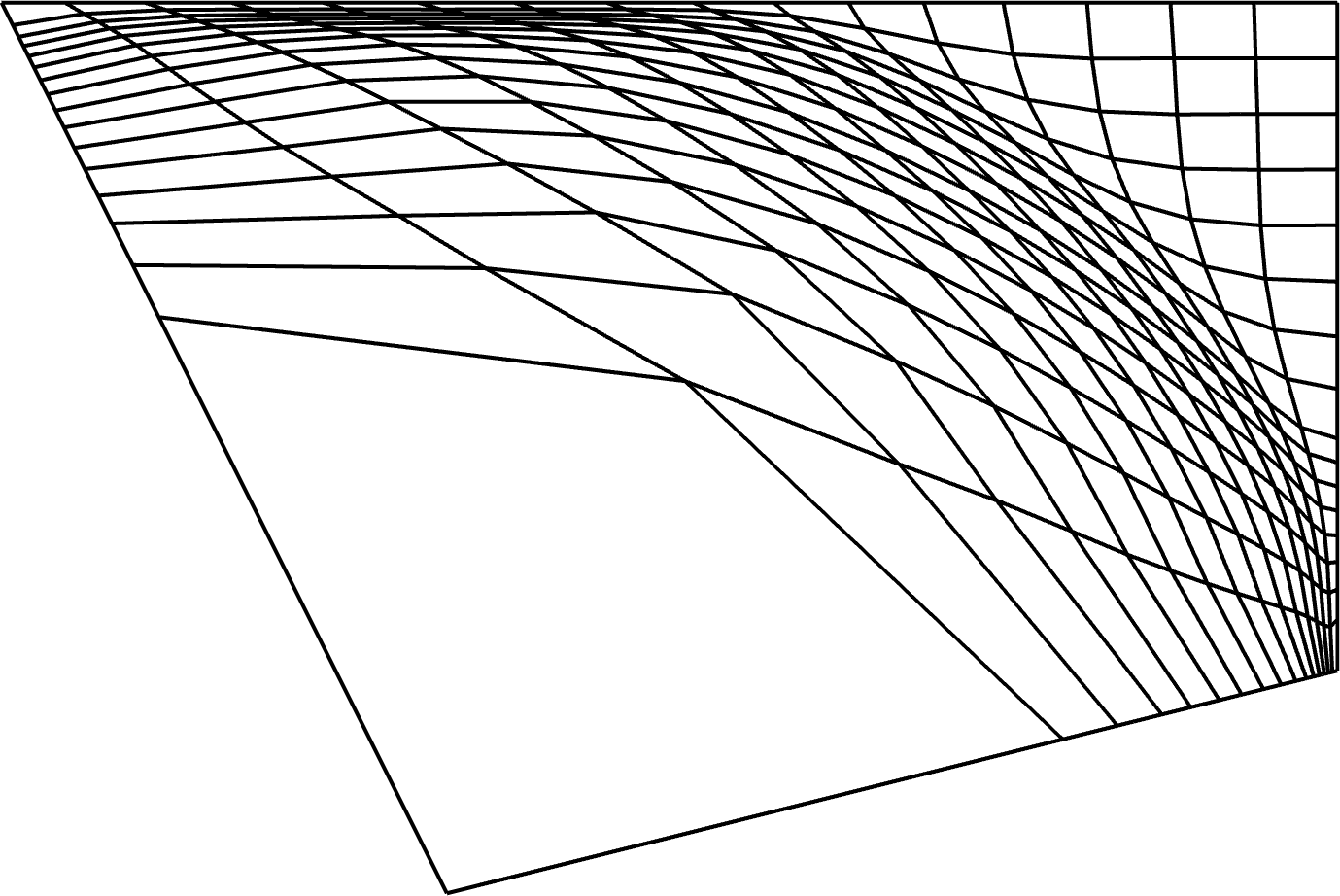} 
		\qquad & \qquad  
		\includegraphics[width=0.30\linewidth]{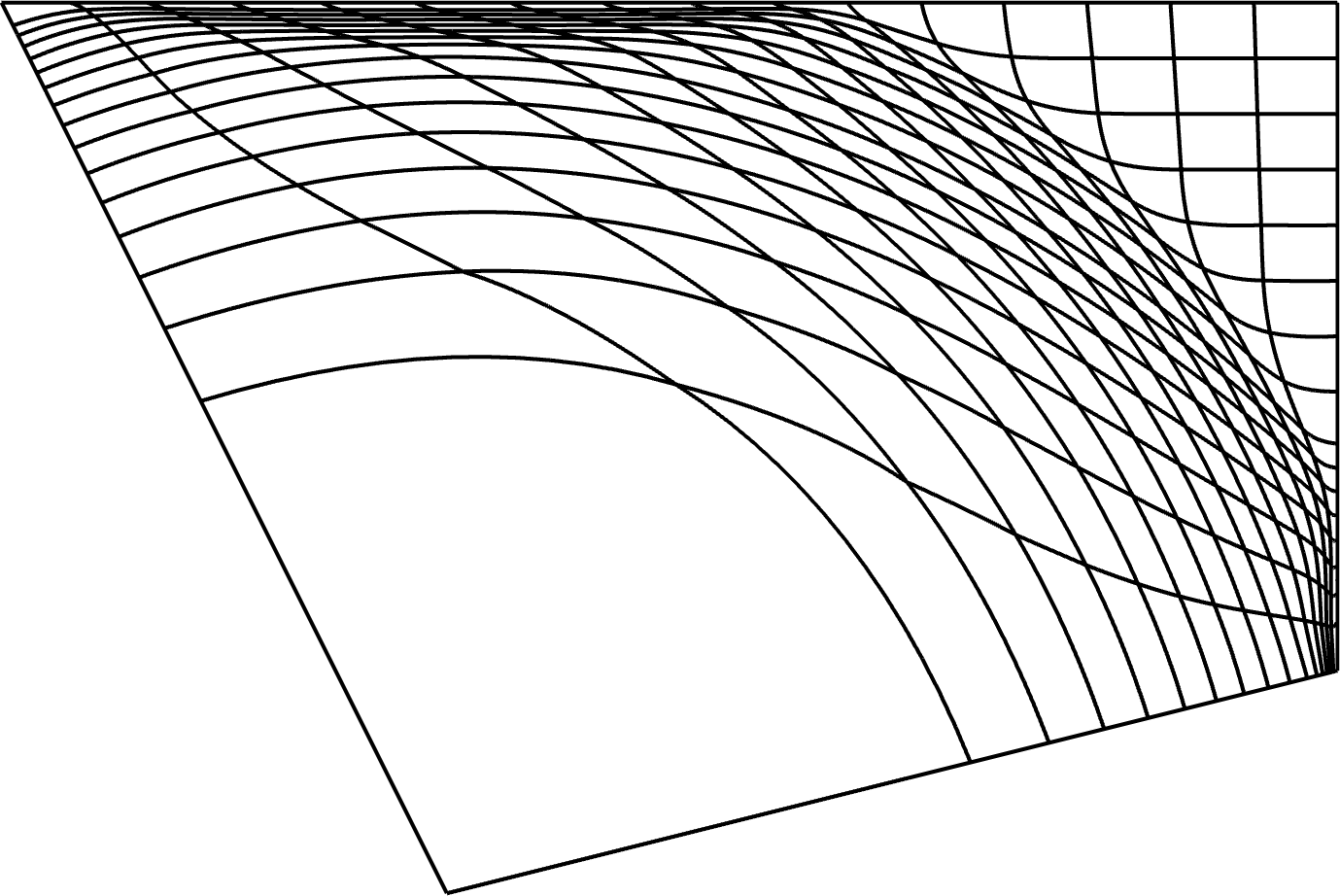}
		\qquad & \qquad  
		\includegraphics[width=0.30\linewidth]{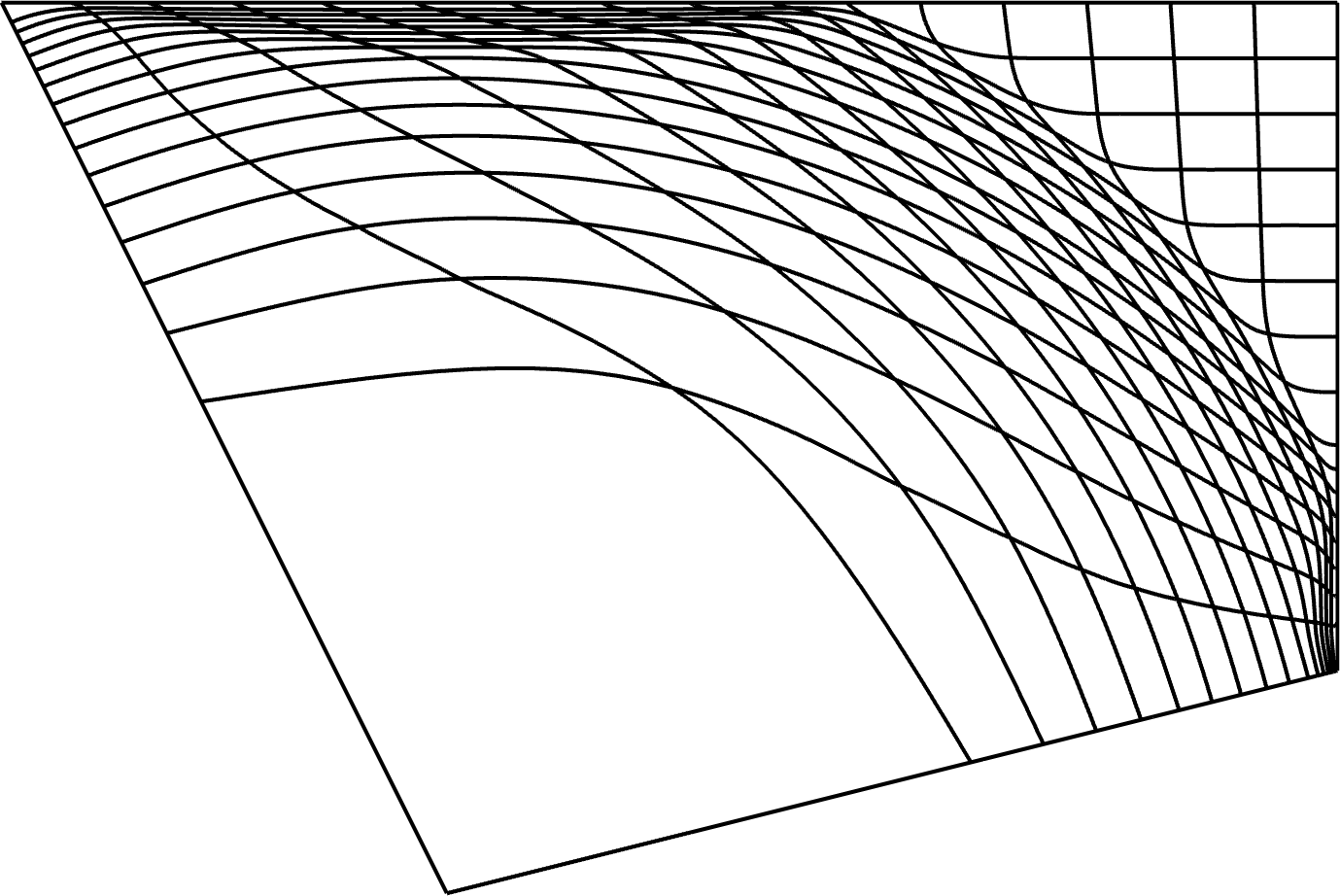}  \\[.0cm]
		$Q_{1}-Q_{0}$ 
		&
		$Q_{2}-Q_{1}$
		&
		$Q_{3}-Q_{2}$
	\end{tabular}
	\caption{Plots of the velocity and density fields in addition to the mesh deformation for the Sedov test in a trapezoidal domain using $Q_{1}-Q_{0}$, $Q_{2}-Q_{1}$, $Q_{3}-Q_{2}$ velocity-energy pairs.}
	\label{fig:TrapezoidSedovresults}
\end{figure}

We perform the Sedov test in a trapezoidal domain and show in Figure~\ref{fig:TrapezoidSedovresults} plots of the velocity and density fields in addition to the mesh deformation at the final time of $t = 1.3$ for the $Q_{1}-Q_{0}$, $Q_{2}-Q_{1}$, $Q_{3}-Q_{2}$ velocity-energy pairs. 
It is worth mentioning that conducting this test with strong wall boundary enforcement is cumbersome, as it would require enforcement of linear constraints on the boundary DOFs.
The weak wall boundary conditions produce the correct shock bounce-back behavior on the top boundary with solutions that remain smooth and do not show any unphysical oscillations.

\subsection{Two-dimensional Sedov explosion in a square with a circular hole}
\label{square_hole_sedov_2d}
\begin{figure}[tb]
	\centering
	\begin{tabular}{ccc}
		\multicolumn{3}{c}{Velocity}  \\
		\includegraphics[width=0.30\linewidth]{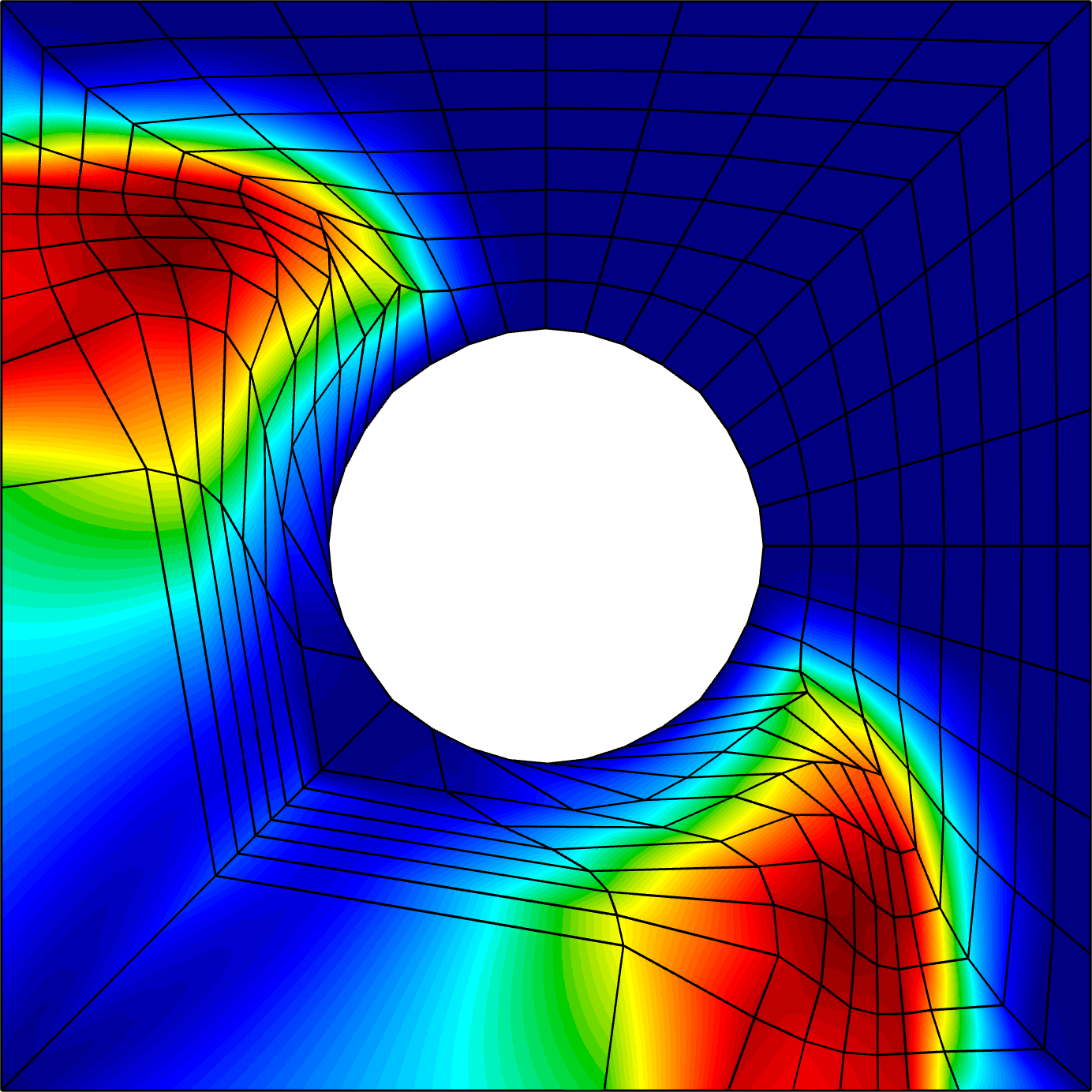} 
		\qquad & \qquad  
		\includegraphics[width=0.30\linewidth]{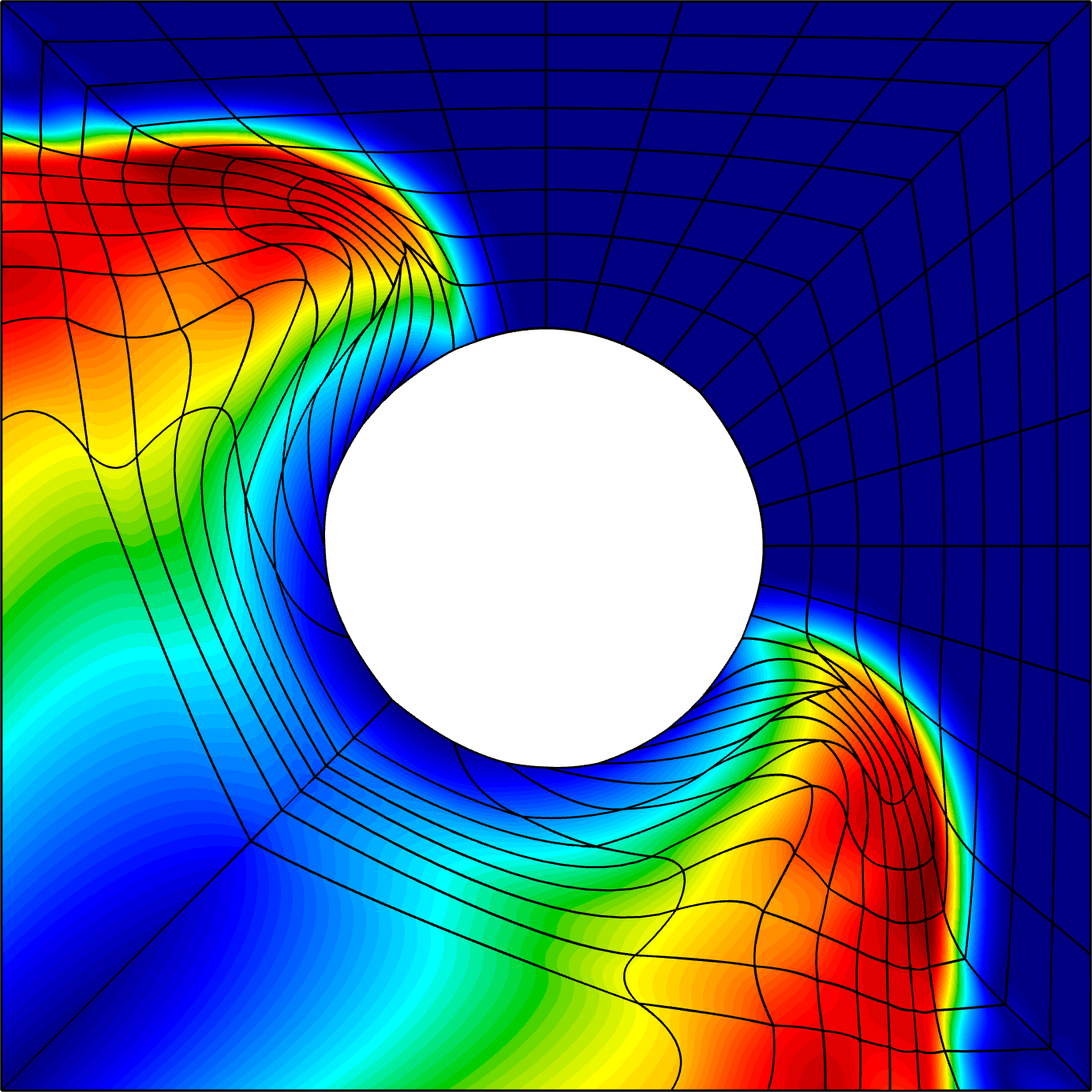}
		\qquad & \qquad  
		\includegraphics[width=0.30\linewidth]{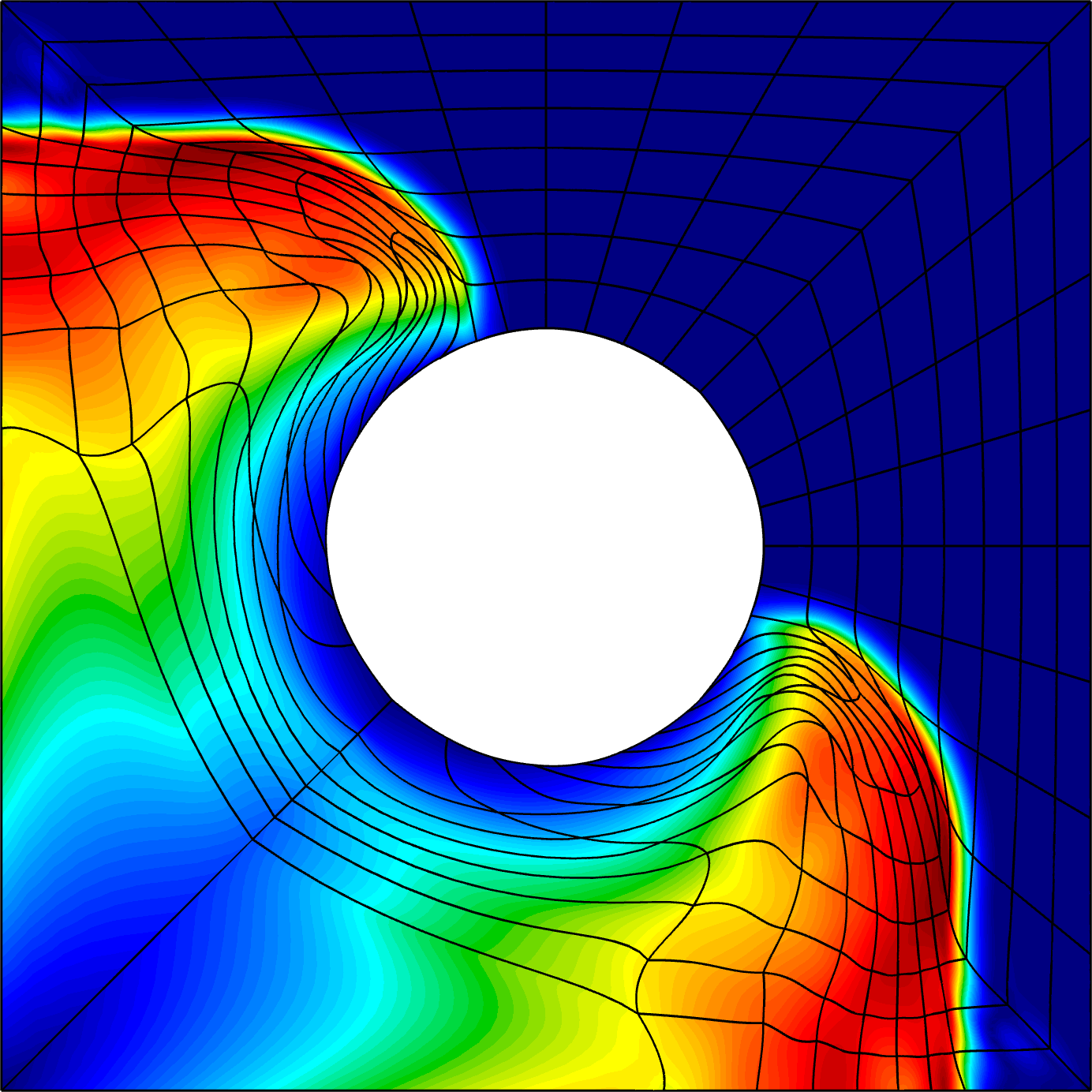} \\[.4cm]
		\multicolumn{3}{c}{Density}  \\
		\includegraphics[width=0.30\linewidth]{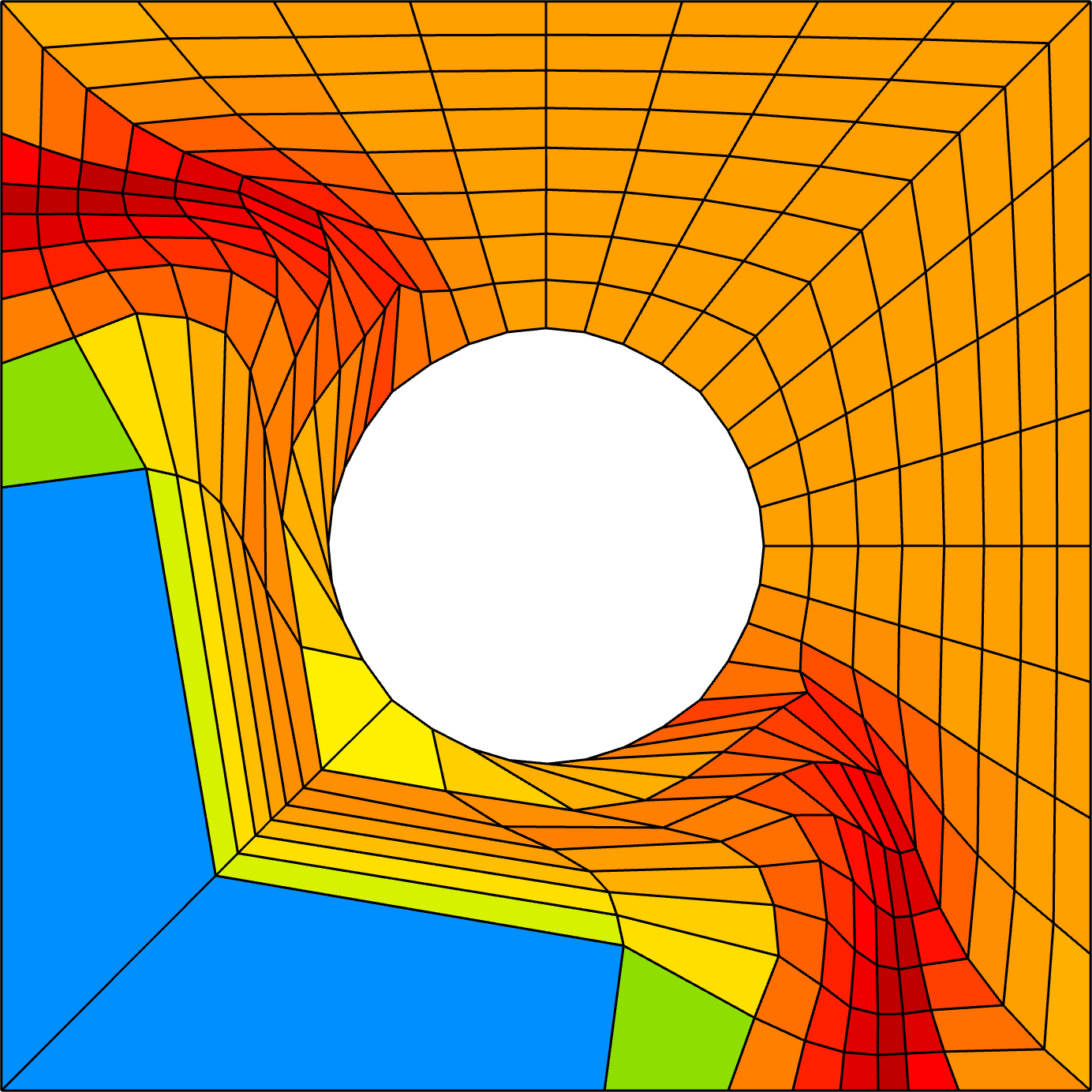} 
		\qquad & \qquad  
		\includegraphics[width=0.30\linewidth]{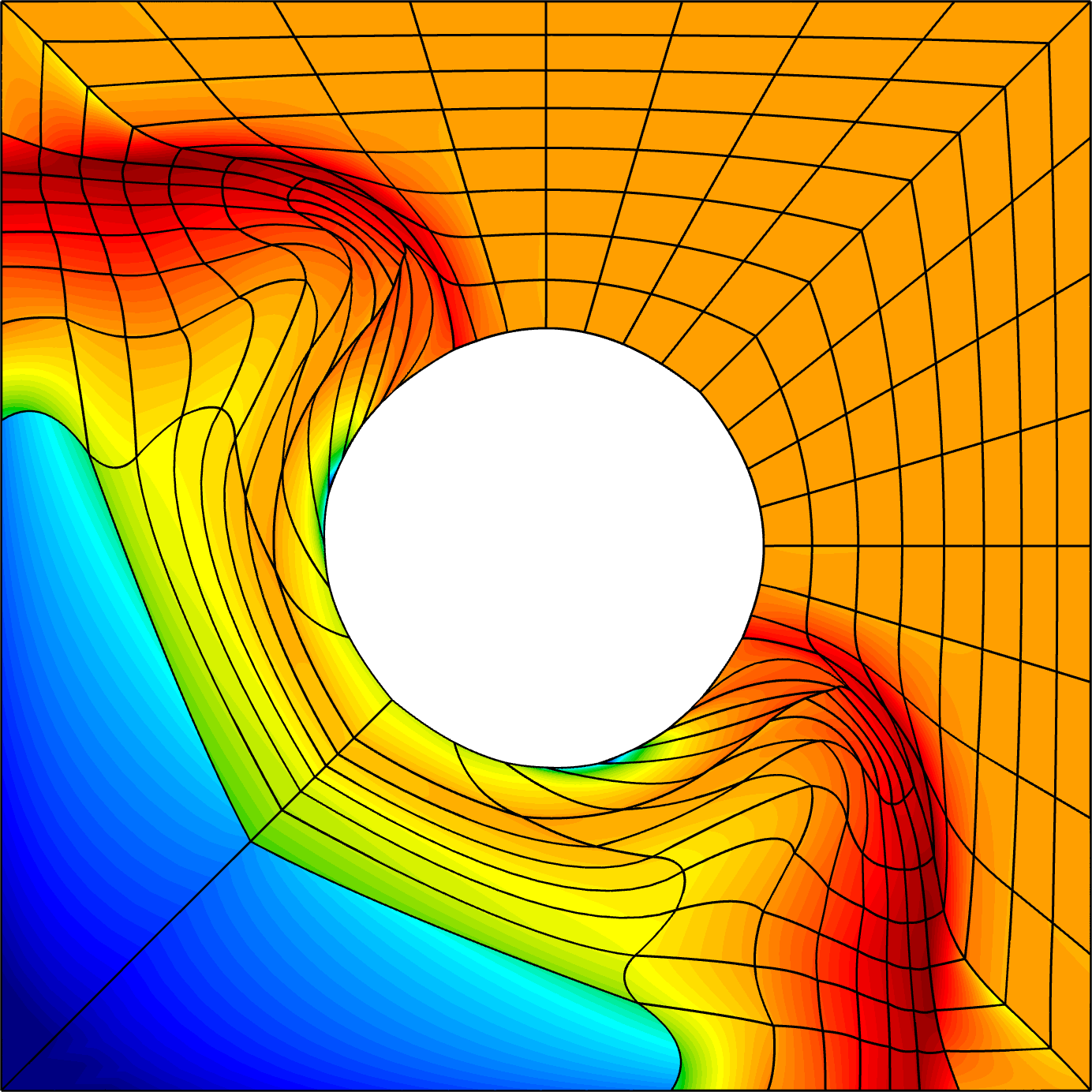}
		\qquad & \qquad  
		\includegraphics[width=0.30\linewidth]{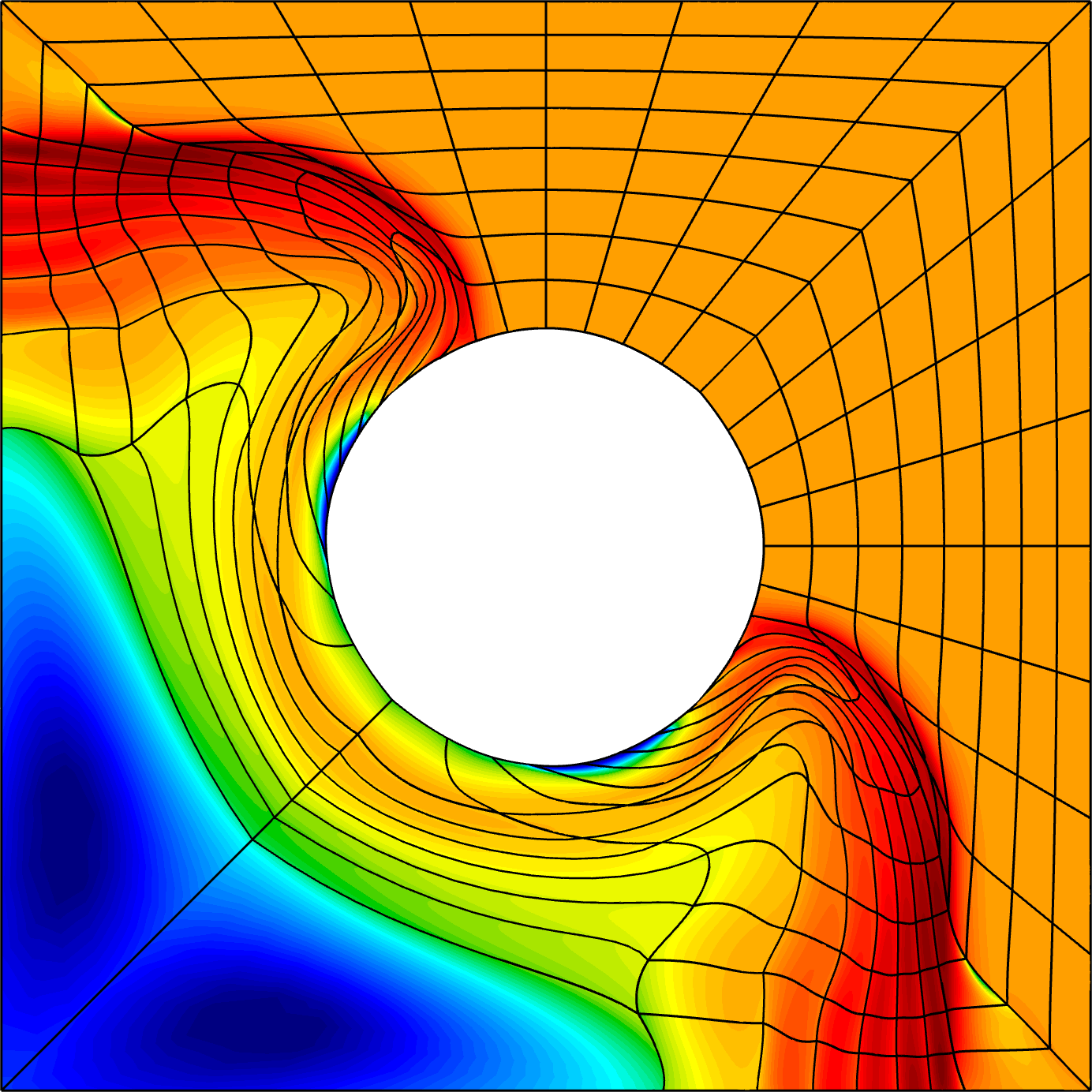}  \\[.4cm]
		\multicolumn{3}{c}{Mesh Deformation}  \\
		\includegraphics[width=0.30\linewidth]{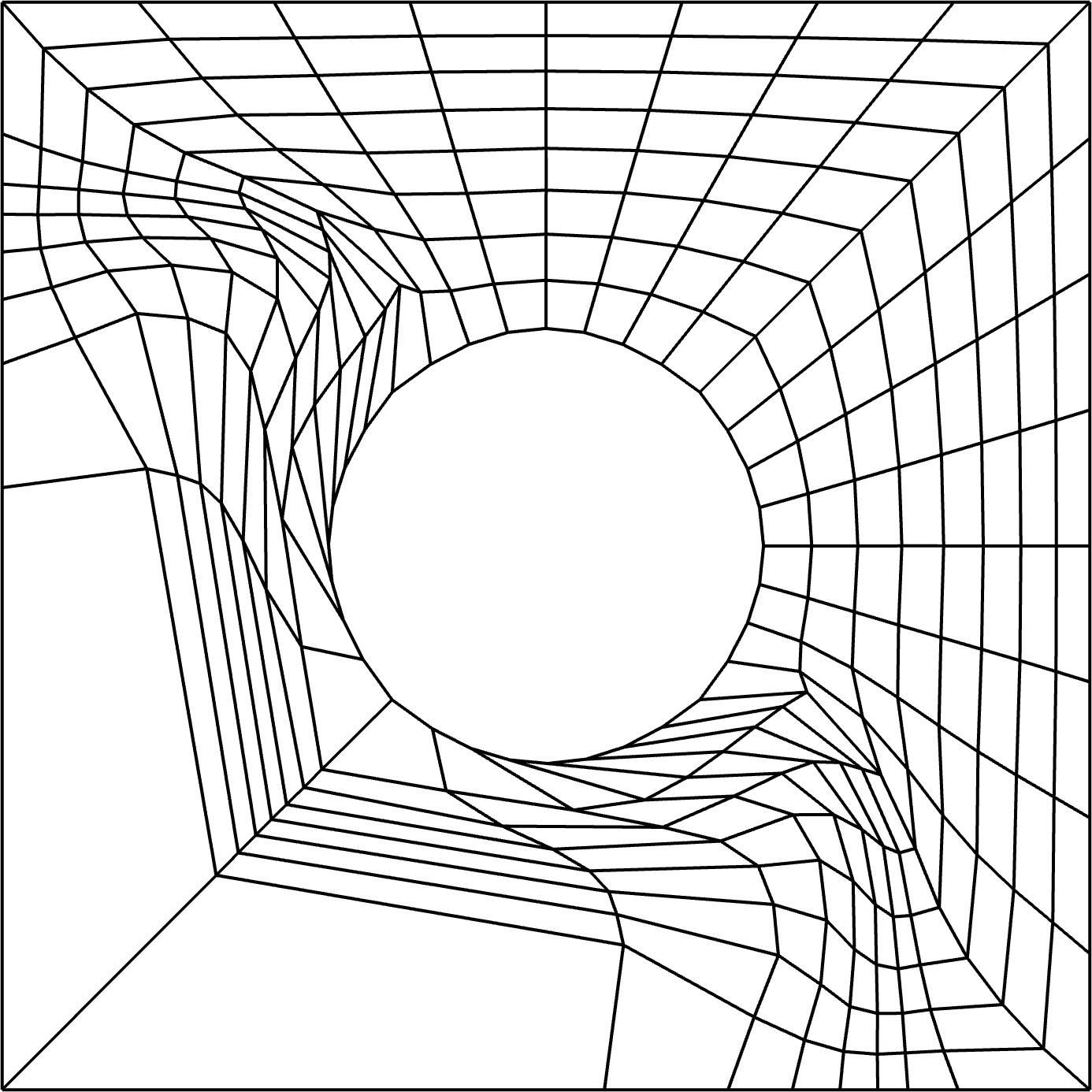} 
		\qquad & \qquad  
		\includegraphics[width=0.30\linewidth]{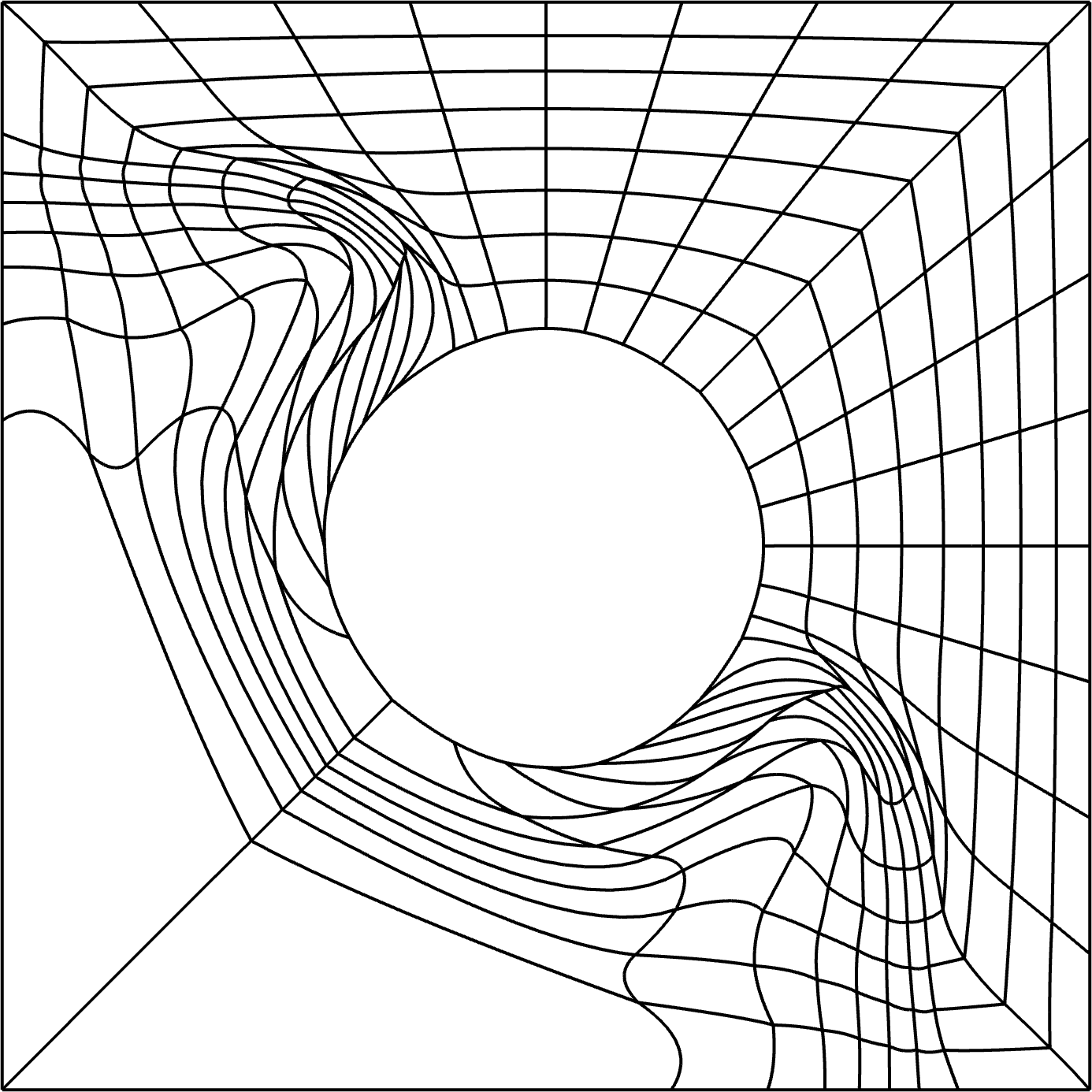}
		\qquad & \qquad  
		\includegraphics[width=0.30\linewidth]{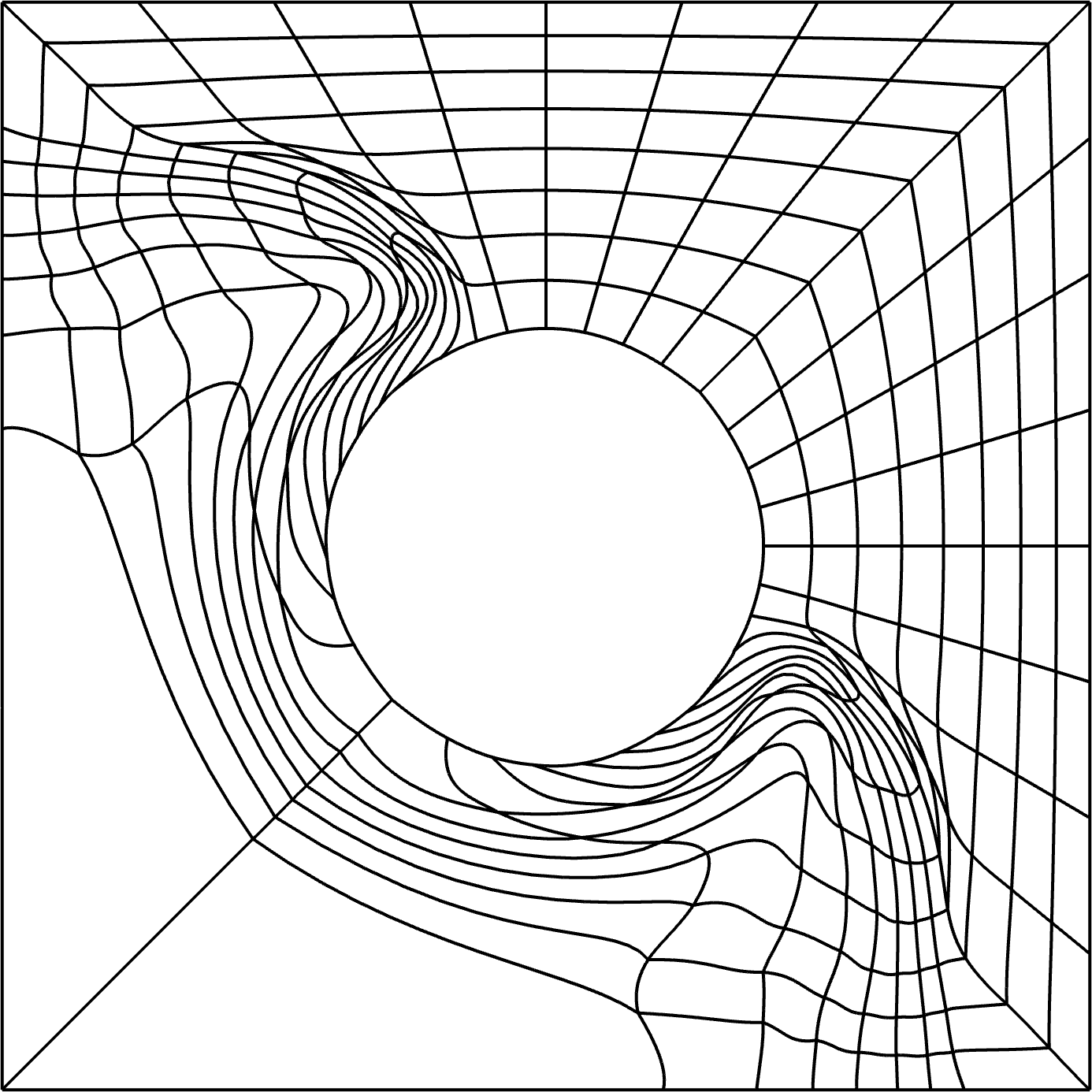}  \\[.0cm]
		$Q_{1}-Q_{0}$ 
		&
		$Q_{2}-Q_{1}$
		&
		$Q_{3}-Q_{2}$
	\end{tabular}
	\caption{Plots of the velocity and density fields in addition to the mesh deformation for the Sedov test in a square domain with a circular hole using the $Q_{1}-Q_{0}$, $Q_{2}-Q_{1}$, $Q_{3}-Q_{2}$ velocity-energy pairs.}
	\label{fig:SquareHoleSedovresults}
\end{figure}

We perform the Sedov test in a unit square with a circular hole and show in Figure~\ref{fig:SquareHoleSedovresults} plots of the velocity and density  fields in addition to the mesh deformation at the final time of $t = 0.8$ for the $Q_{1}-Q_{0}$, $Q_{2}-Q_{1}$, $Q_{3}-Q_{2}$ velocity-energy pairs.
It is worth mentioning that conducting this test with strong wall boundary enforcement is very cumbersome, due to the curvature of the internal surface, while it is seamless with our proposed weak form.
The weak wall boundary enforcement produces the expected response and
smooth solutions without any unphysical oscillations.
This test demonstrates the robustness of the proposed algorithm, which allows the shock wave to propagate around the curved obstacle under extreme levels of mesh deformation, even for the $Q_{3}-Q_{2}$ discretization.

\subsection{Two-dimensional Sedov explosion in a square with a square hole}
\label{square_square_hole_sedov_2d}
\begin{figure}[tb]
	\centering
	\begin{tabular}{ccc}
		\multicolumn{3}{c}{Velocity}  \\
		\includegraphics[width=0.30\linewidth]{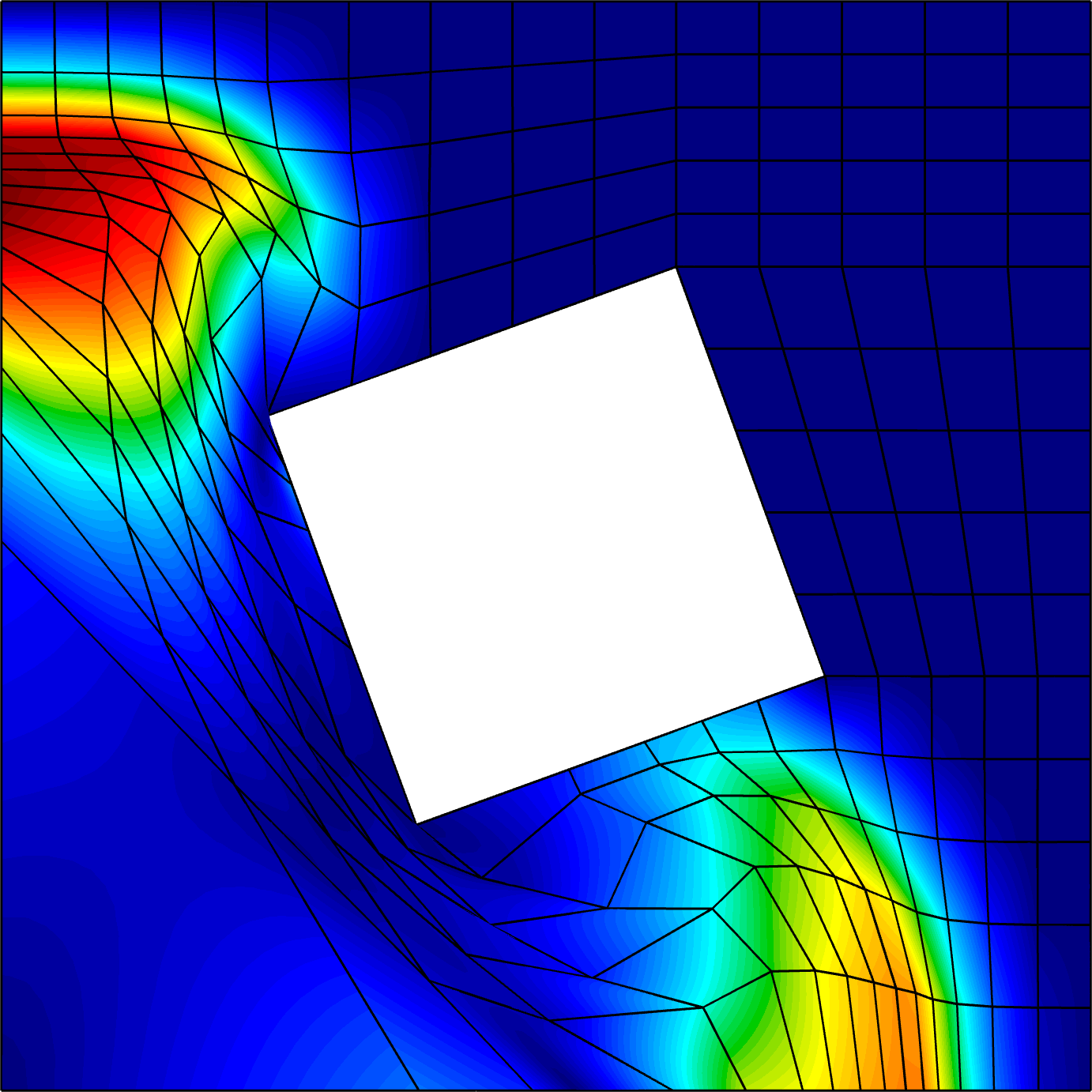} 
		\qquad & \qquad  
		\includegraphics[width=0.30\linewidth]{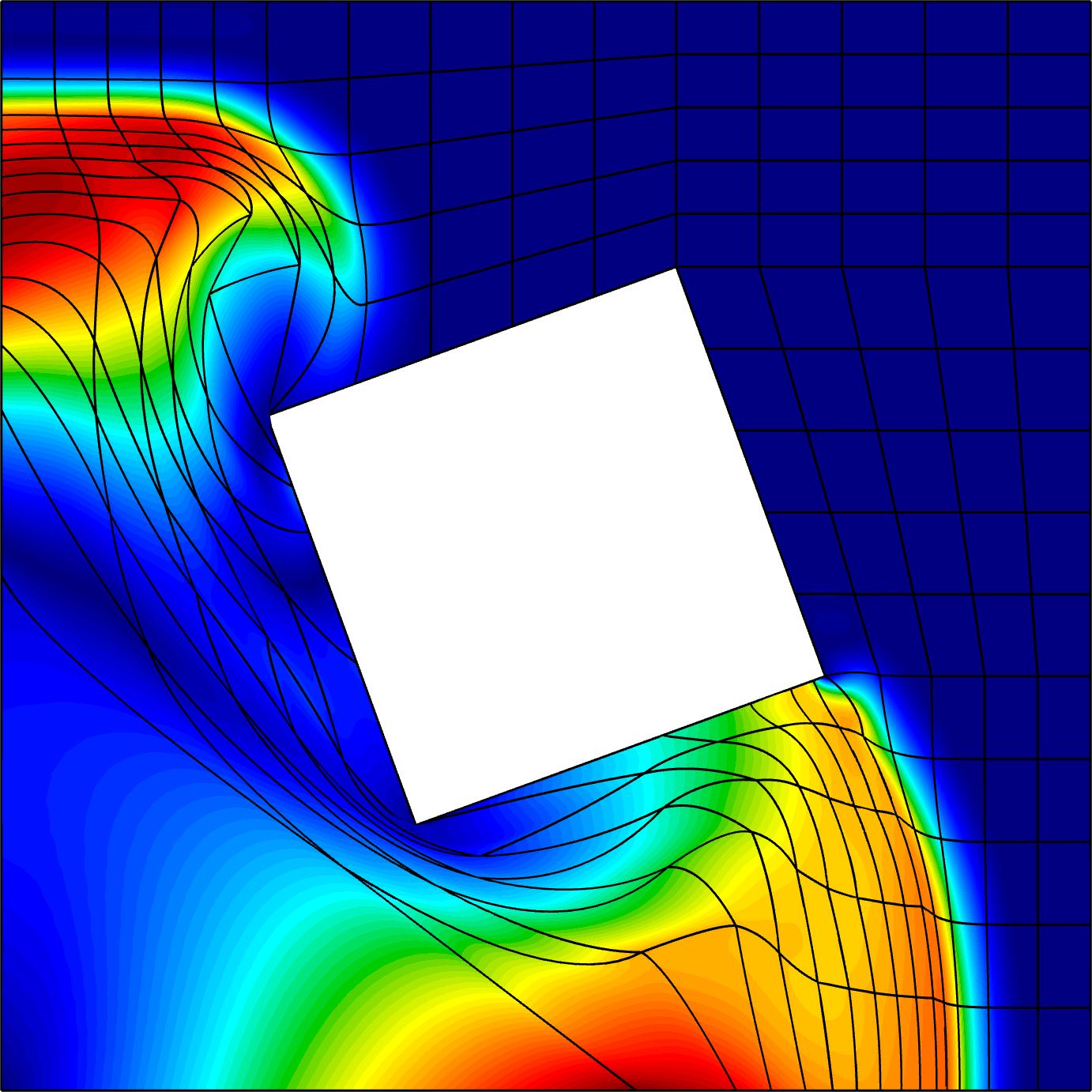}
		\qquad & \qquad  
		\includegraphics[width=0.30\linewidth]{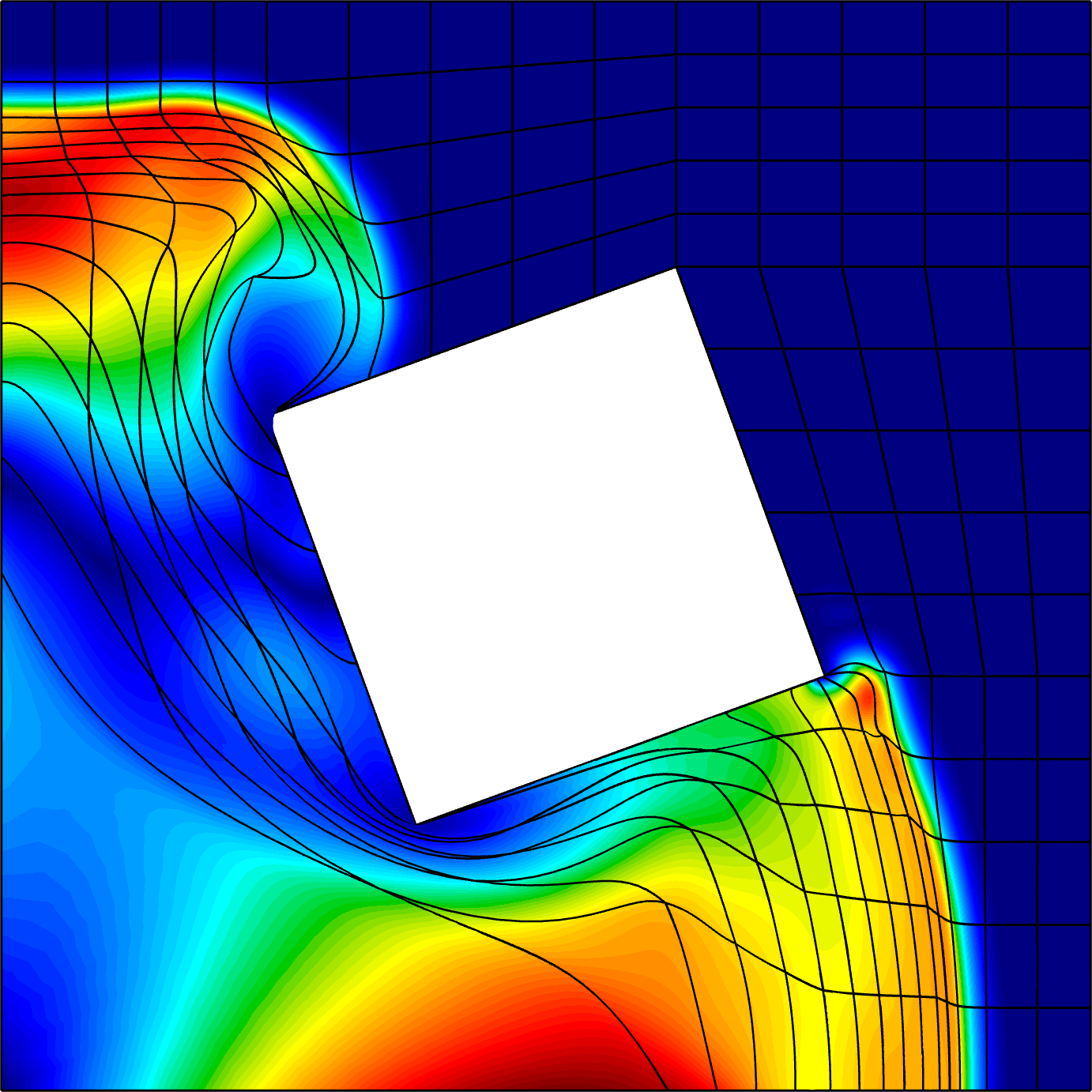} \\[.4cm]
		\multicolumn{3}{c}{Density}  \\
		\includegraphics[width=0.30\linewidth]{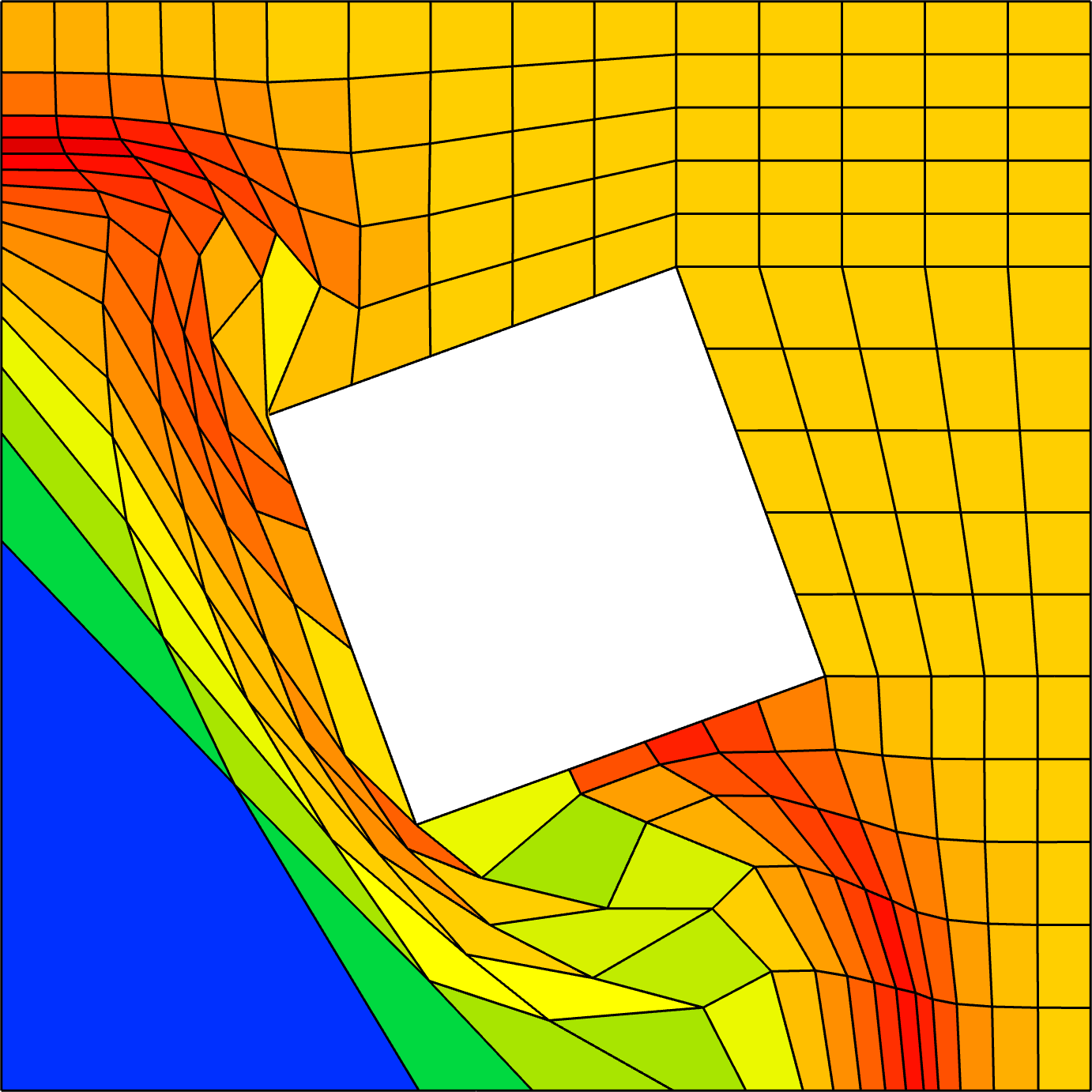} 
		\qquad & \qquad  
		\includegraphics[width=0.30\linewidth]{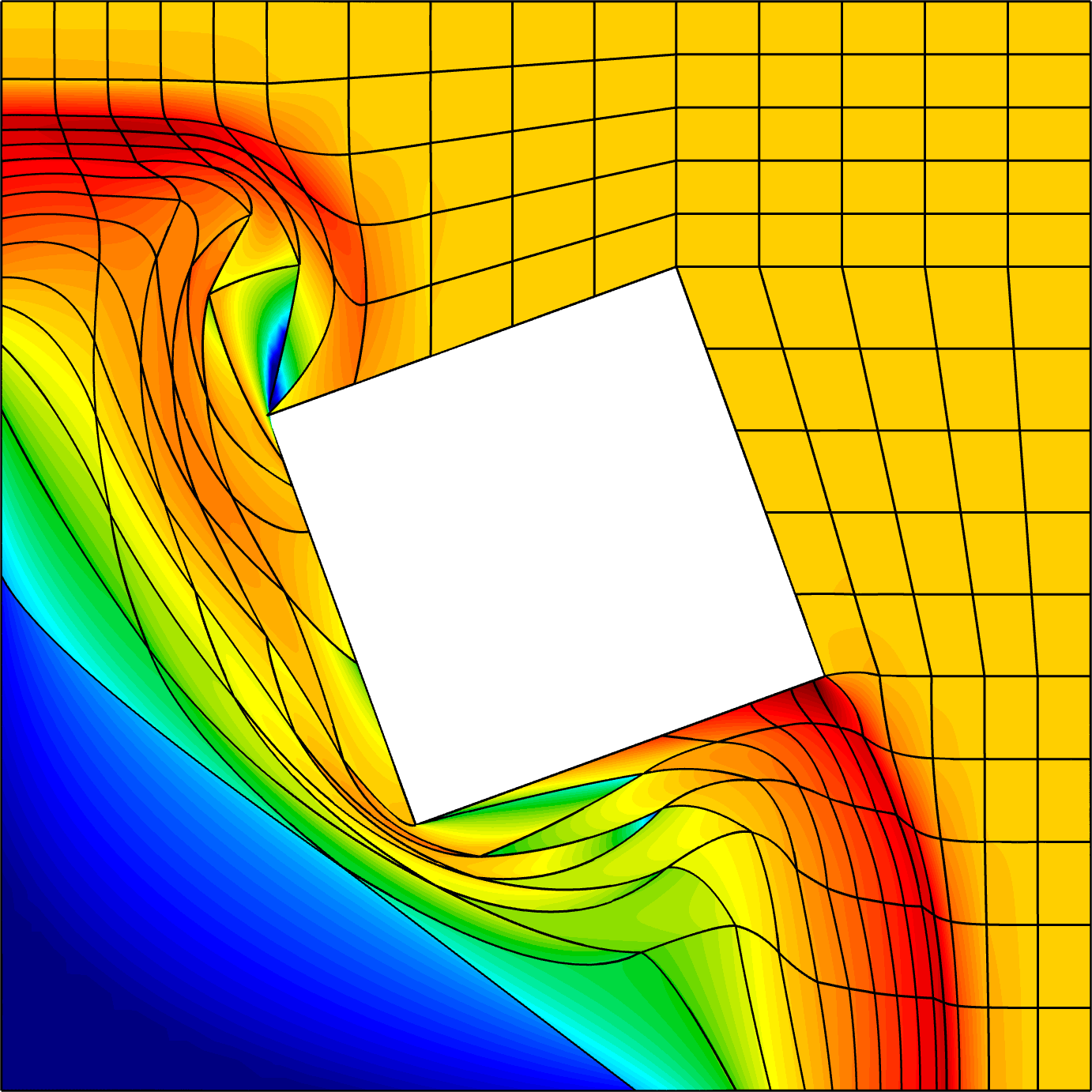}
		\qquad & \qquad  
		\includegraphics[width=0.30\linewidth]{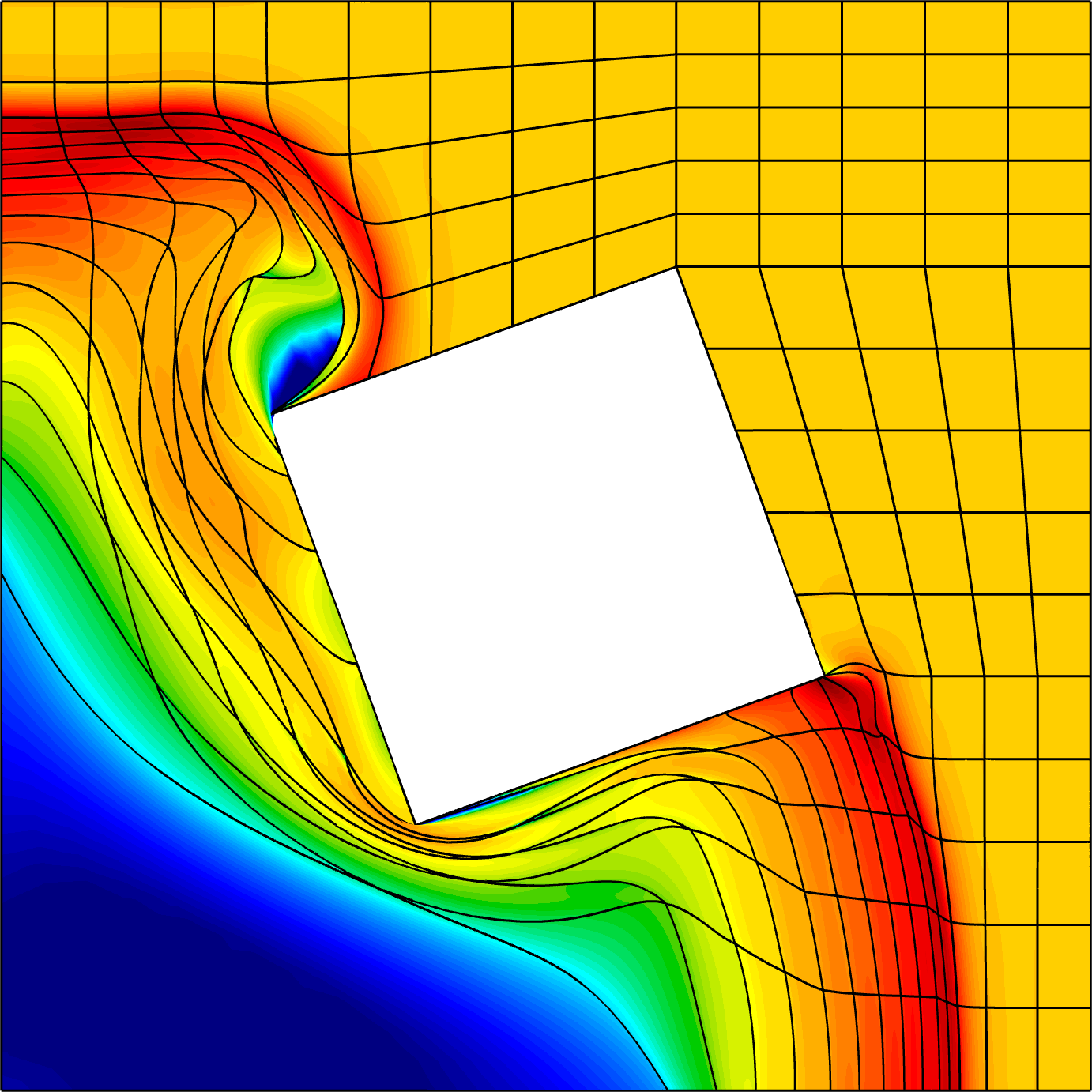}  \\[.4cm]
		\multicolumn{3}{c}{Mesh Deformation}  \\
		\includegraphics[width=0.30\linewidth]{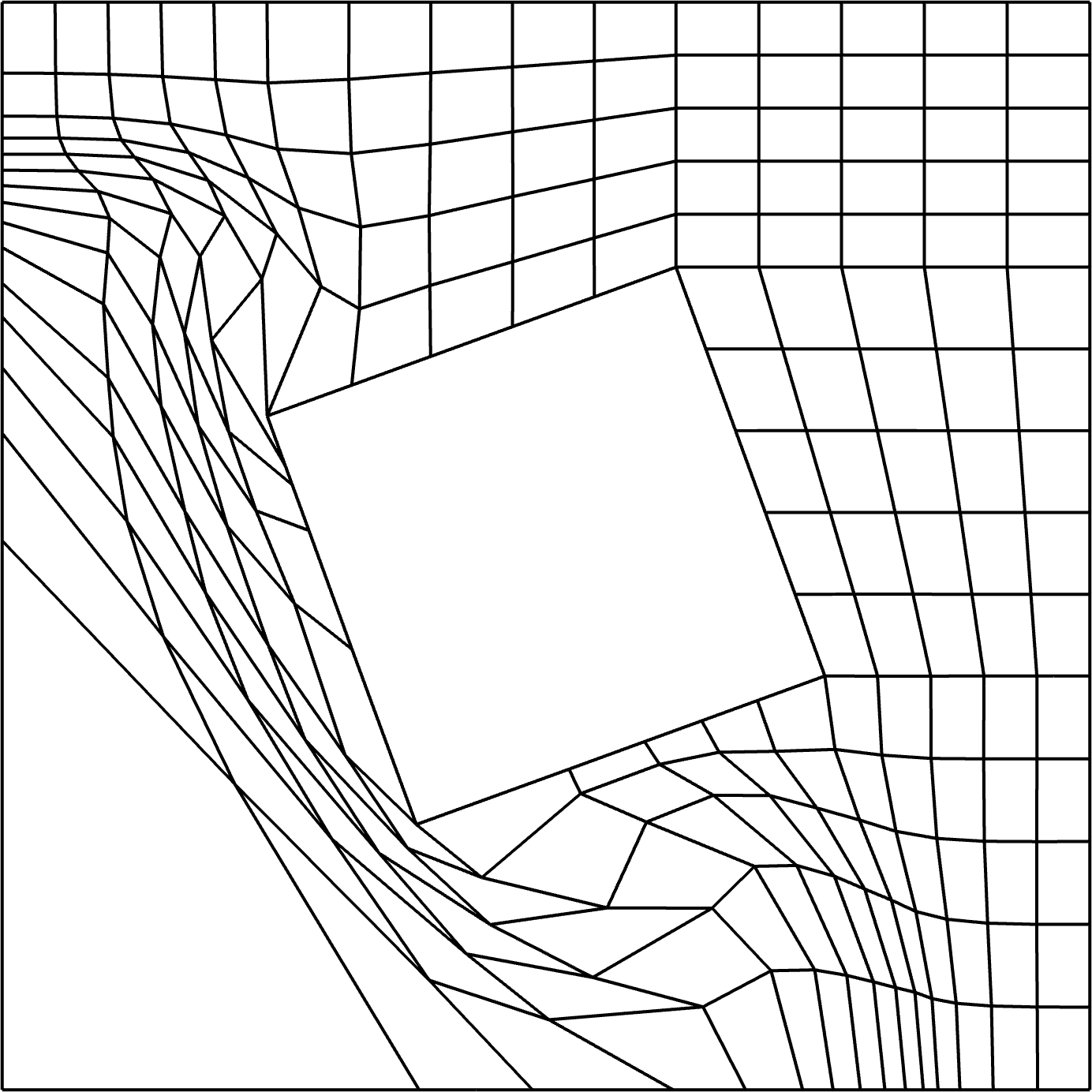} 
		\qquad & \qquad  
		\includegraphics[width=0.30\linewidth]{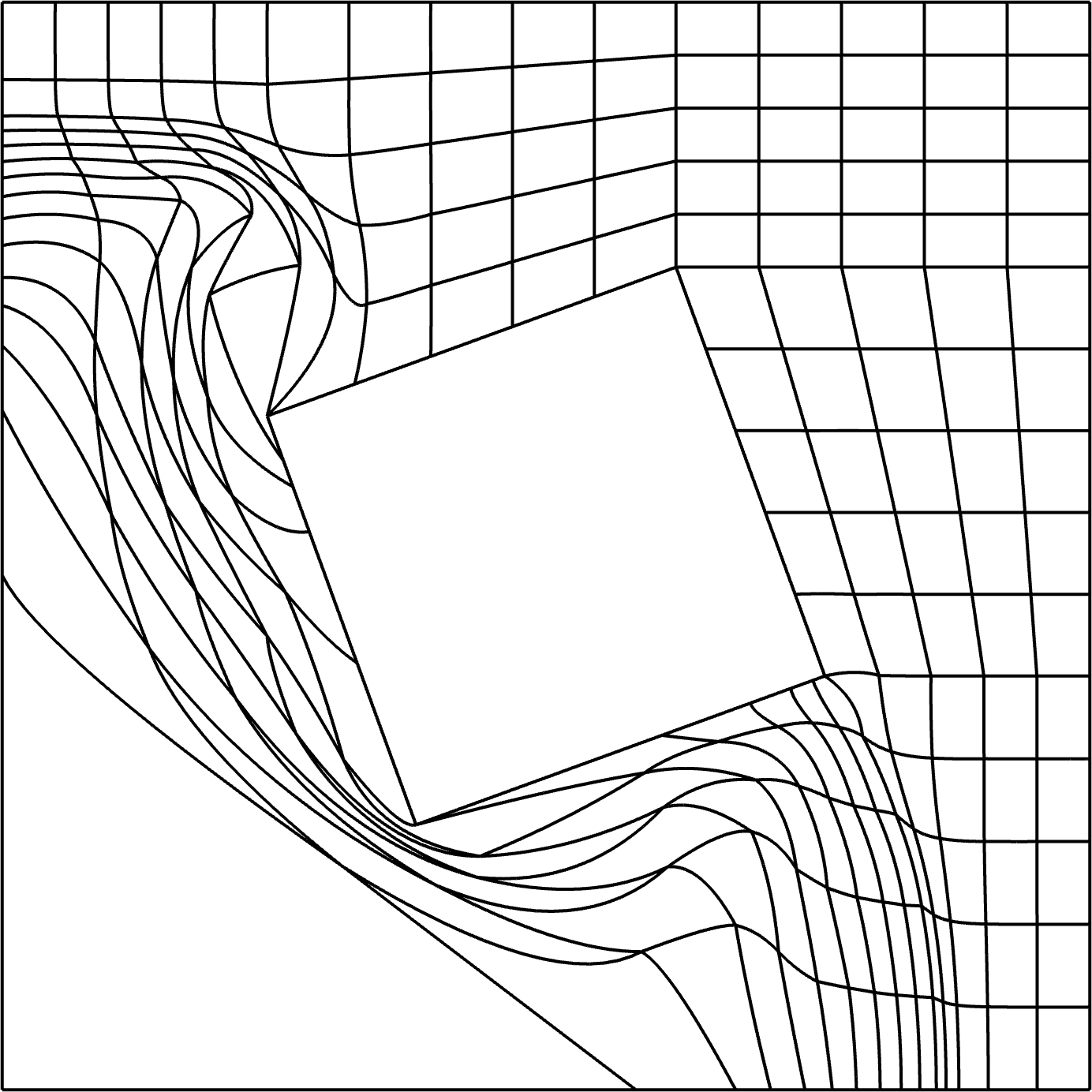}
		\qquad & \qquad  
		\includegraphics[width=0.30\linewidth]{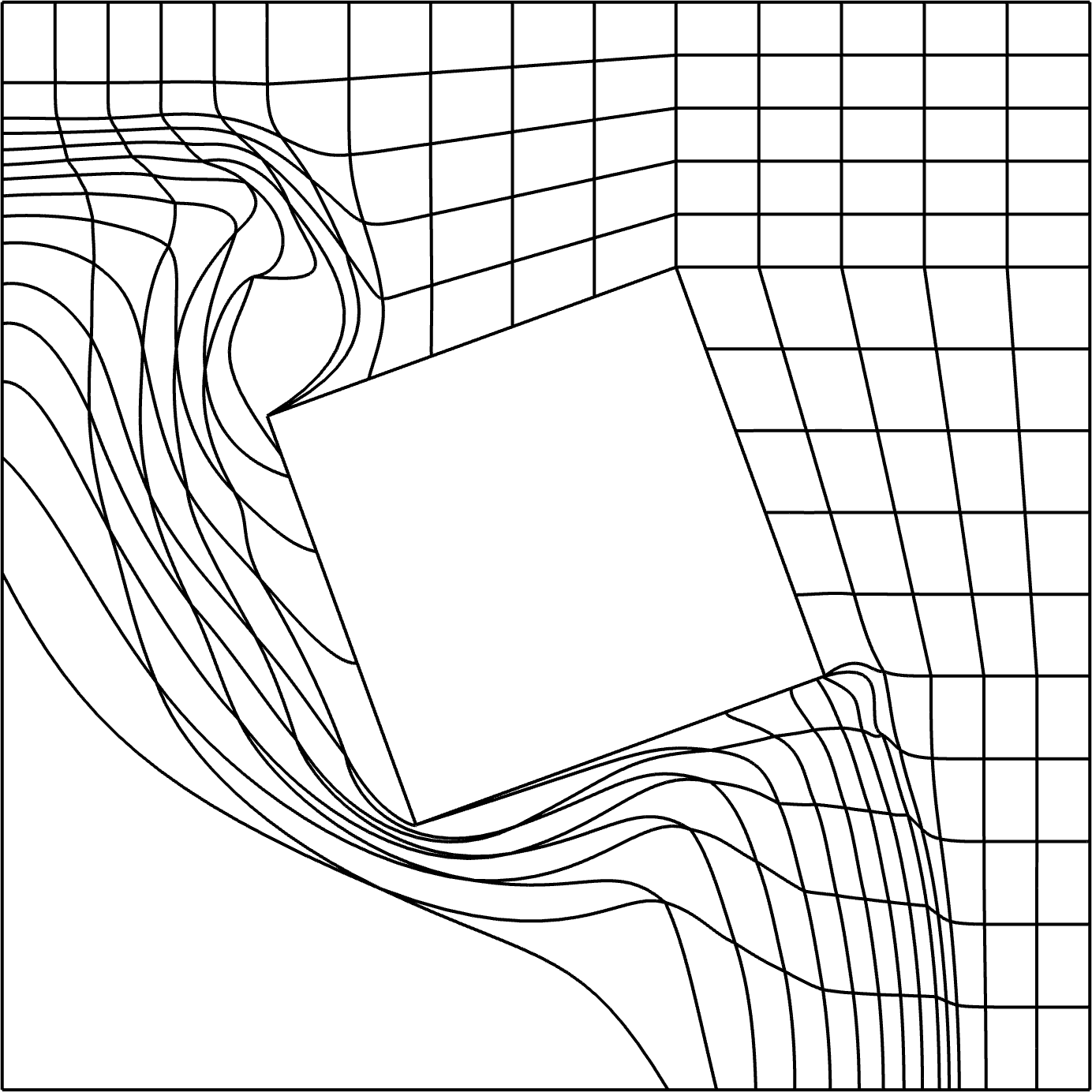}  \\[.0cm]
		$Q_{1}-Q_{0}$ 
		&
		$Q_{2}-Q_{1}$
		&
		$Q_{3}-Q_{2}$
	\end{tabular}
	\caption{Plots of the velocity and density fields in addition to the mesh deformation for the Sedov test in a square domain with a square hole using the $Q_{1}-Q_{0}$, $Q_{2}-Q_{1}$, $Q_{3}-Q_{2}$ velocity-energy pairs.}
	\label{fig:SquareSquareHoleSedovresults}
\end{figure}

We perform the Sedov test in a unit square with a square hole rotated by $20$ degrees counter-clockwise and show plots of the velocity and density  fields in addition to the mesh deformation at the final time of $t = 0.8$ for the $Q_{1}-Q_{0}$, $Q_{2}-Q_{1}$, $Q_{3}-Q_{2}$ velocity-energy pairs as shown in Figure~\ref{fig:SquareSquareHoleSedovresults}.
The weak wall boundary conditions produce the expected response and smooth solutions without any unphysical oscillations.
The proposed method demonstrates its robustness by capturing the propagation of the shock wave past the sharp top left corner of the obstacle.

\subsection{Two-dimensional Sedov explosion in a disc}
\label{circle_sedov_2d}
\begin{figure}[tb]
	\centering
	\begin{tabular}{cccc}
		\text{}
		\qquad & \qquad
		\text{Velocity} 
		\qquad & \qquad
		\text{Density} 
		\qquad & \qquad
		\text{Mesh Deformation} 
		\\[.4cm]
		$t = 0.01$
		\qquad & \qquad
		\raisebox{-.5\height}{\includegraphics[width=0.205\linewidth]{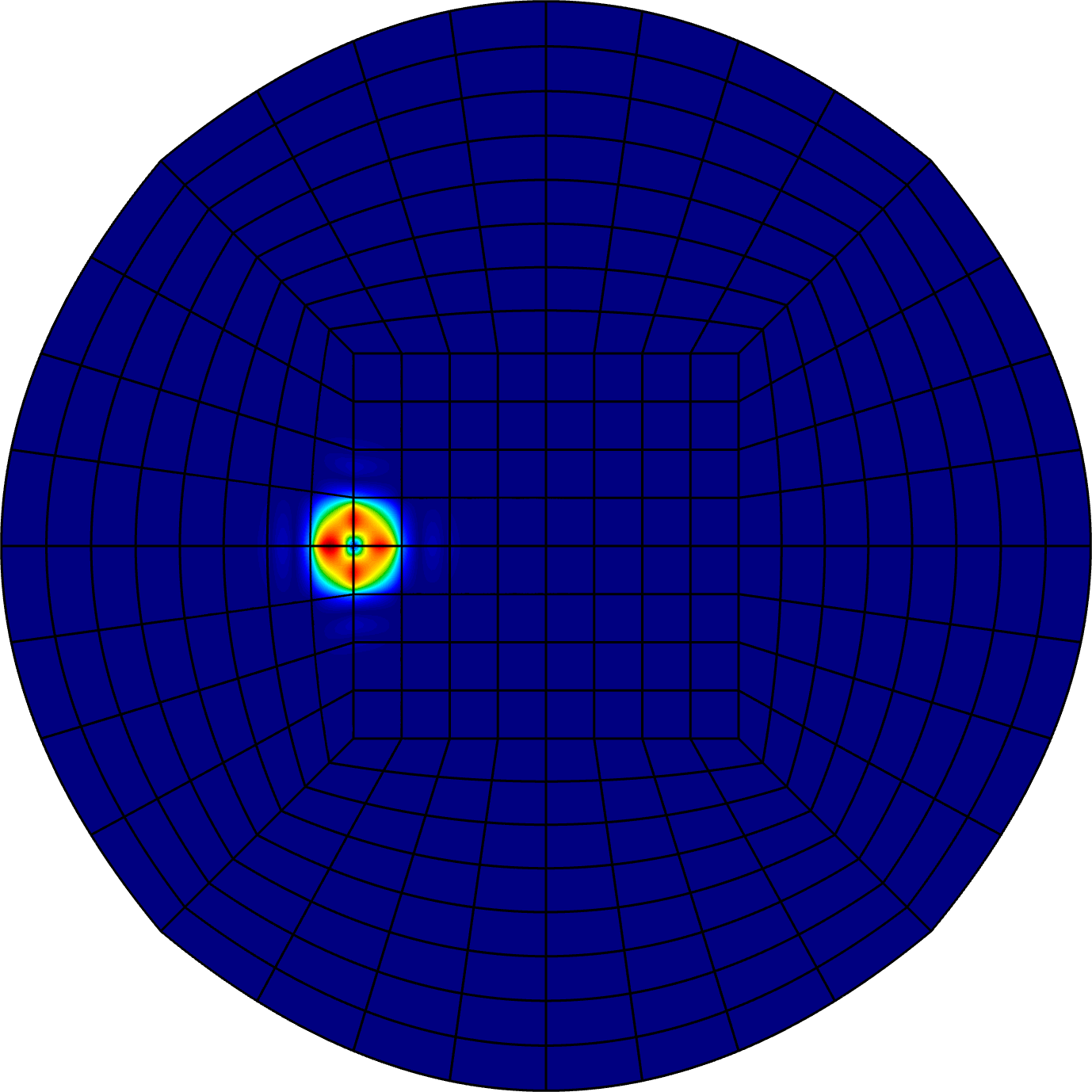}}
		\qquad & \qquad
		\raisebox{-.5\height}{\includegraphics[width=0.205\linewidth]{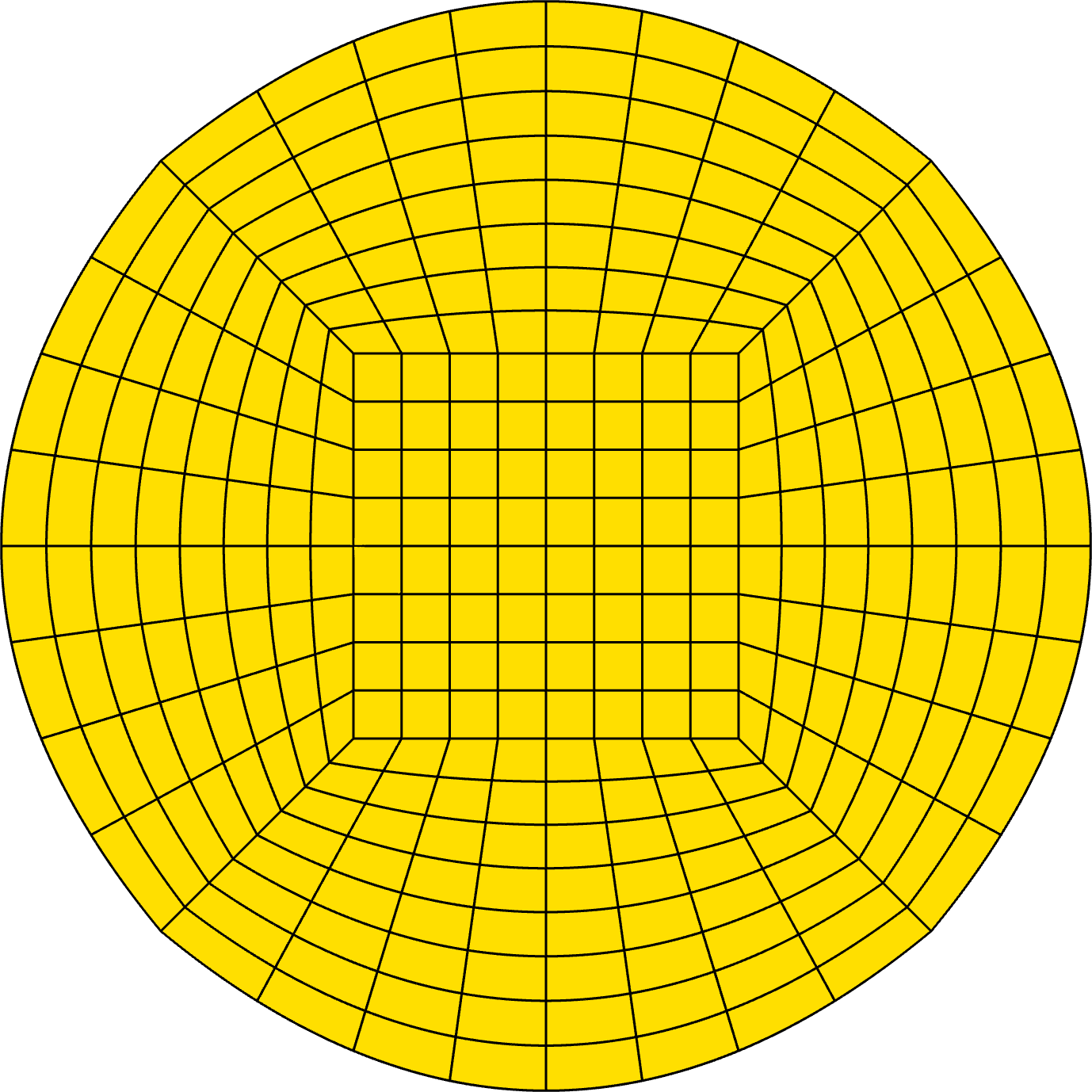}}
		\qquad & \qquad 
		\raisebox{-.5\height}{\includegraphics[width=0.205\linewidth]{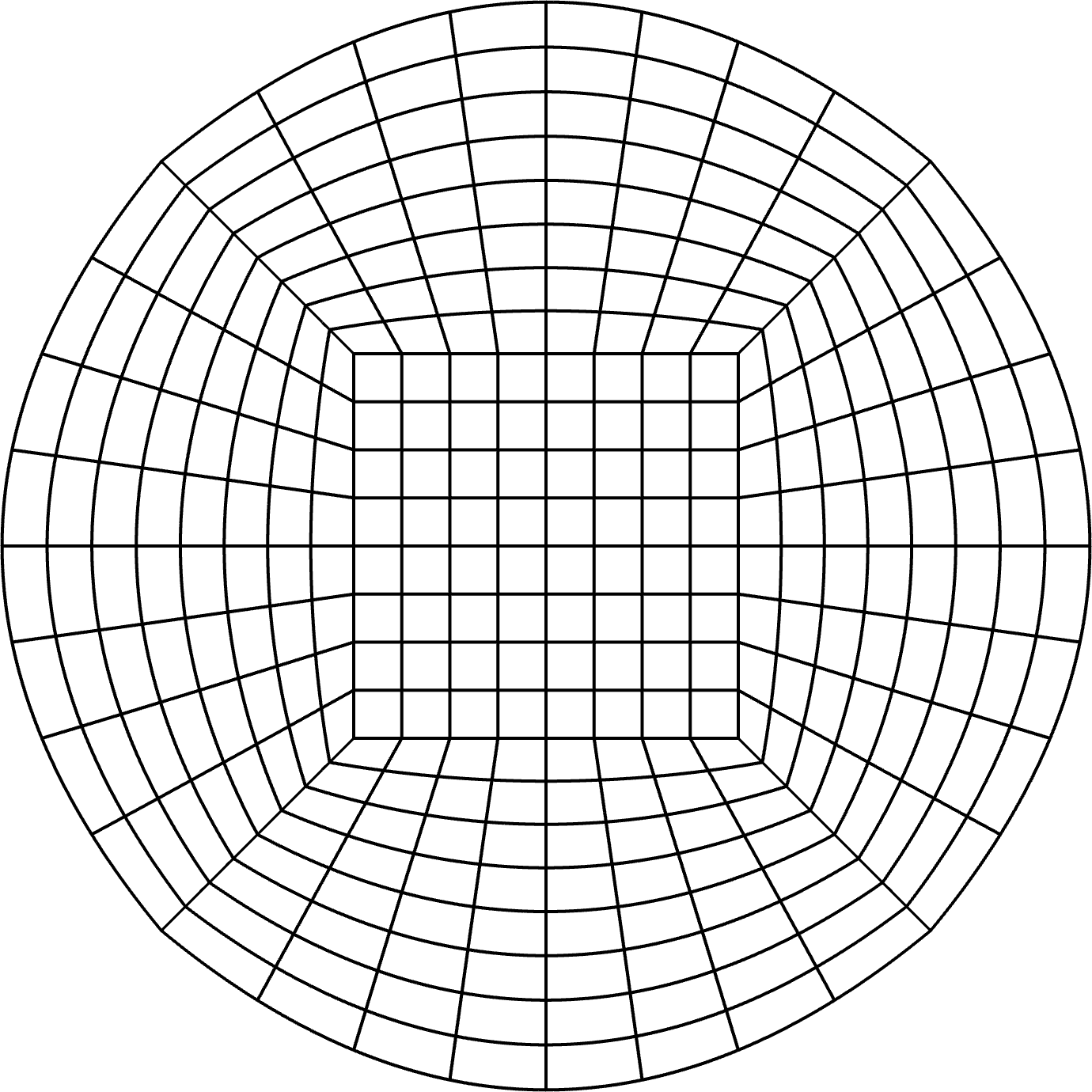}} 
		\\[.4cm]
		$t = 4.0$
		\qquad & \qquad
		\raisebox{-.5\height}{\includegraphics[width=0.205\linewidth]{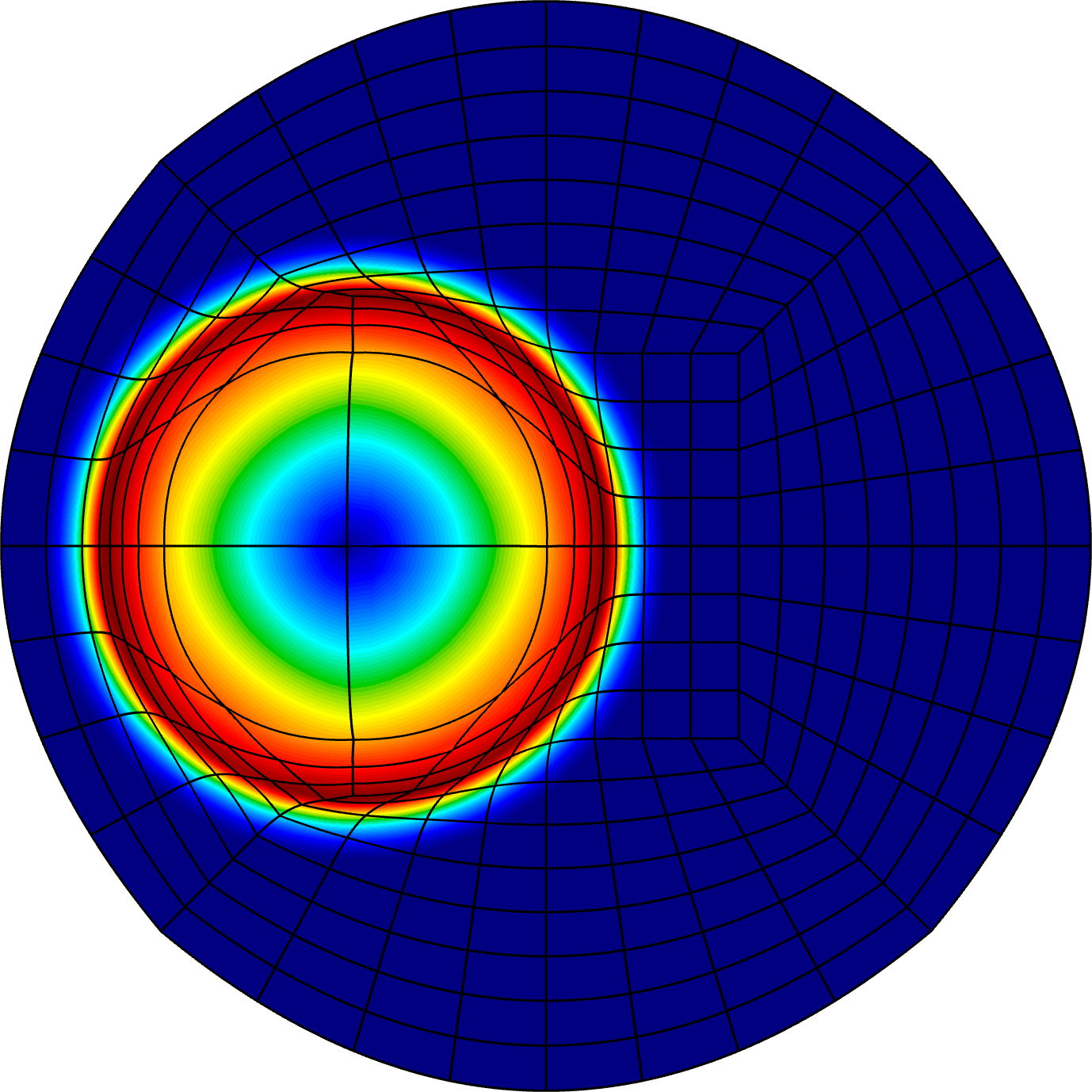}}
			\qquad & \qquad
		\raisebox{-.5\height}{\includegraphics[width=0.205\linewidth]{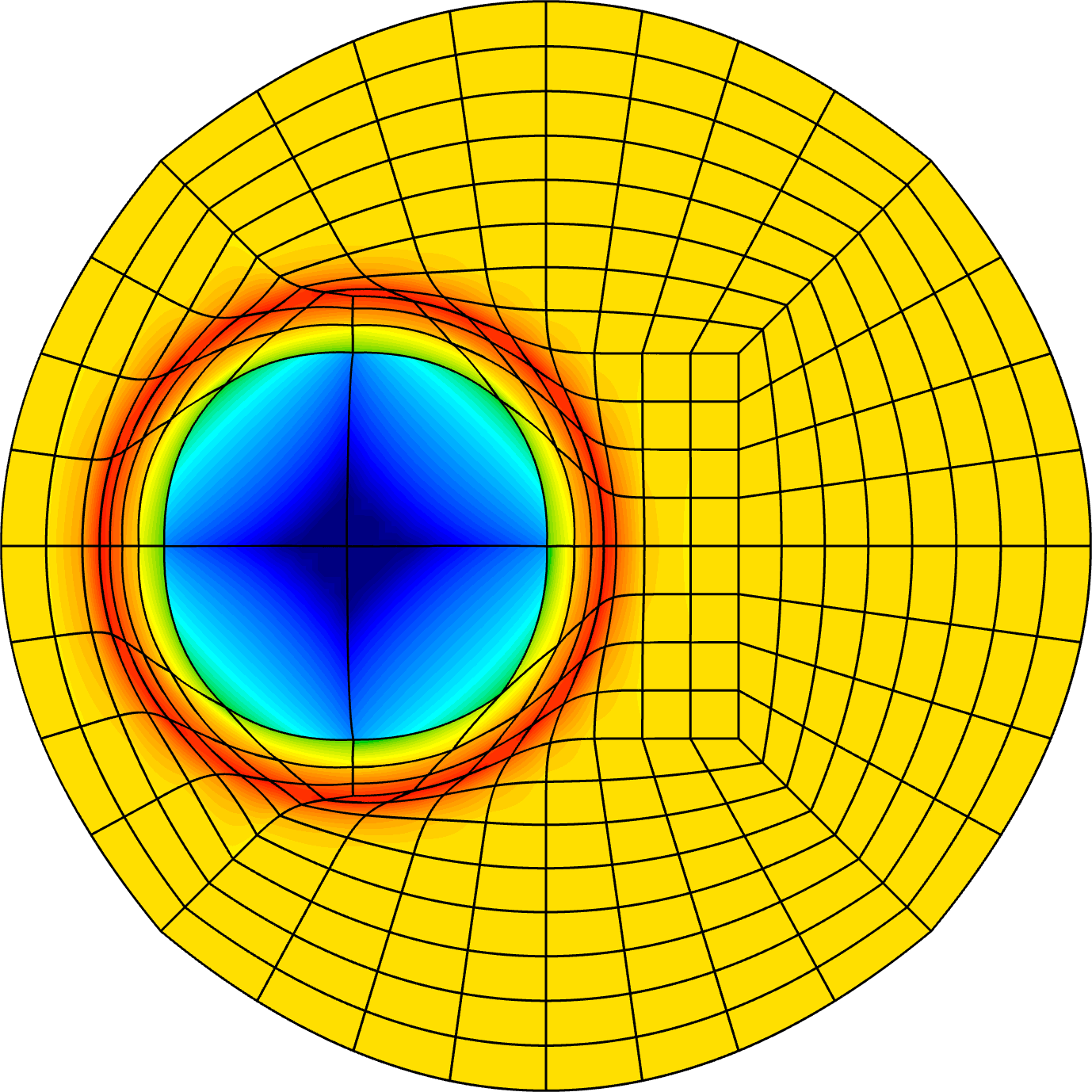}}
		  	\qquad & \qquad 
		\raisebox{-.5\height}{\includegraphics[width=0.205\linewidth]{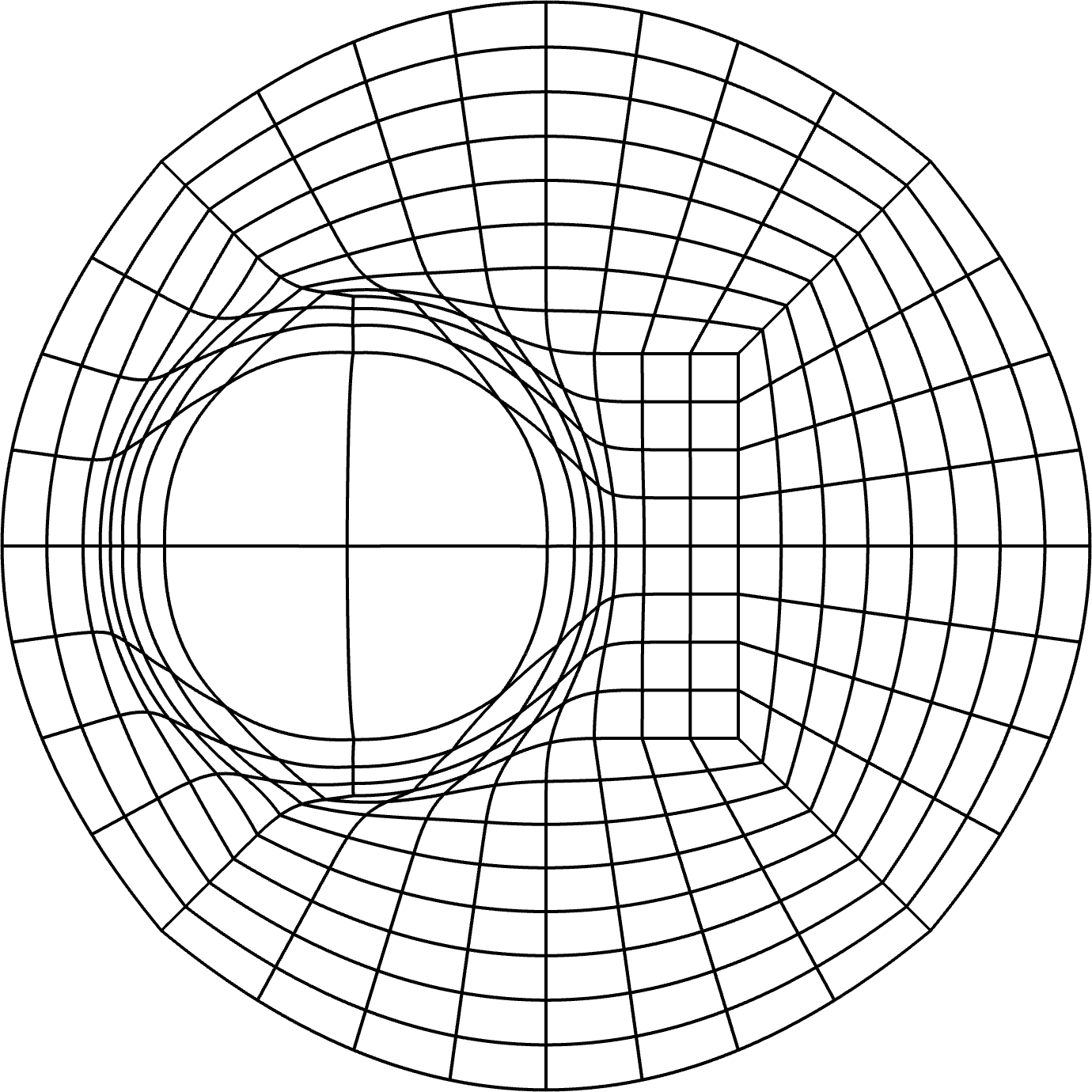}} 
		 \\[.4cm]
		 $t = 8.0$
		 \qquad & \qquad
		 \raisebox{-.5\height}{\includegraphics[width=0.205\linewidth]{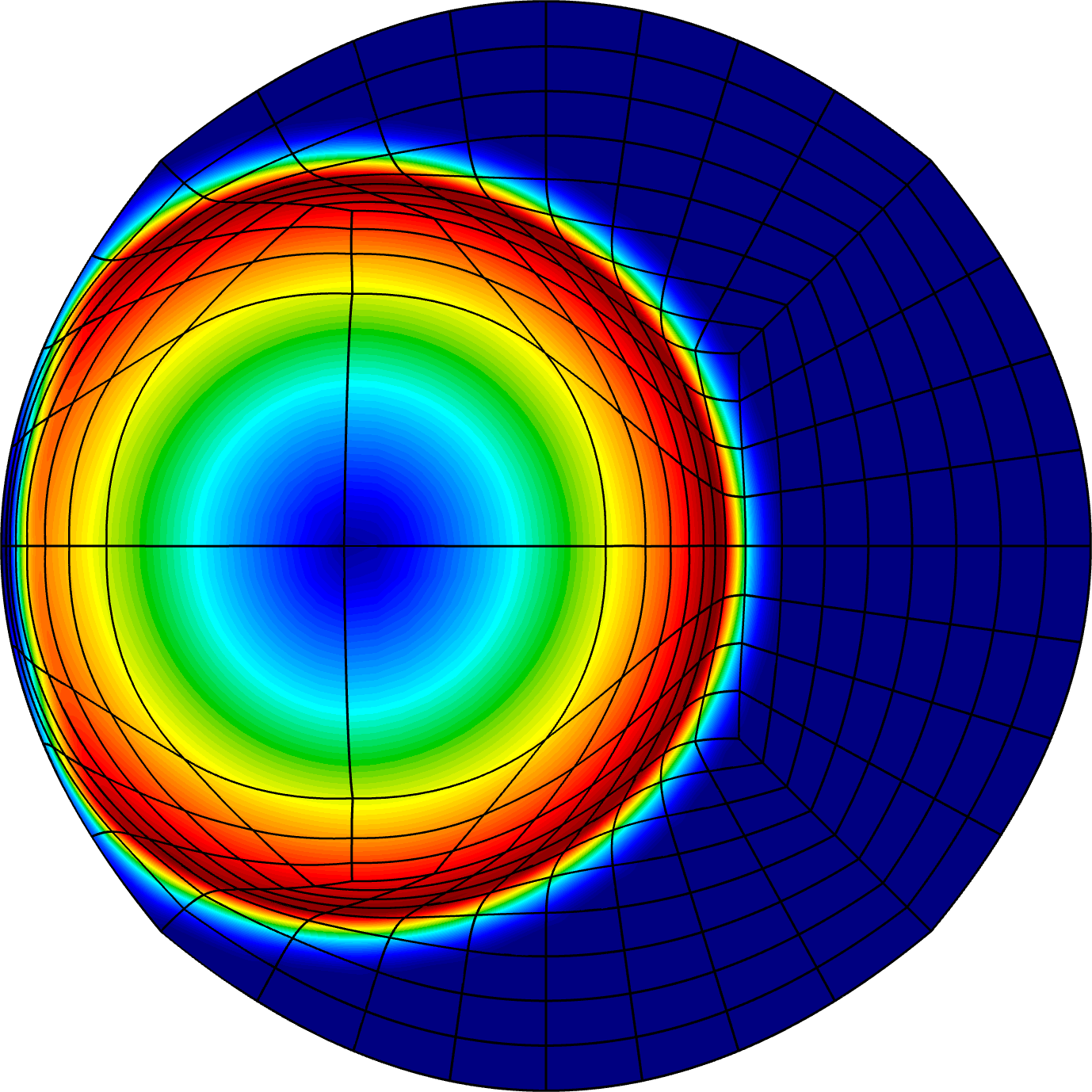}}
		 \qquad & \qquad
		 \raisebox{-.5\height}{\includegraphics[width=0.205\linewidth]{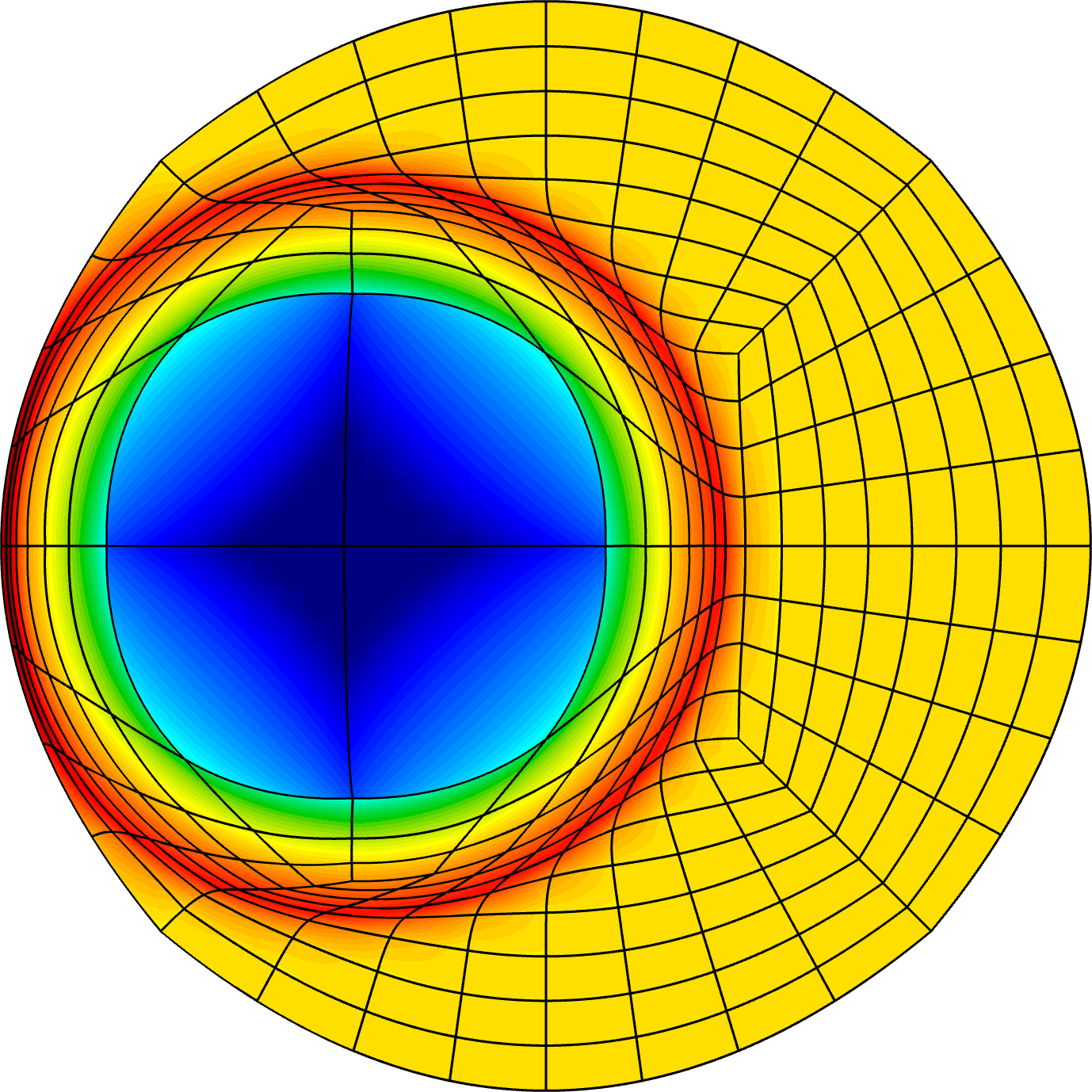}}
		 \qquad & \qquad 
		 \raisebox{-.5\height}{\includegraphics[width=0.205\linewidth]{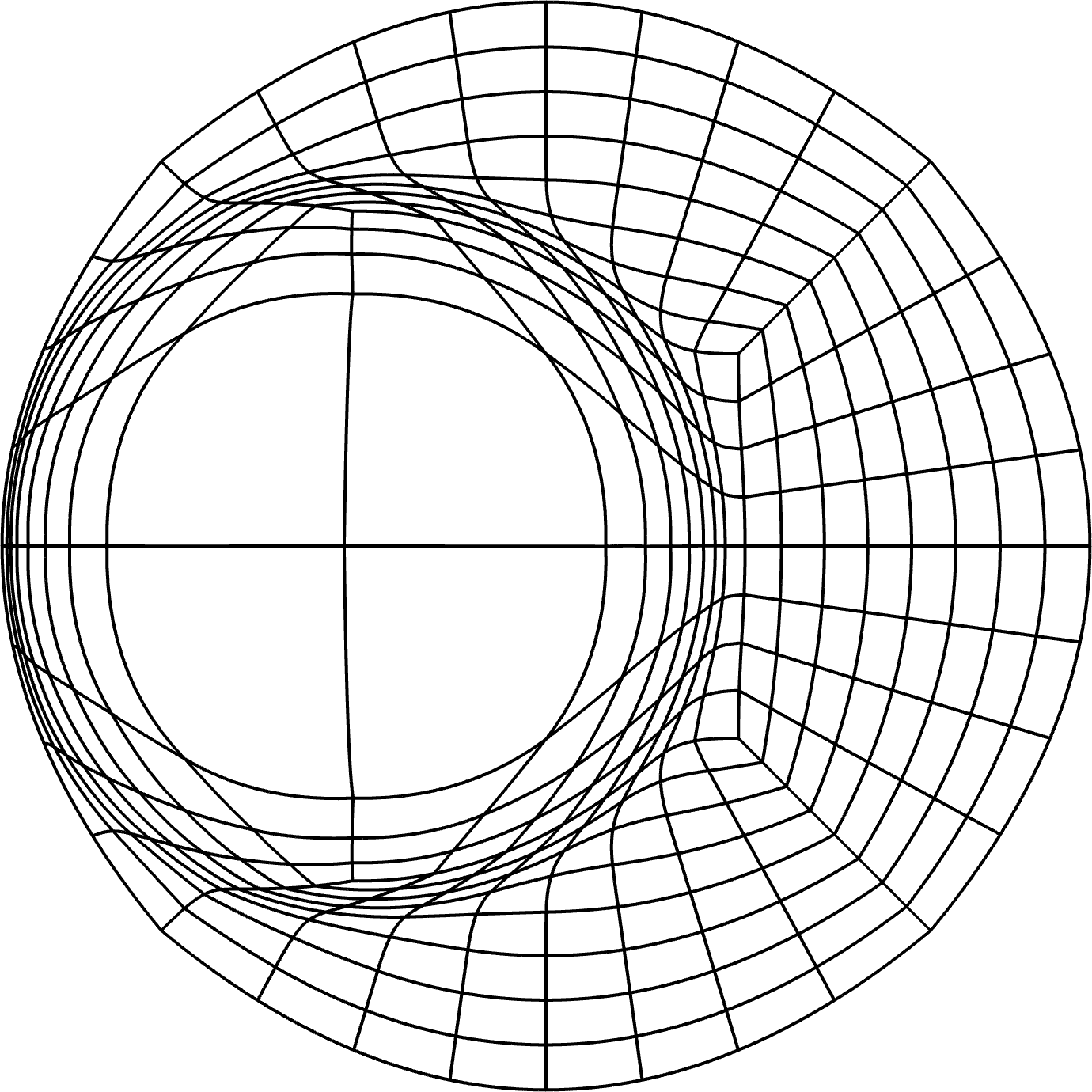} }
		 	 \\[.4cm]
		 	 		 $t = 12.0$
		 	 \qquad & \qquad
		 	 \raisebox{-.5\height}{\includegraphics[width=0.205\linewidth]{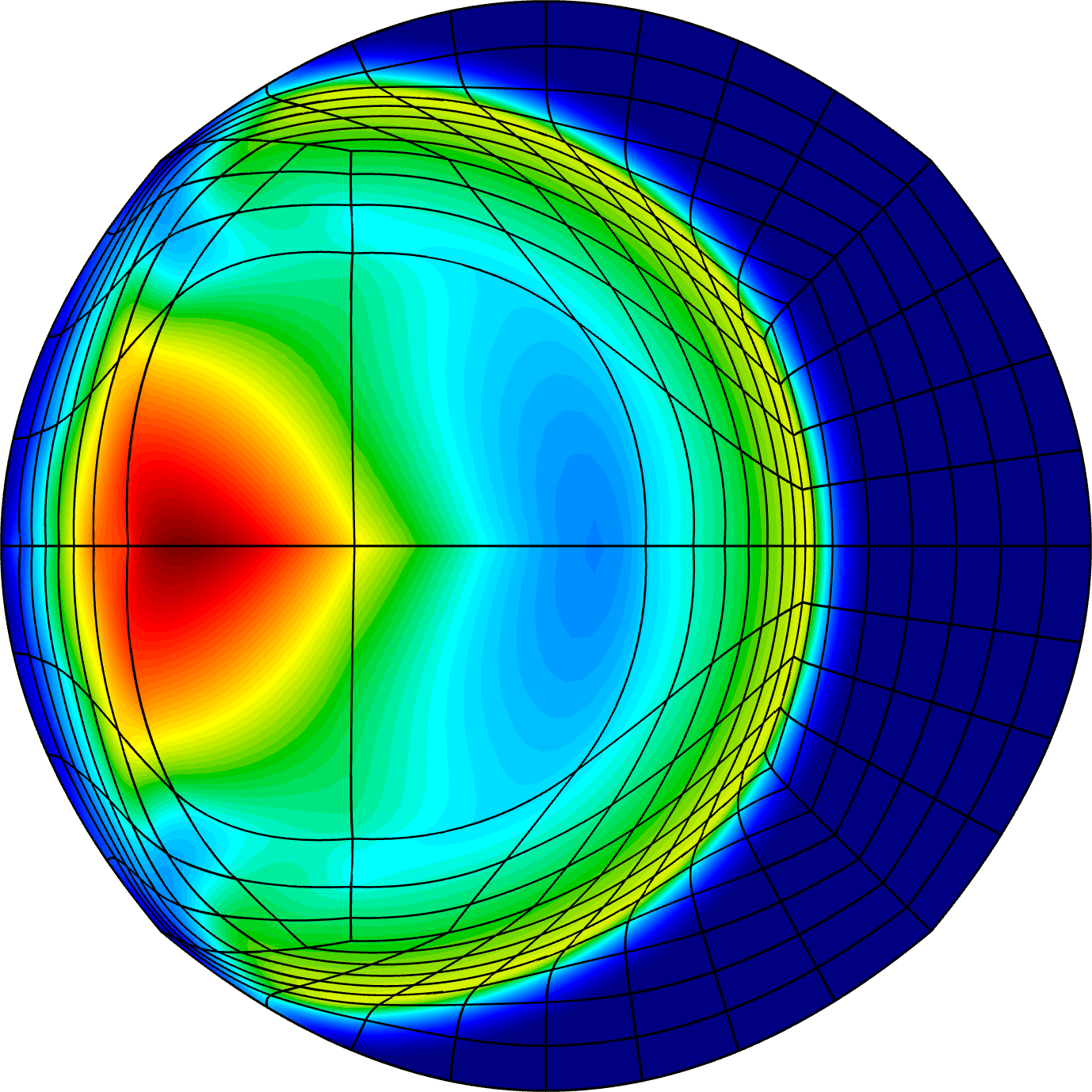}}
		 	 \qquad & \qquad
		 	 \raisebox{-.5\height}{\includegraphics[width=0.205\linewidth]{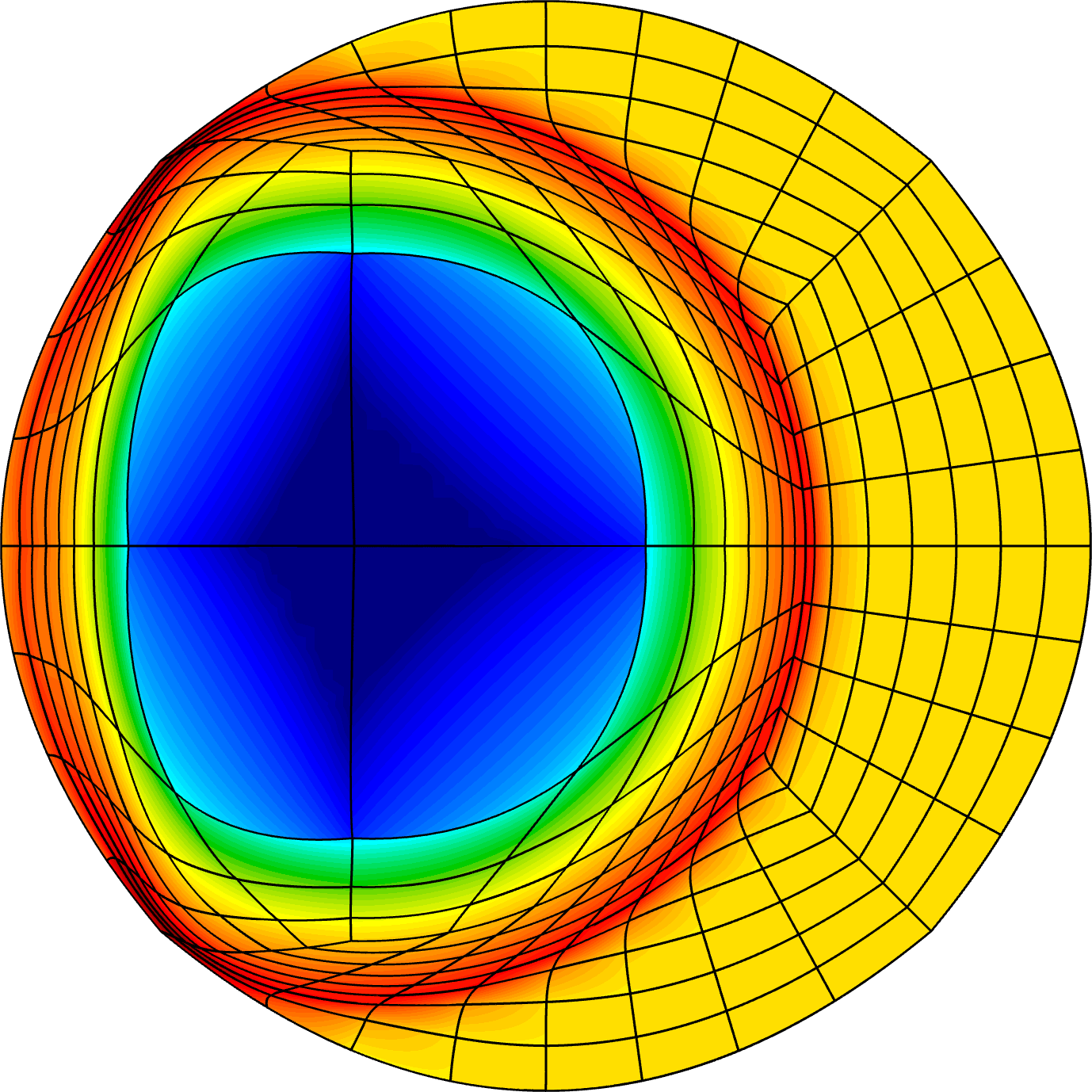}}
		 	 \qquad & \qquad 
		 	 \raisebox{-.5\height}{\includegraphics[width=0.205\linewidth]{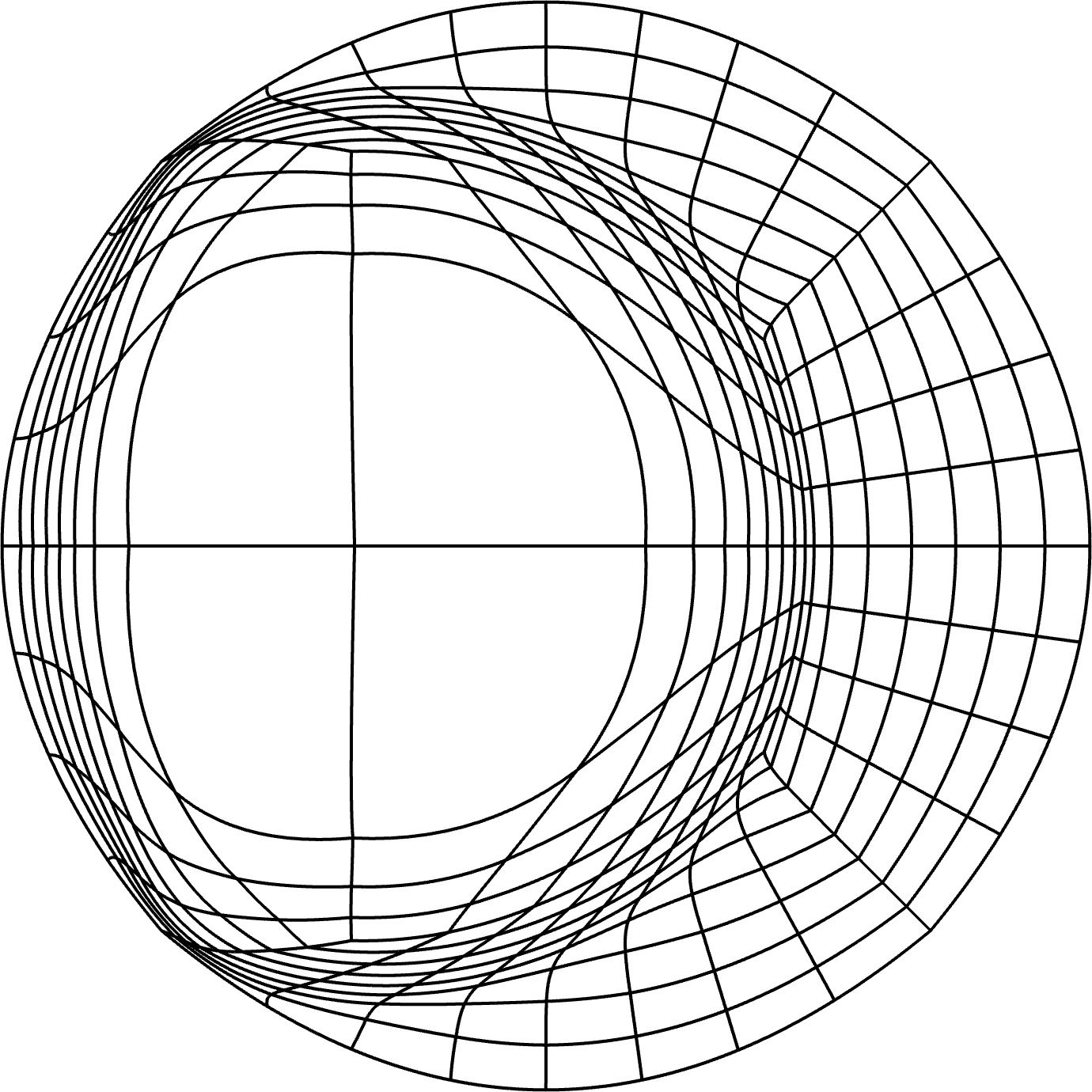} }
		 	 \\[.4cm]
		 	  	 		 $t = 16.0$
		 	 \qquad & \qquad
		 	 \raisebox{-.5\height}{\includegraphics[width=0.205\linewidth]{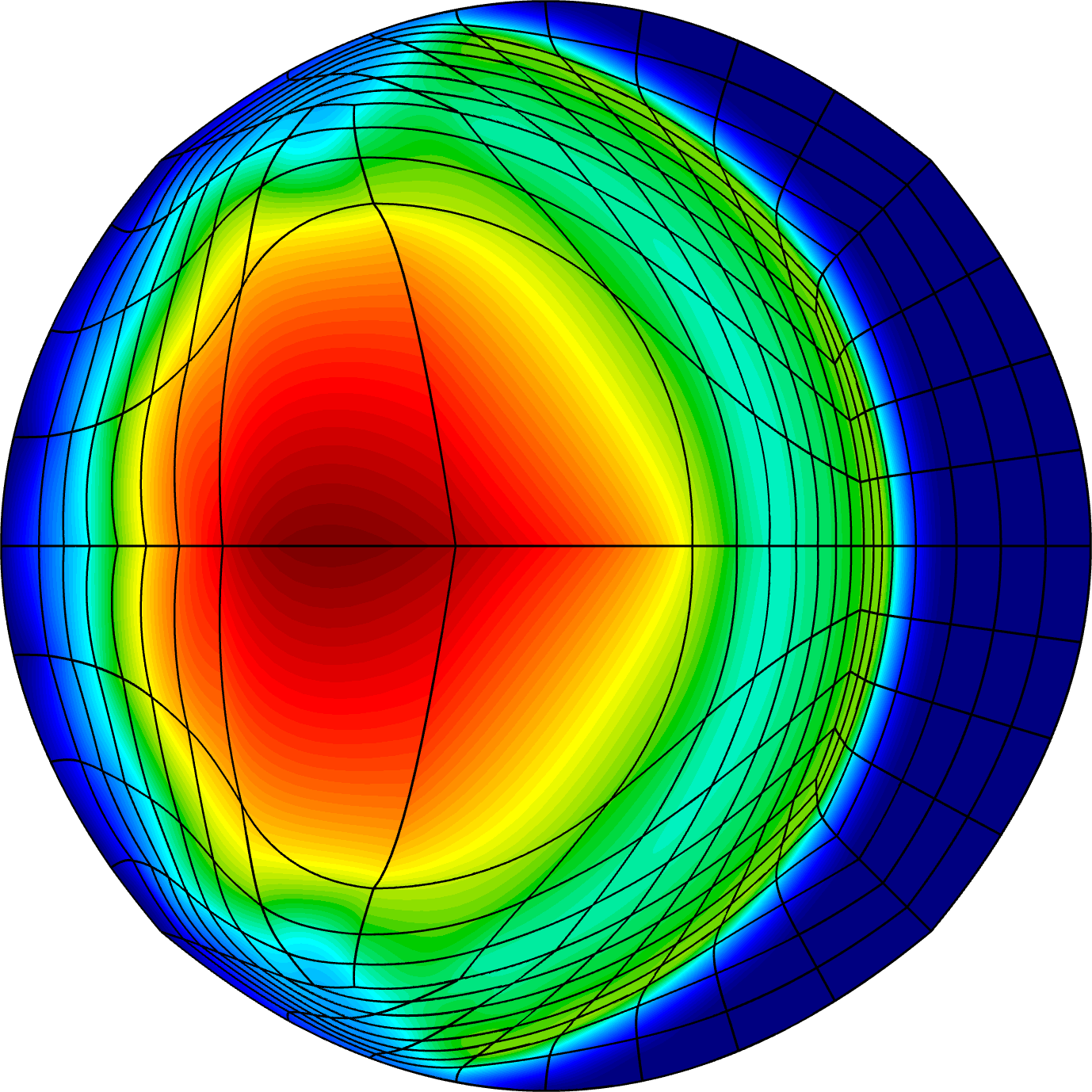}}
		 	 \qquad & \qquad
		 	 \raisebox{-.5\height}{\includegraphics[width=0.205\linewidth]{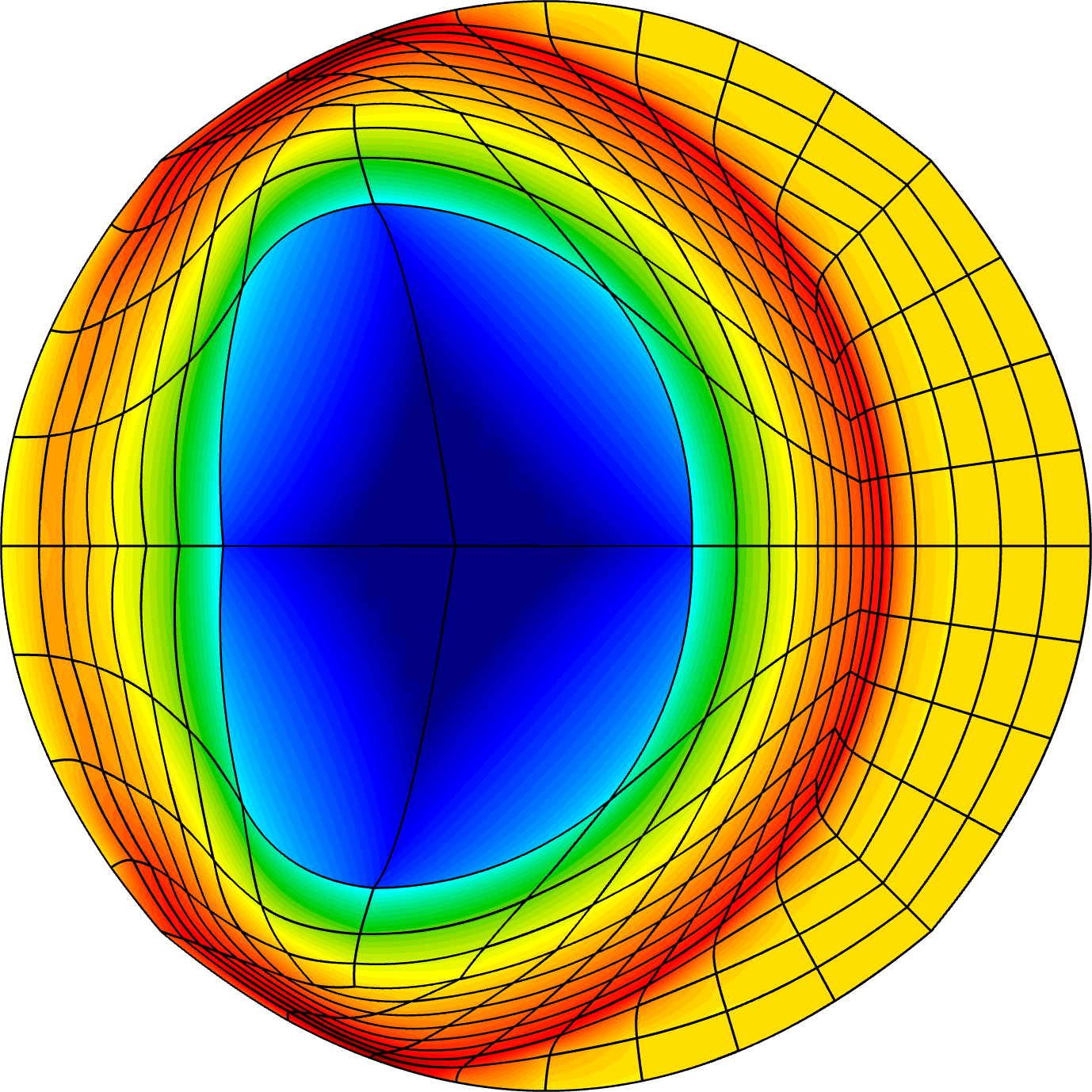}}
		 	 \qquad & \qquad 
		 	 \raisebox{-.5\height}{\includegraphics[width=0.205\linewidth]{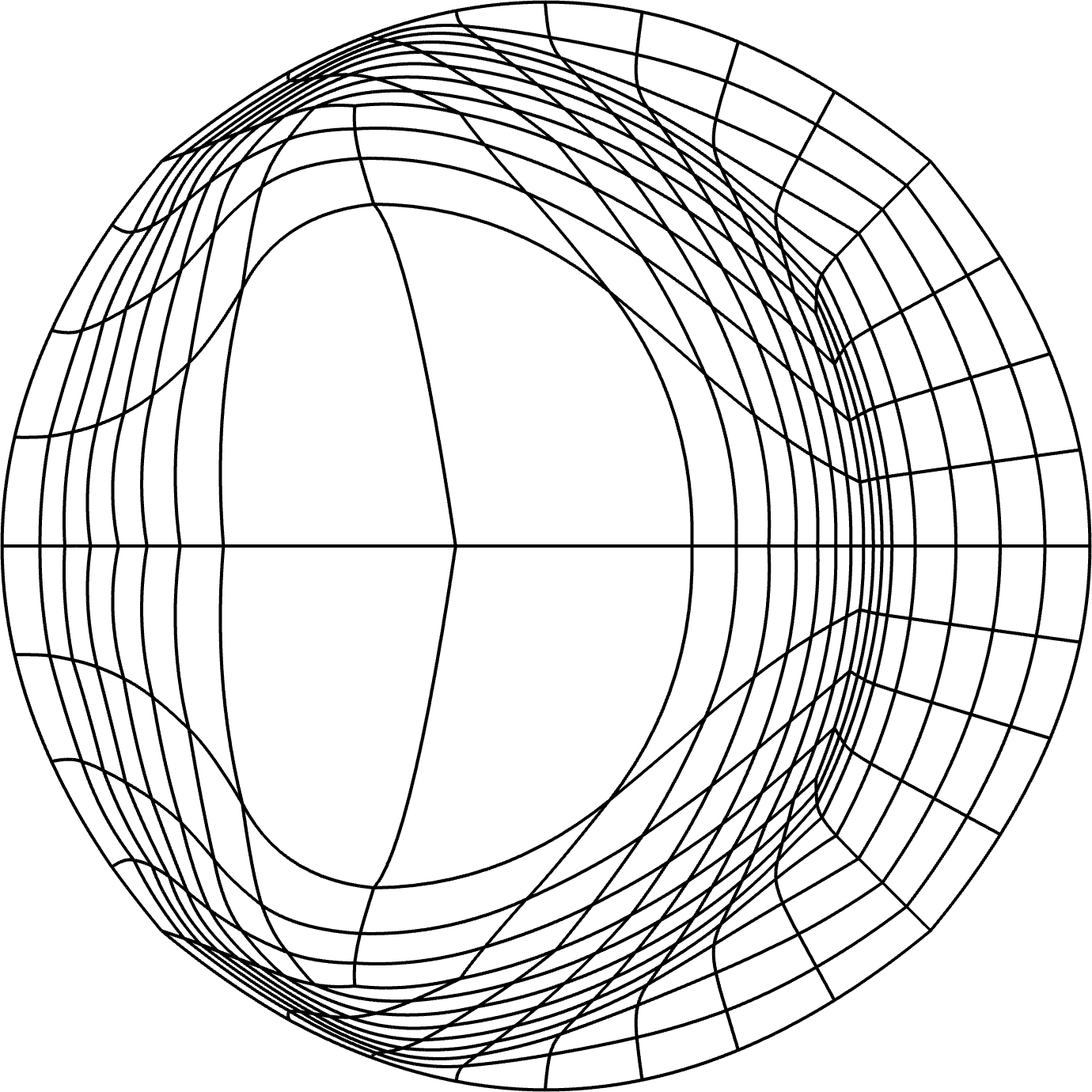} }
		 	 \\[.4cm]
		 	 $t = 20.0$
		 	 \qquad & \qquad
\raisebox{-.5\height}{\includegraphics[width=0.205\linewidth]{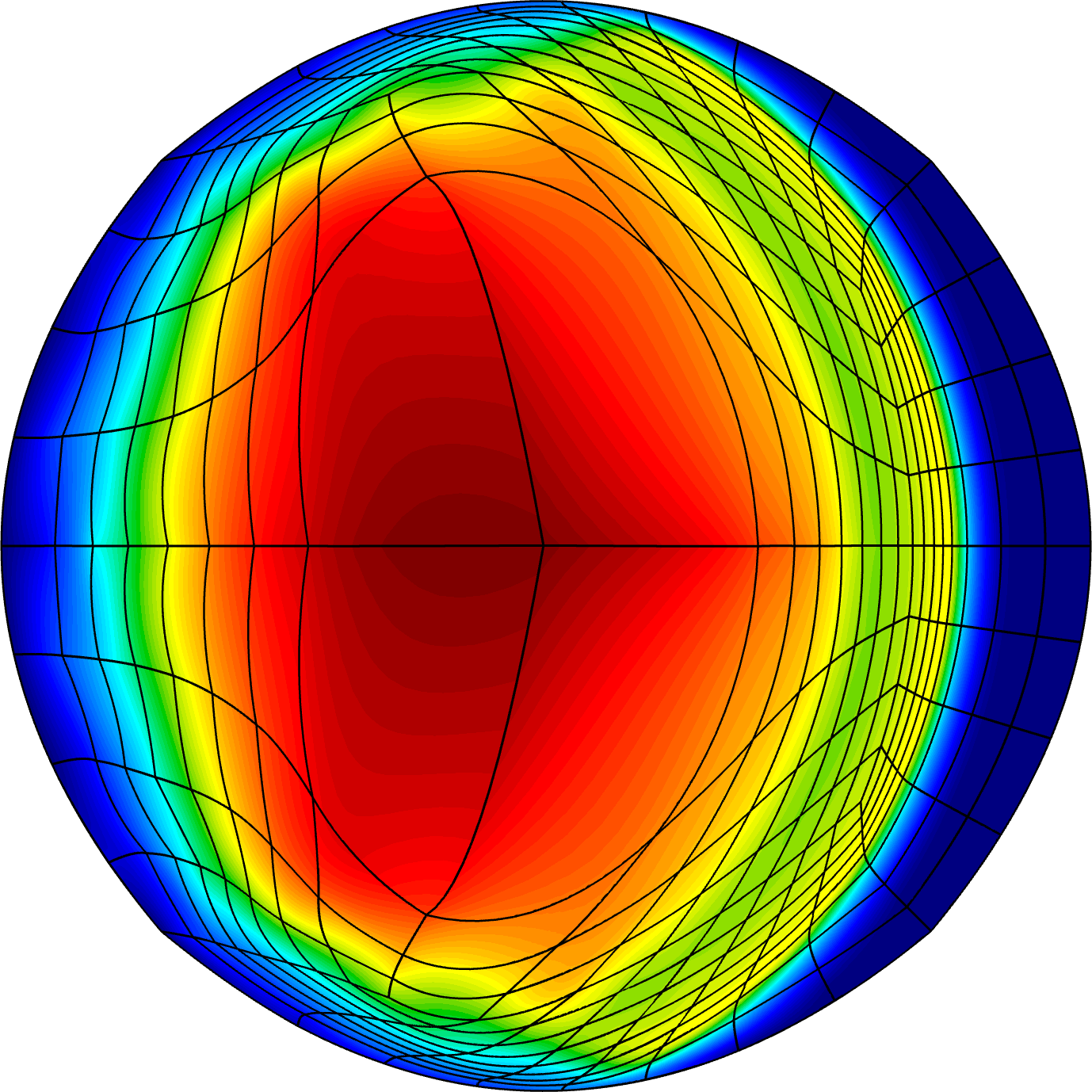}}
\qquad & \qquad
\raisebox{-.5\height}{\includegraphics[width=0.205\linewidth]{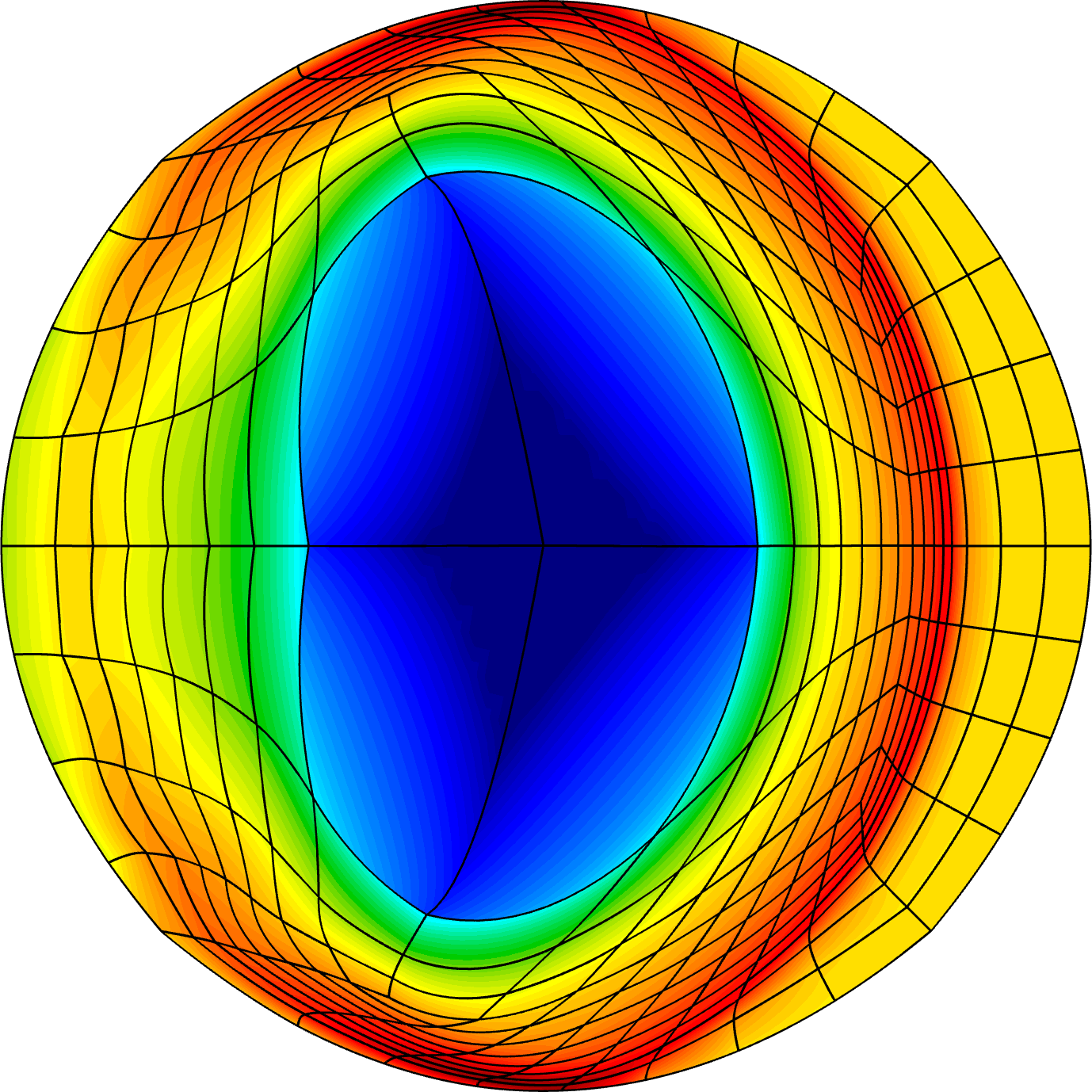}}
\qquad & \qquad 
\raisebox{-.5\height}{\includegraphics[width=0.205\linewidth]{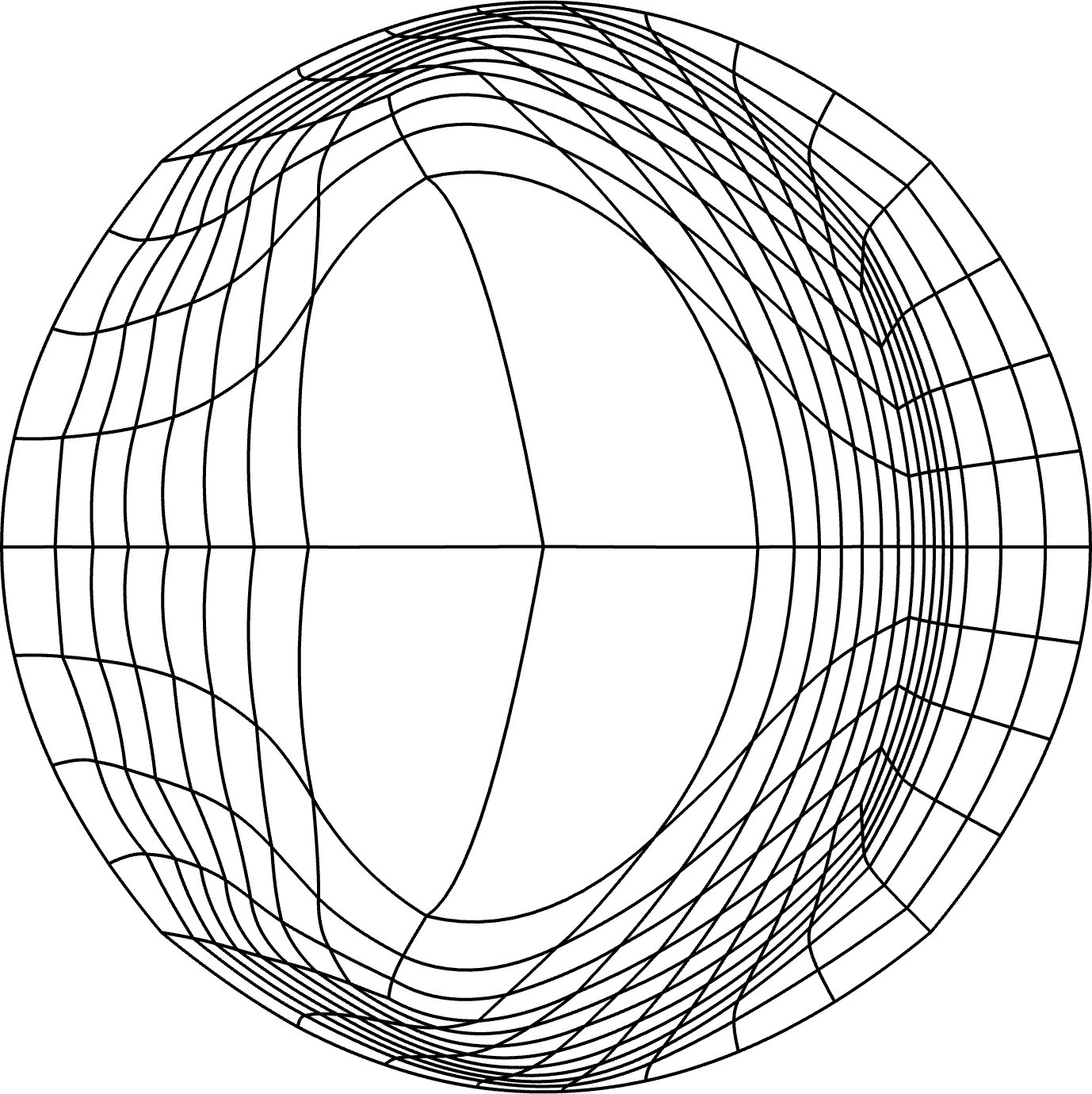} }
	  \\[.0cm]
	\end{tabular}
	\caption{Plots of the velocity and density fields in addition to the mesh deformation for the Sedov test in a circular domain using a $Q_{2}-Q_{1}$ velocity-energy pair.}
	\label{fig:CircleSedovresults}
\end{figure}

We perform the Sedov test in a circular domain and show plots of the velocity and density fields in addition to the mesh deformation at the different points-in-time until $t = 20.0$ for $Q_{2}-Q_{1}$ velocity-energy pair as shown in Figure~\ref{fig:CircleSedovresults}.
Observe that the shock slides correctly along the outside curved boundary without shape distortions.
Similarly, the curved shape of the boundary is not affected by the strong shock bounce.
The solution appears physically correct without any unphysical oscillations.

\subsection{Three-dimensional Sedov explosion in a cube}
\label{cube_sedov_3d}
\begin{figure}[tb]
	\centering
	\begin{tabular}{ccc}
		\multicolumn{3}{c}{Velocity}  \\
		\includegraphics[width=0.275\linewidth]{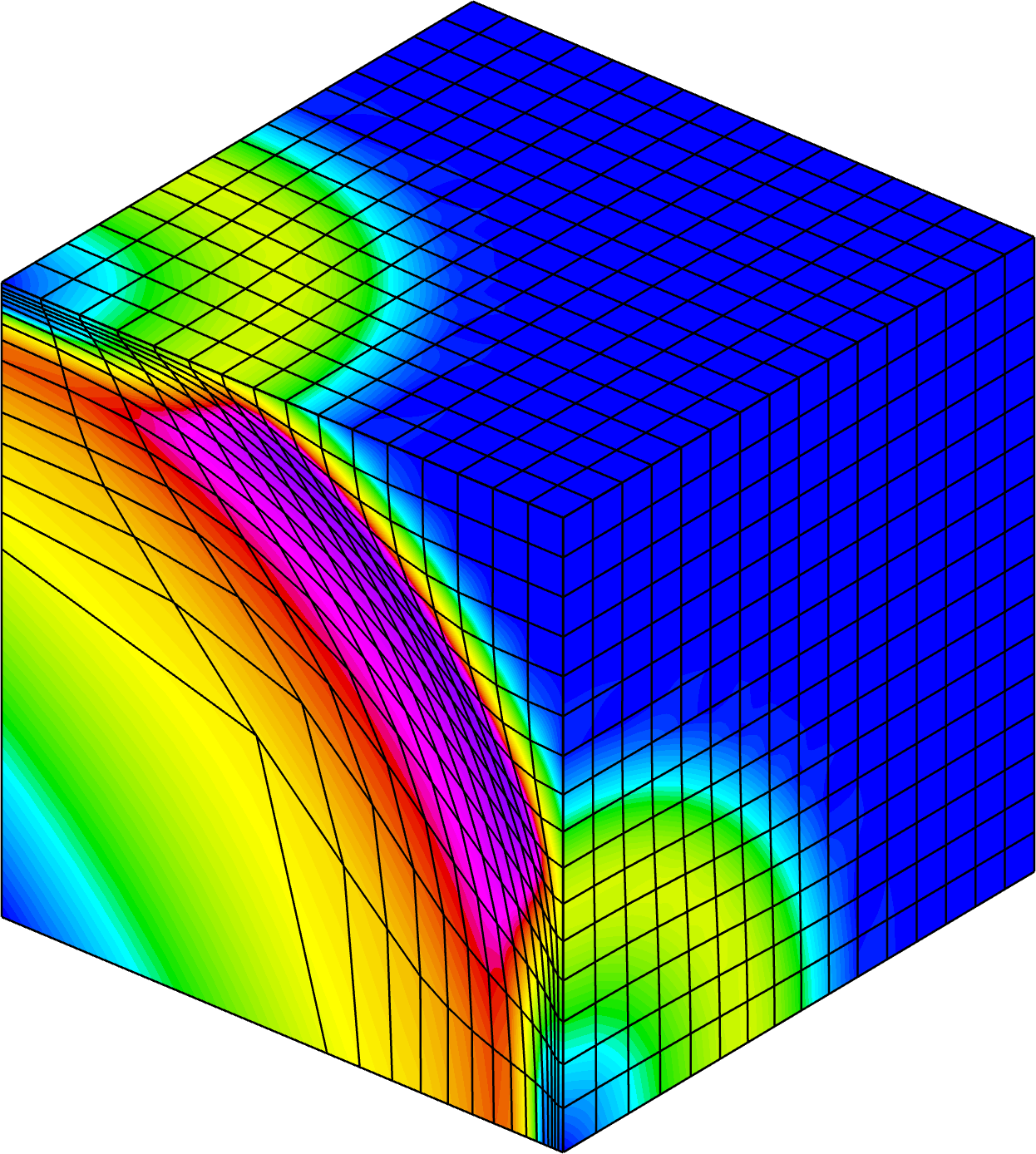} 
		\qquad & \qquad  
		\includegraphics[width=0.275\linewidth]{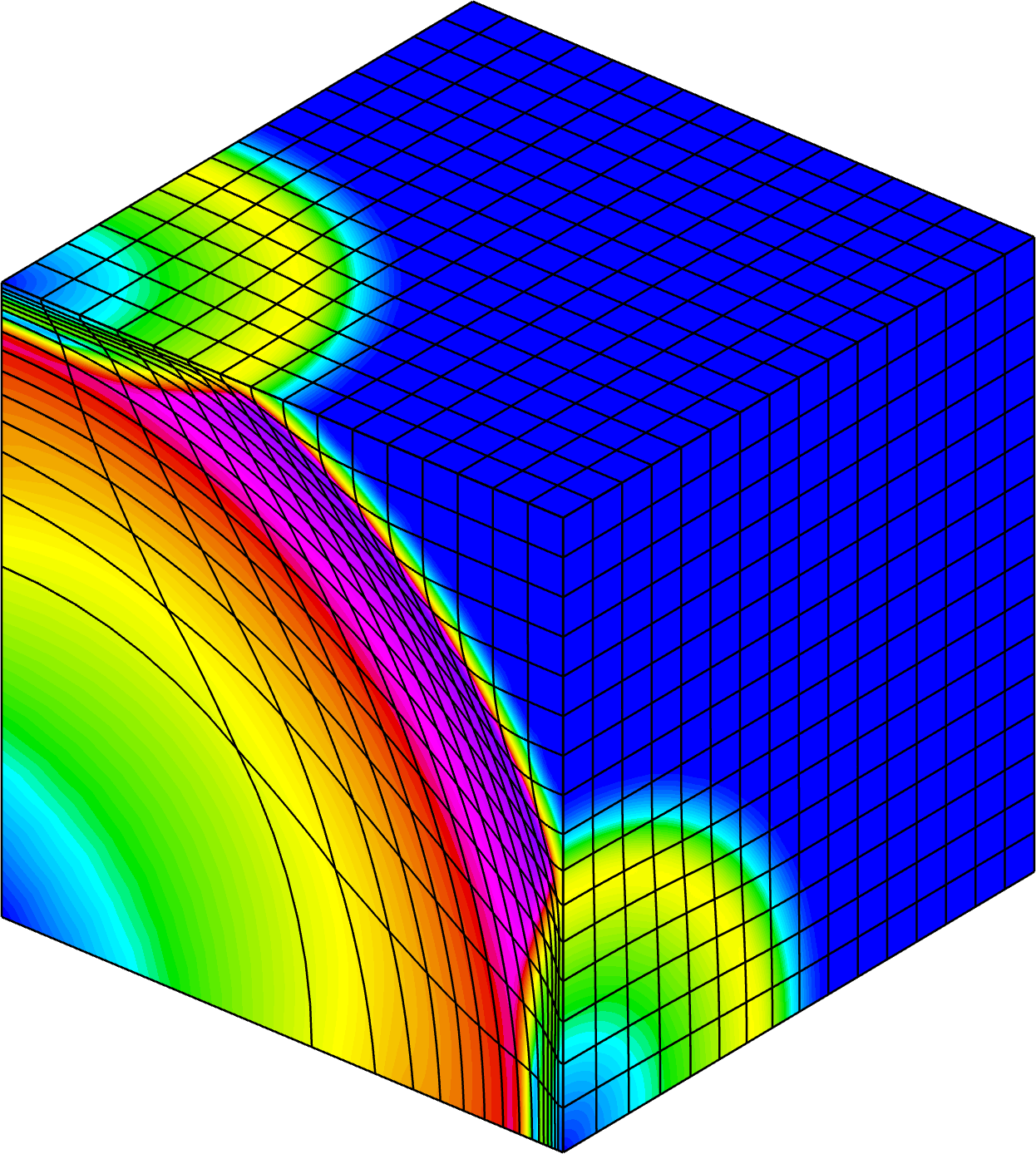}
		\qquad & \qquad  
		\includegraphics[width=0.275\linewidth]{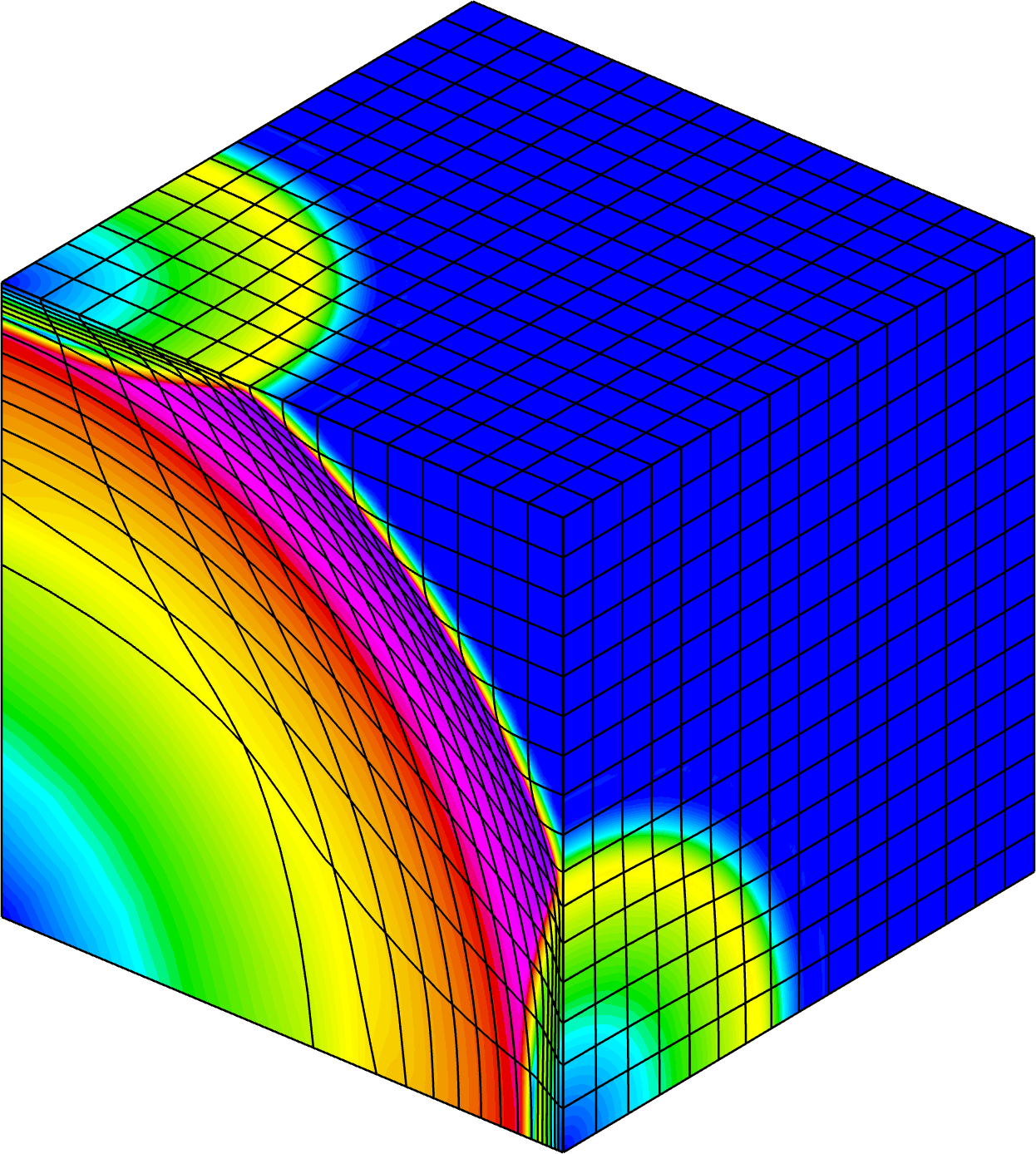} \\[.4cm]
		\includegraphics[width=0.275\linewidth]{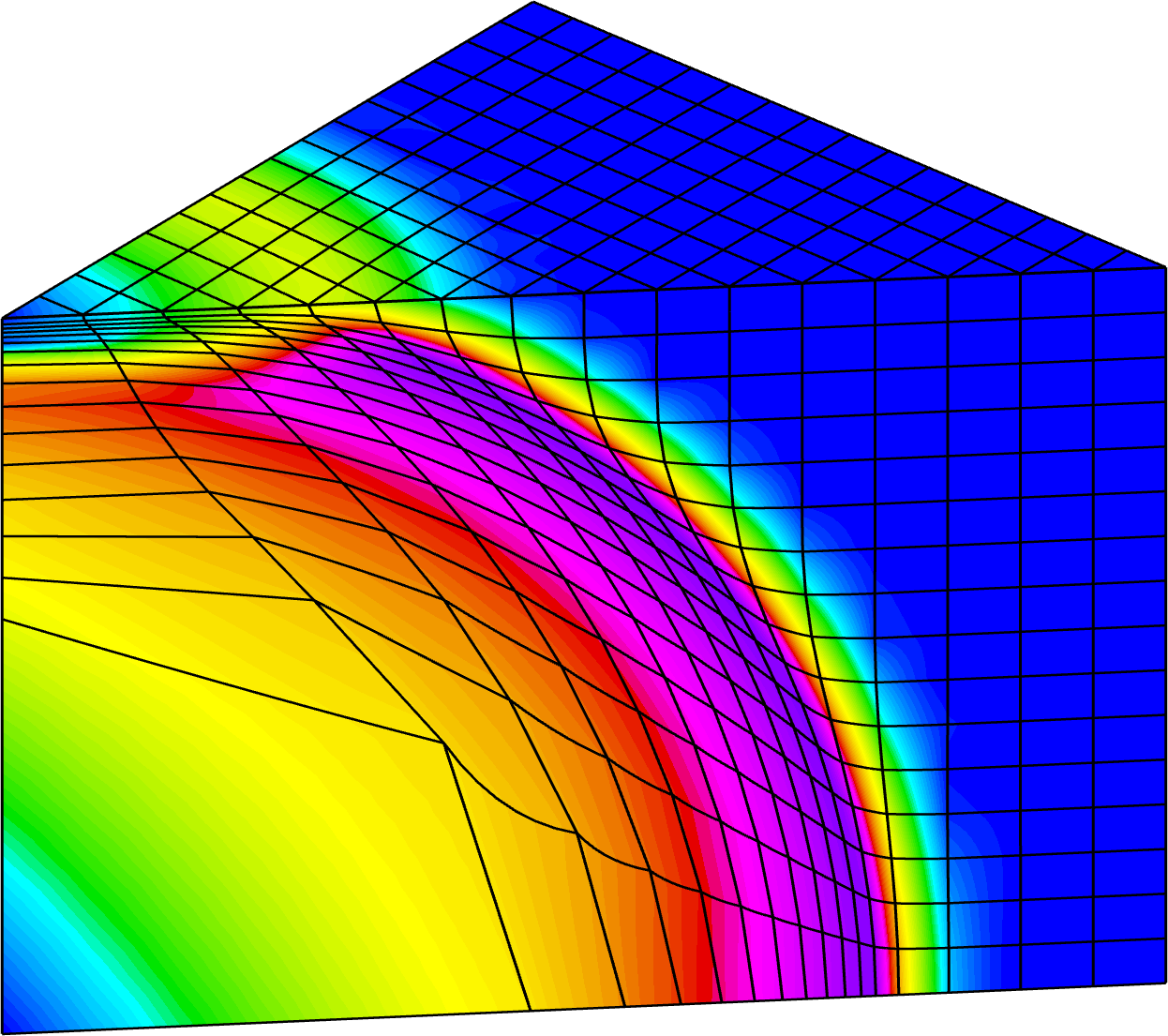} 
		\qquad & \qquad  
		\includegraphics[width=0.275\linewidth]{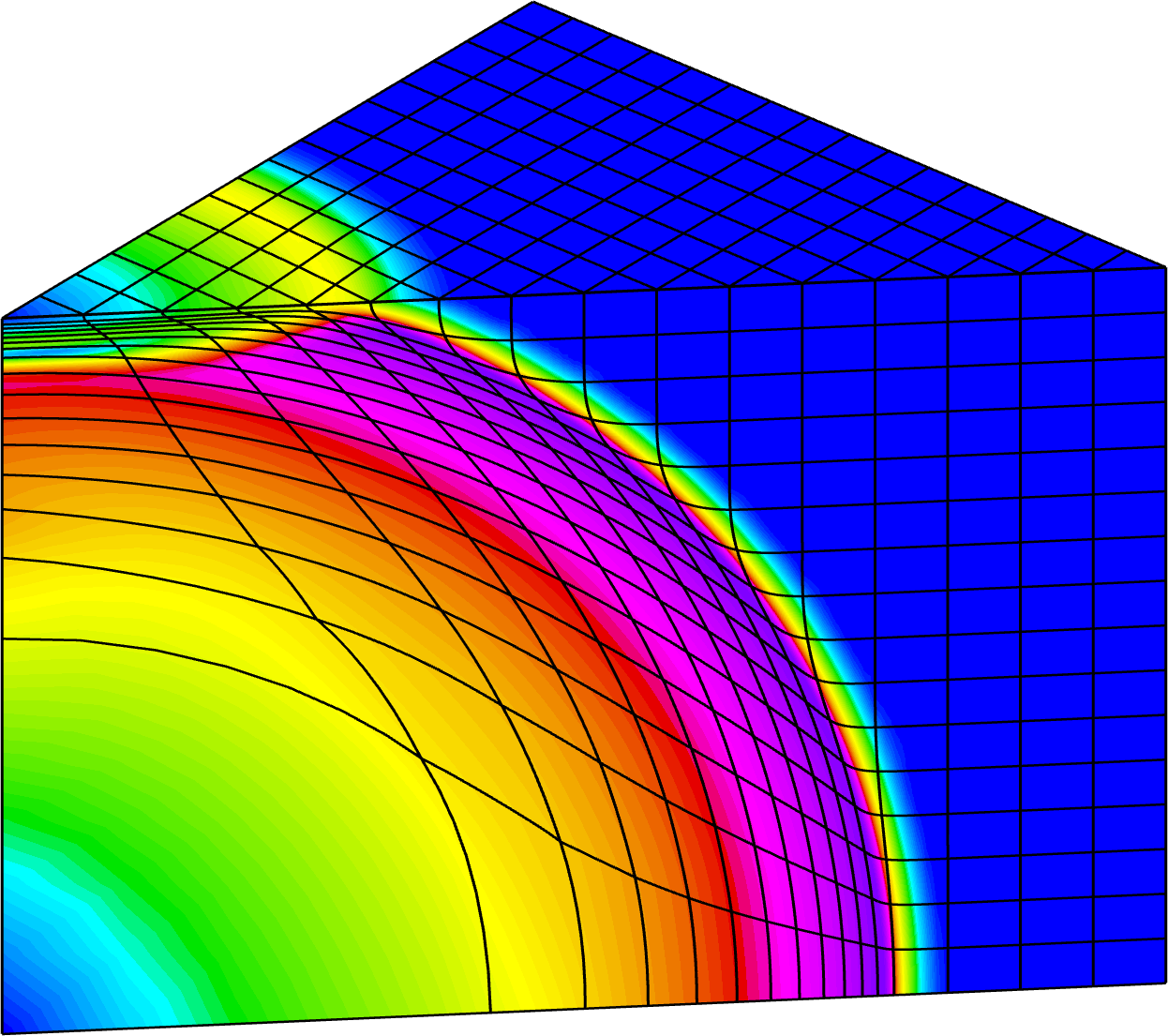}
		\qquad & \qquad  
		\includegraphics[width=0.275\linewidth]{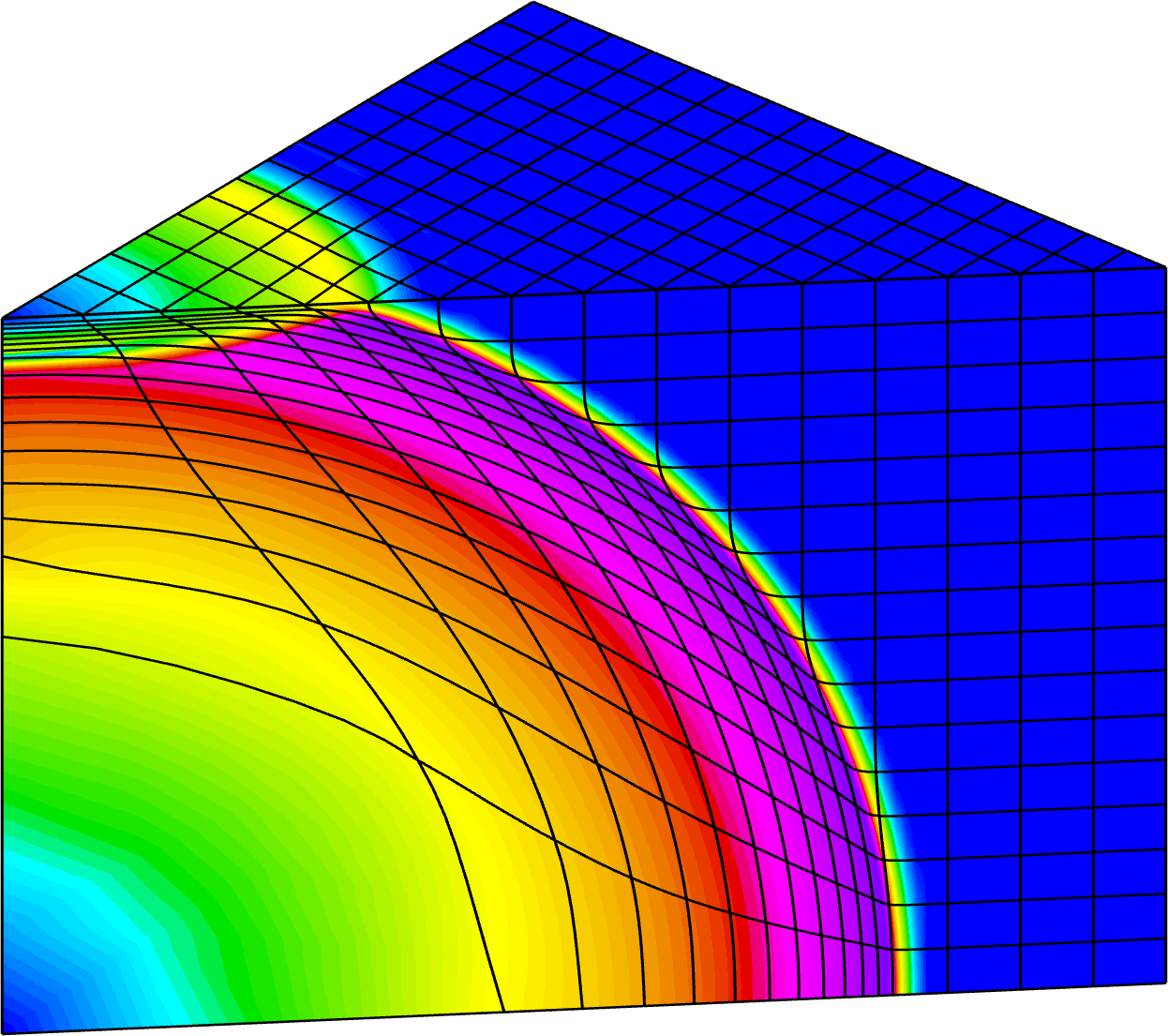}\\[.4cm]
		\multicolumn{3}{c}{Density}  \\
		\includegraphics[width=0.275\linewidth]{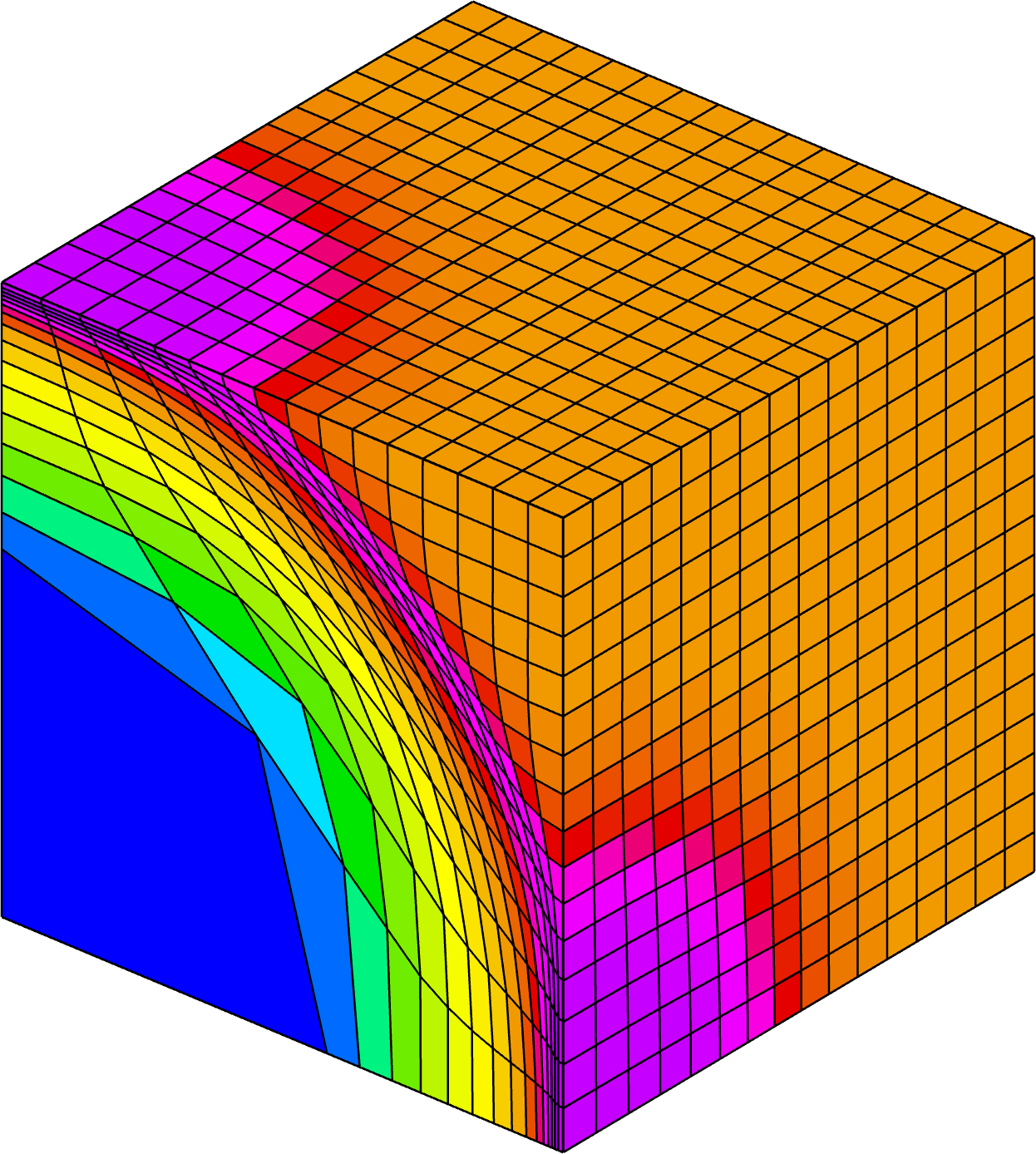} 
		\qquad & \qquad  
		\includegraphics[width=0.275\linewidth]{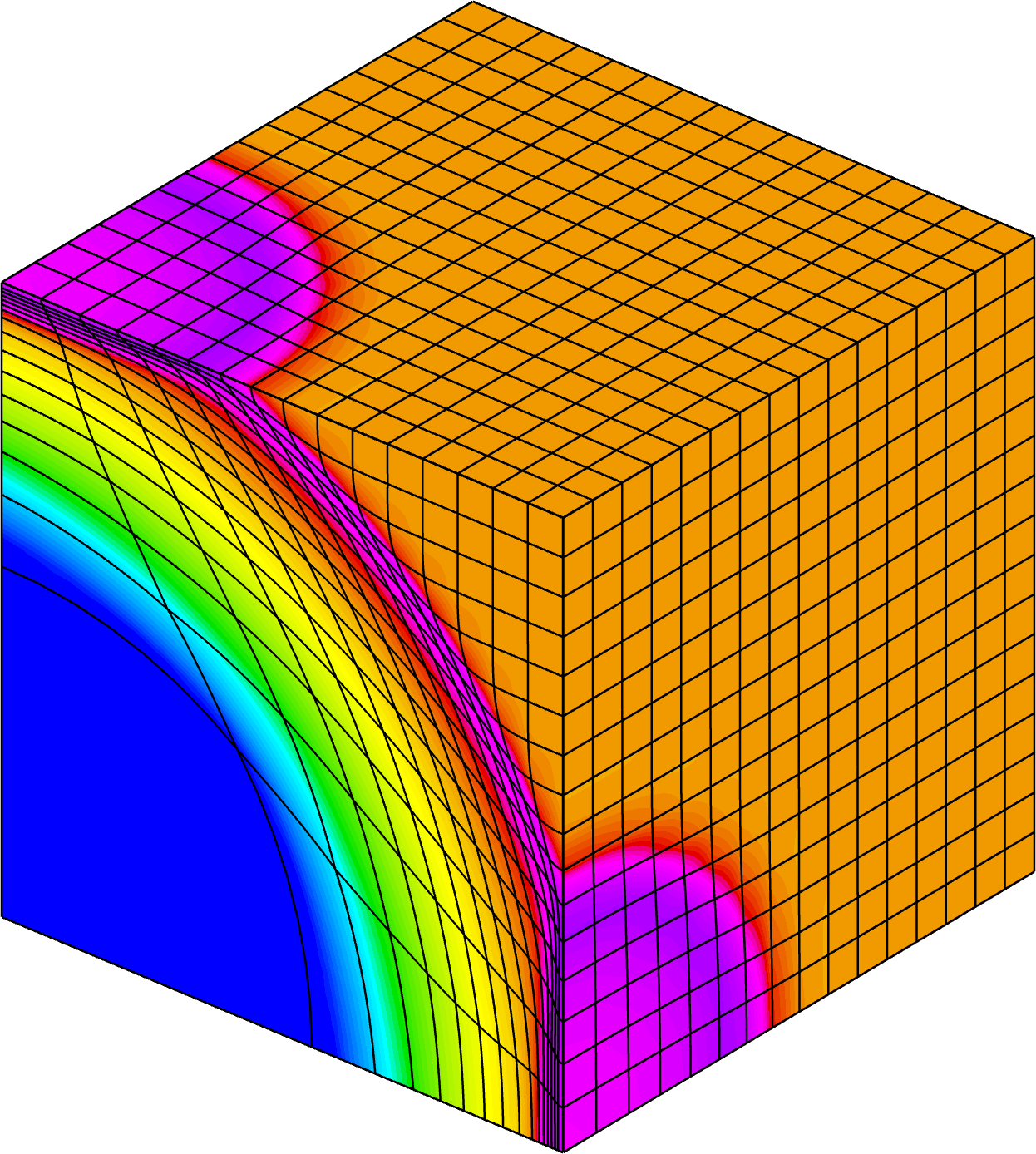}
		\qquad & \qquad  
		\includegraphics[width=0.275\linewidth]{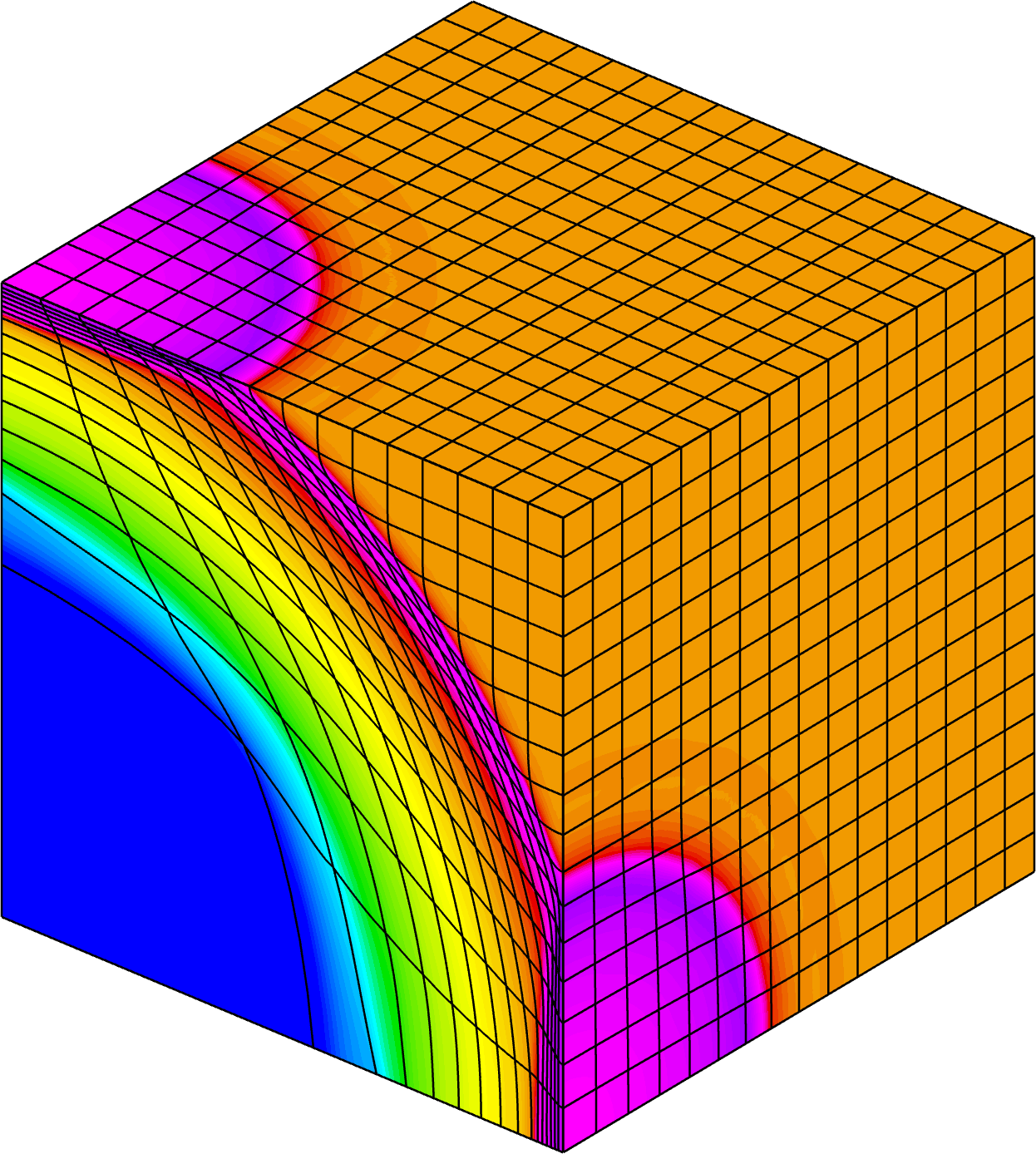}  \\[.4cm]
		\includegraphics[width=0.275\linewidth]{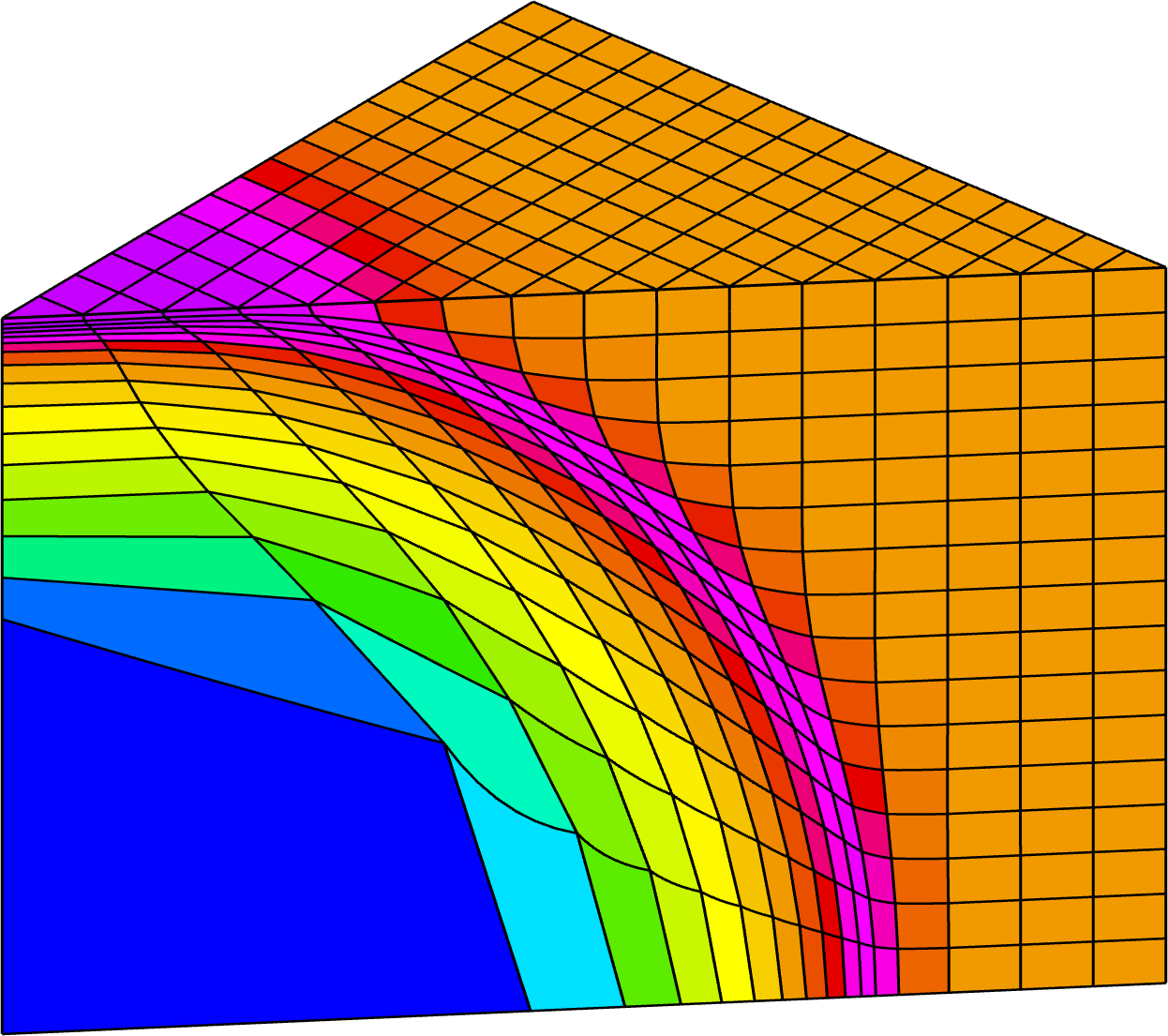} 
		\qquad & \qquad  
		\includegraphics[width=0.275\linewidth]{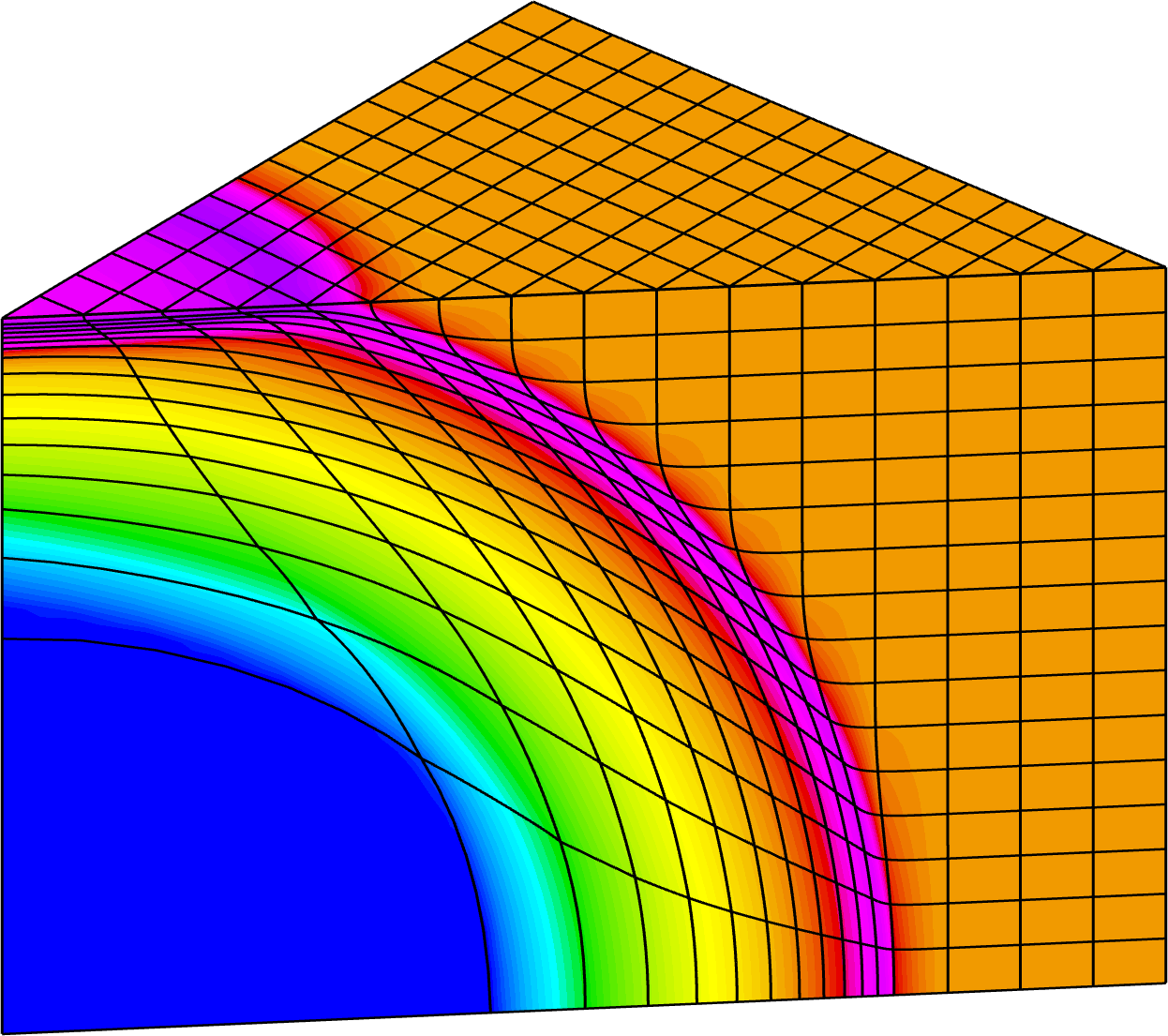}
		\qquad & \qquad  
		\includegraphics[width=0.275\linewidth]{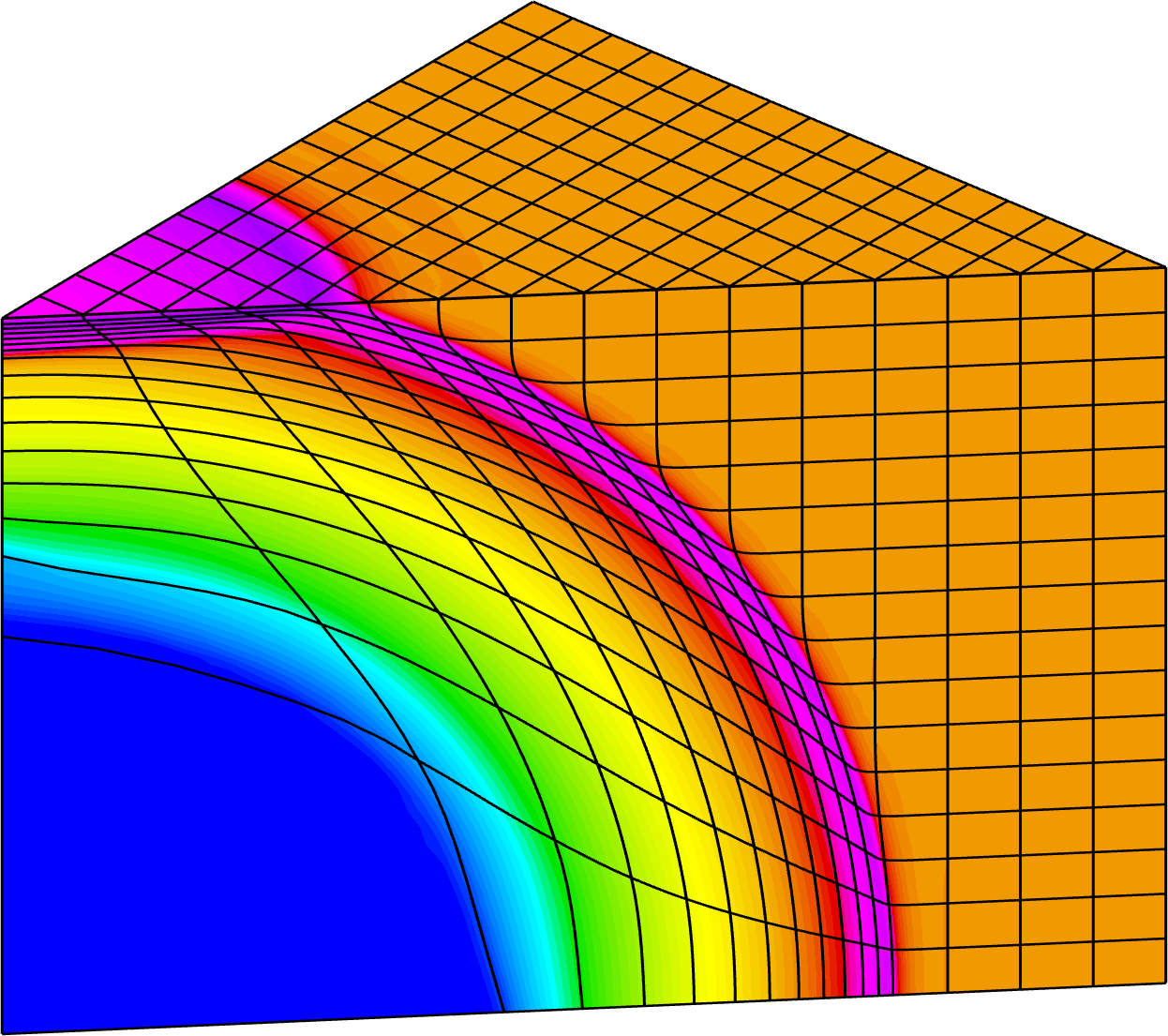} \\[0cm]
		$Q_{1}-Q_{0}$ 
		&
		$Q_{2}-Q_{1}$
		&
		$Q_{3}-Q_{2}$
	\end{tabular}
	\caption{Plots of the velocity and density fields in addition to the mesh deformation for the three-dimensional Sedov test using $Q_{1}-Q_{0}$, $Q_{2}-Q_{1}$, $Q_{3}-Q_{2}$ velocity-energy pairs. Various viewpoints and cuts are presented.}
	\label{fig:CubeSedovresults}
\end{figure}

We consider a $[0, 1] \times [0, 1] \times [0,1]$ domain and a final time $t=0.8$. In Figure~\ref{fig:CubeSedovresults} we show plots of the velocity and density fields in different cross-sections at $t = 0.8$ for the $Q_{1}-Q_{0}$, $Q_{2}-Q_{1}$, $Q_{3}-Q_{2}$ velocity-energy pairs. 
The weak wall boundary conditions produce solutions
indistinguishable from those obtained with strong enforcement.

\subsection{Three-dimensional Sedov explosion in a cube with a spherical hole}
\label{cube_hole_sedov_3d}
We perform the Sedov test in a unit cube with a spherical hole and show plots of the velocity and density  fields in different cross-sections at the final time of $t = 0.8$ for the $Q_{1}-Q_{0}$, $Q_{2}-Q_{1}$, $Q_{3}-Q_{2}$ velocity-energy pairs as shown in Figure~\ref{fig:SquareHoleSedovresults}.
We observe similar behavior as in the corresponding two-dimensional tests, confirming that the method is directly applicable to simulations in complex three-dimensional domains.
\begin{figure}[tb]
	\centering
	\begin{tabular}{ccc}
		\multicolumn{3}{c}{Velocity}  \\
		\includegraphics[width=0.230\linewidth]{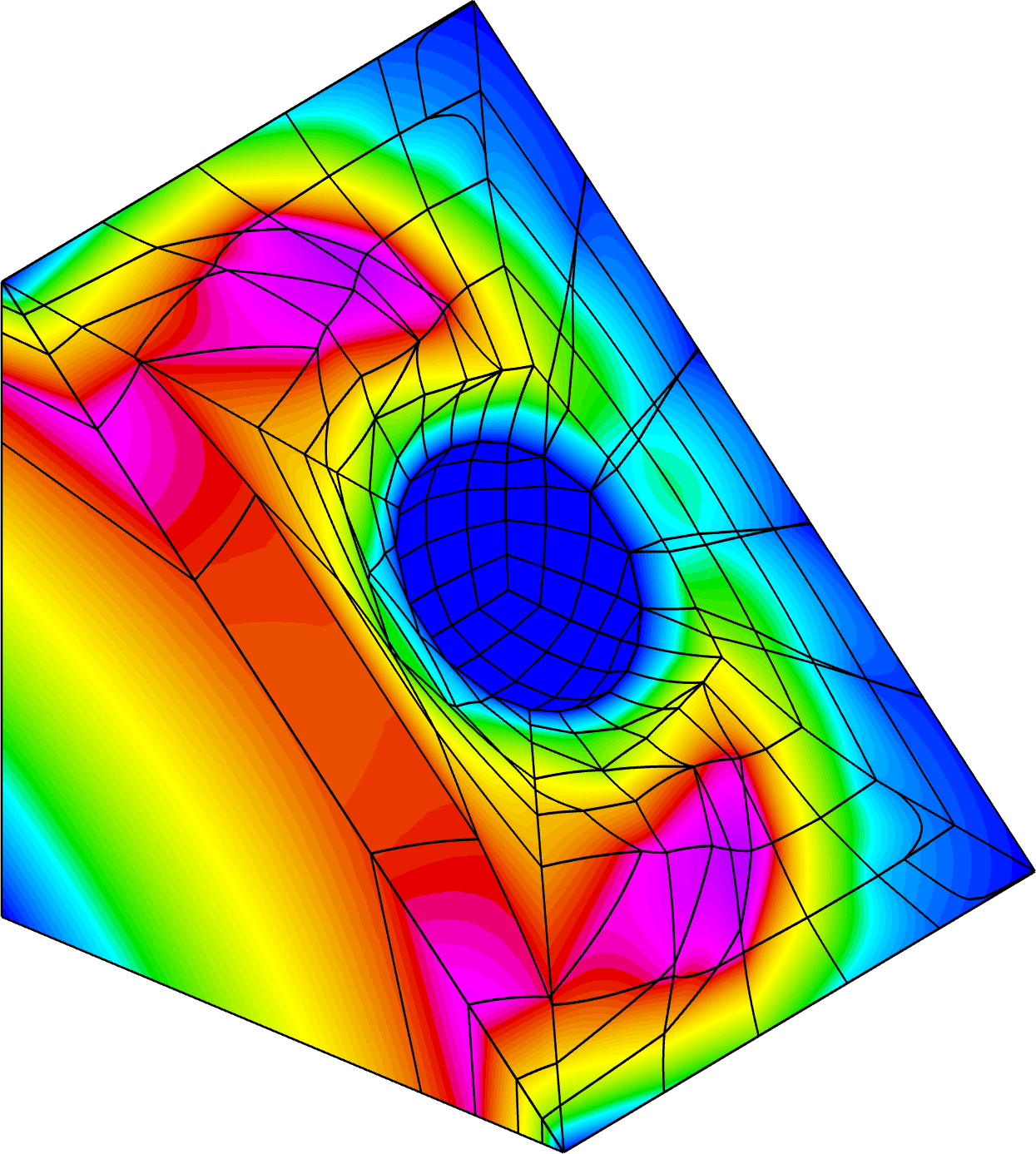} 
		\qquad & \qquad  
		\includegraphics[width=0.230\linewidth]{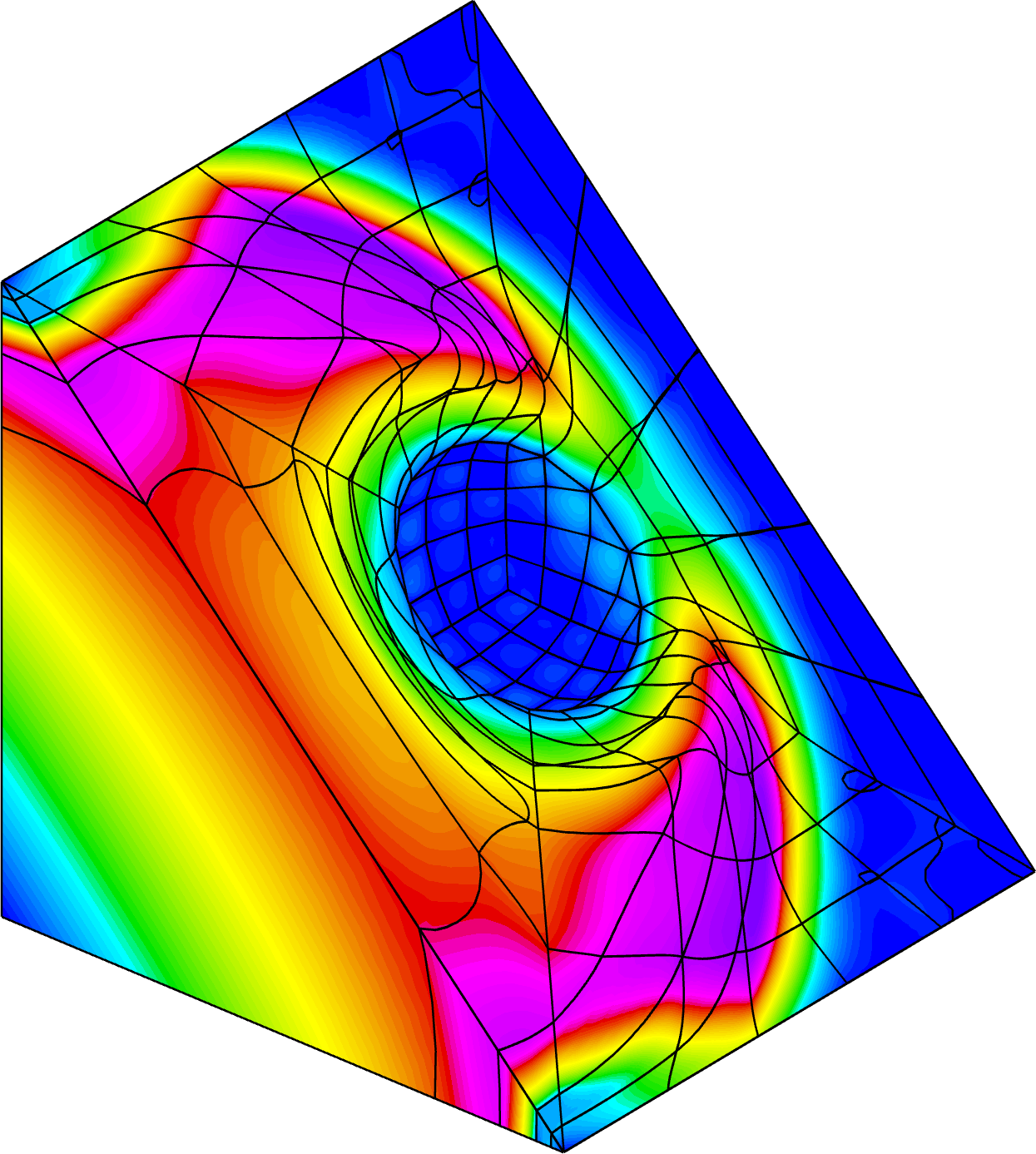}
		\qquad & \qquad  
		\includegraphics[width=0.230\linewidth]{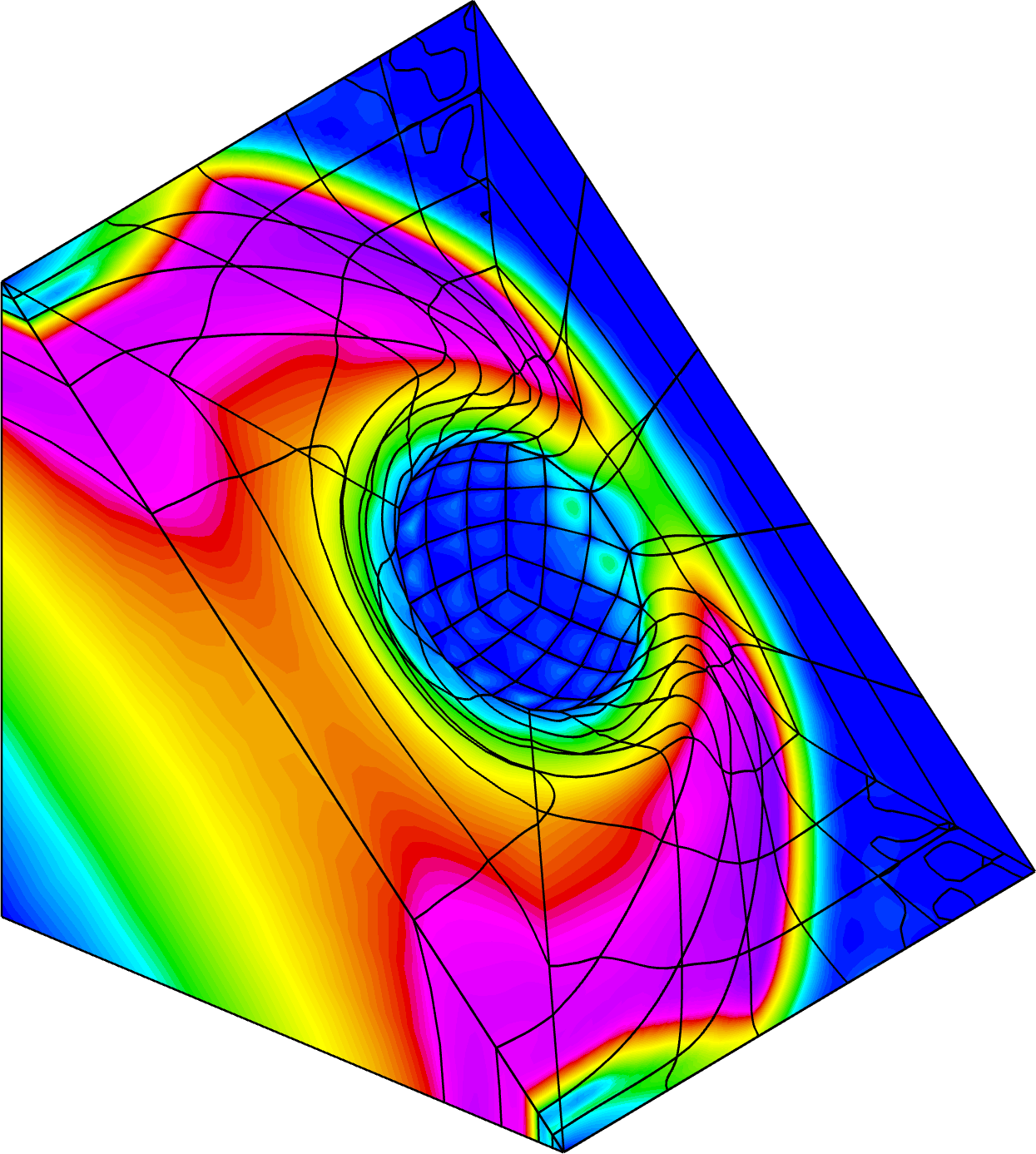} \\[.4cm]
		\includegraphics[width=0.230\linewidth]{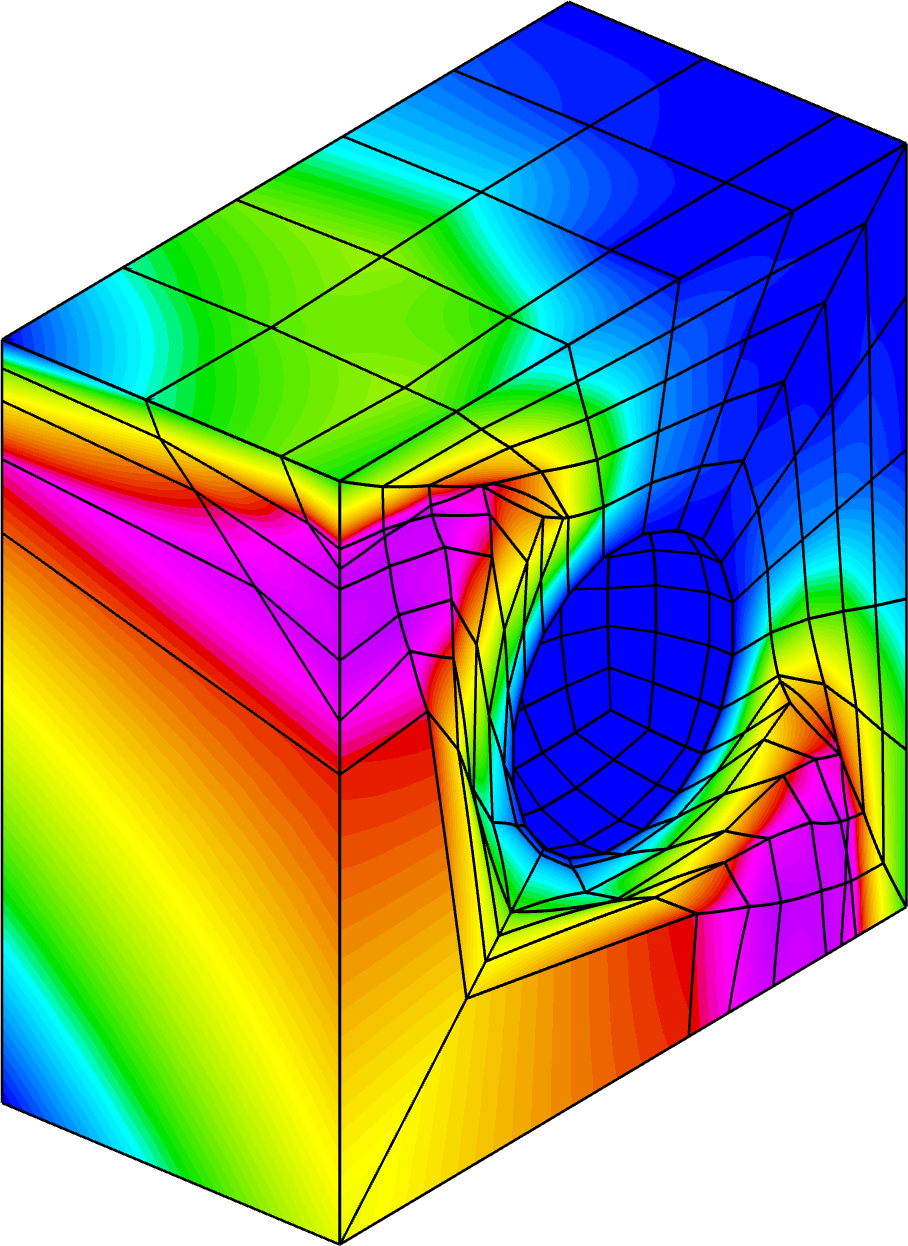} 
		\qquad & \qquad  
		\includegraphics[width=0.230\linewidth]{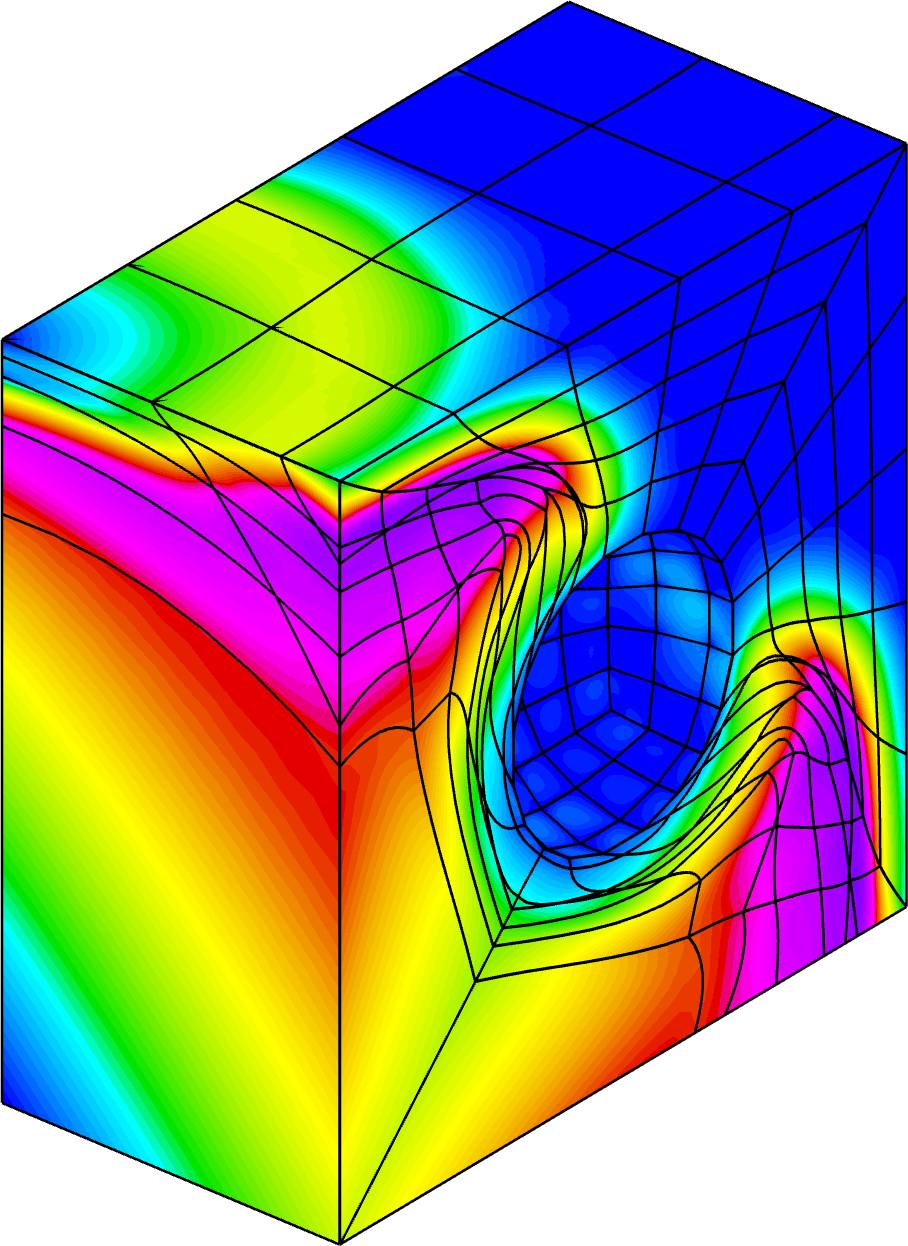}
		\qquad & \qquad  
		\includegraphics[width=0.230\linewidth]{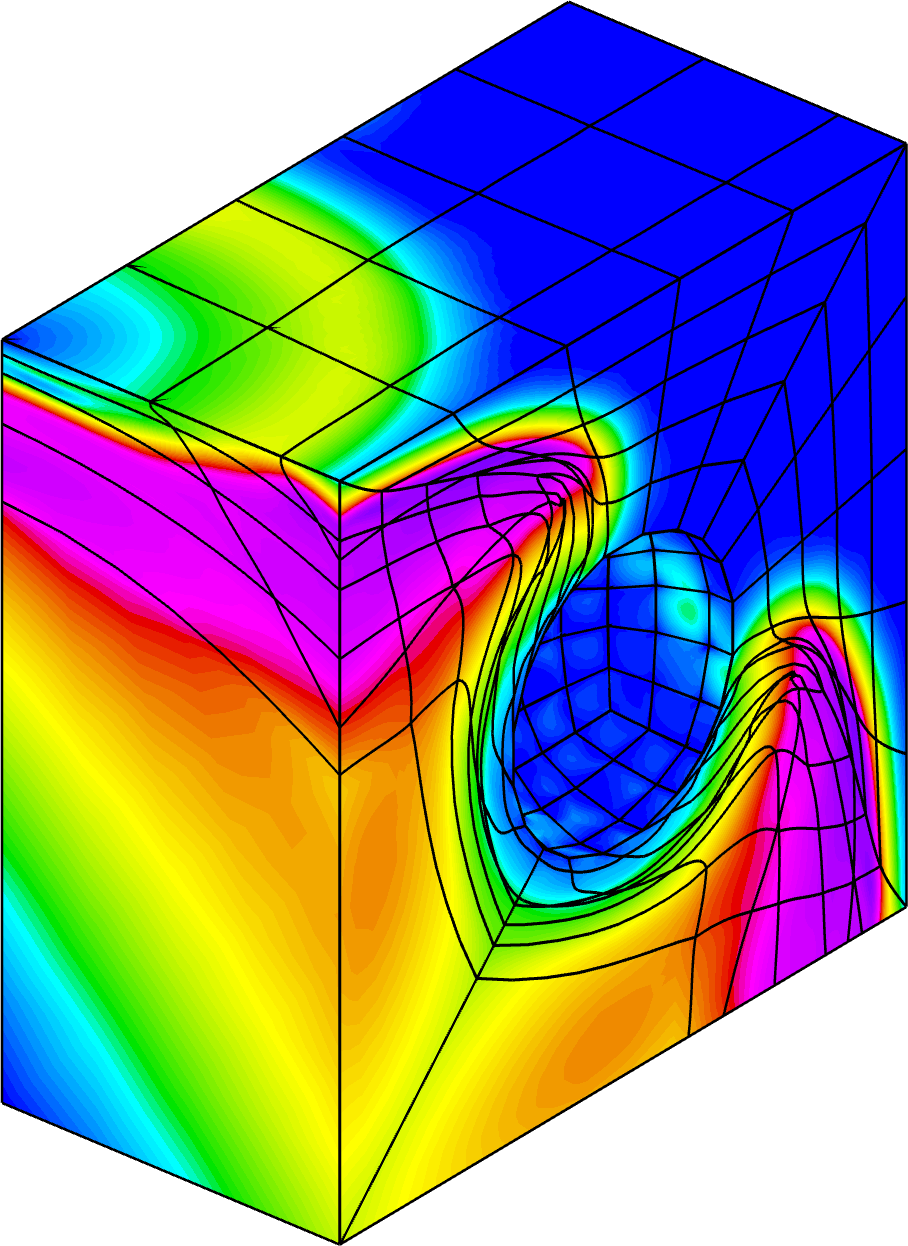}\\[.4cm]
		\multicolumn{3}{c}{Density}  \\
		\includegraphics[width=0.230\linewidth]{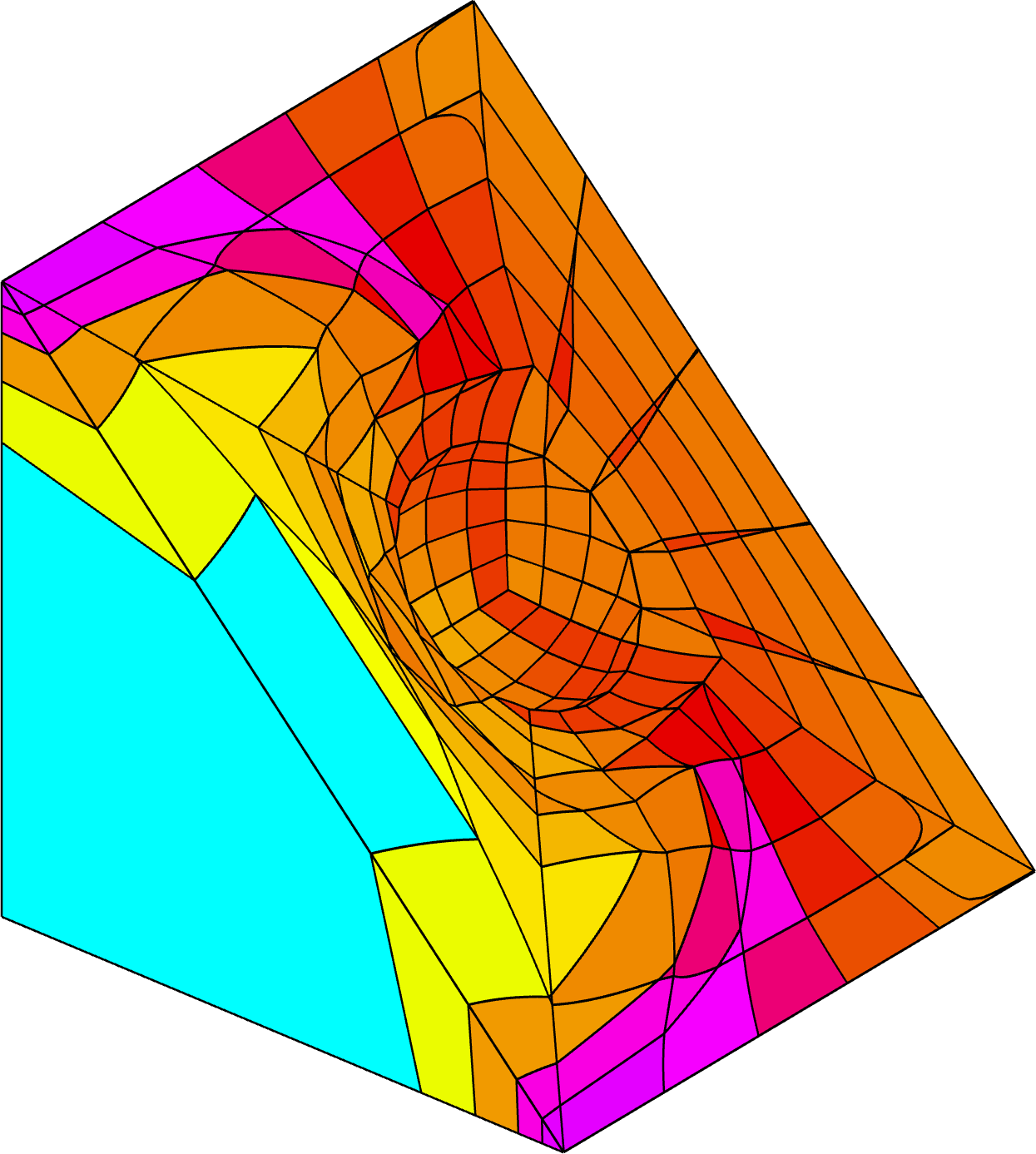} 
		\qquad & \qquad  
		\includegraphics[width=0.230\linewidth]{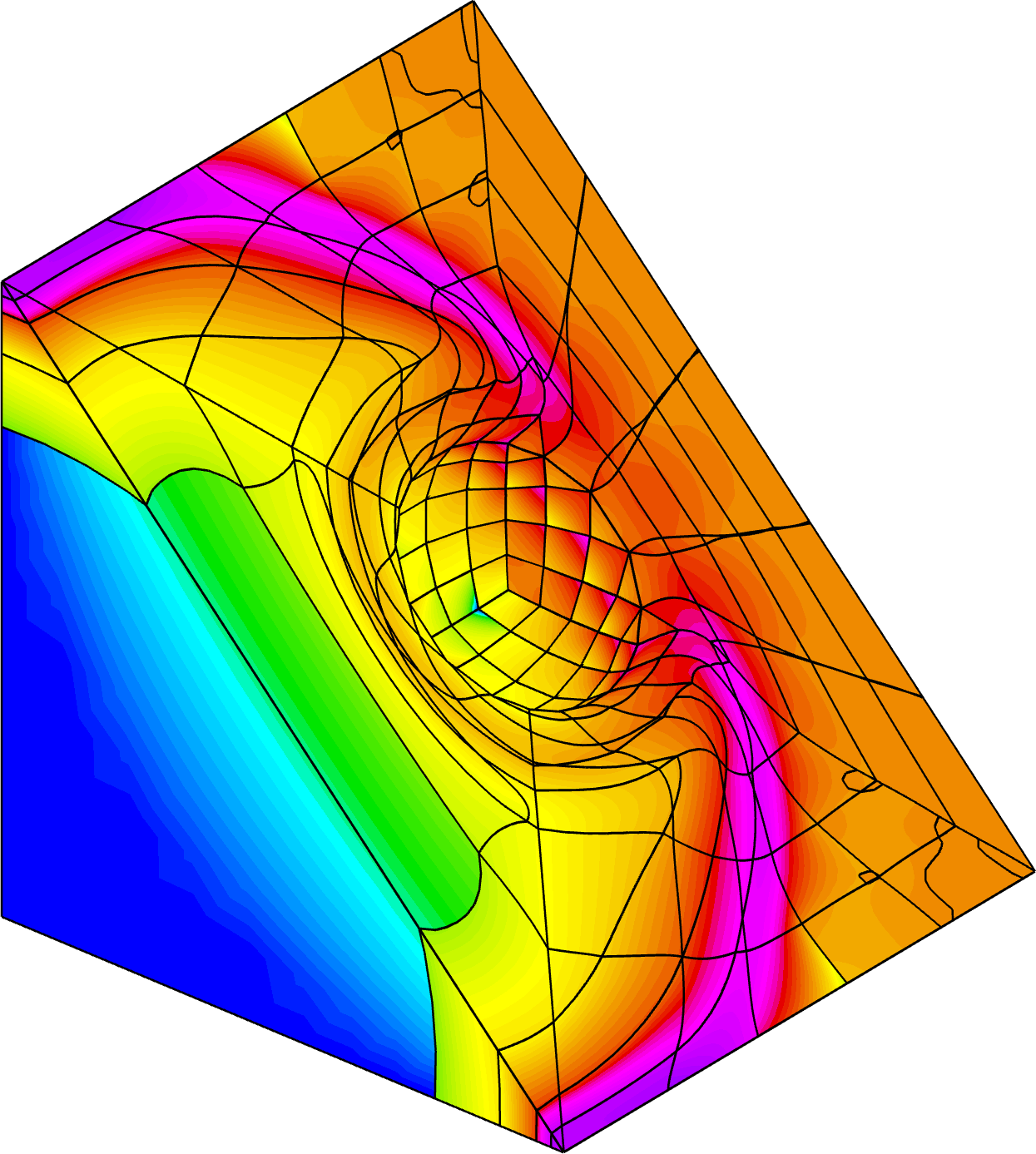}
		\qquad & \qquad  
		\includegraphics[width=0.230\linewidth]{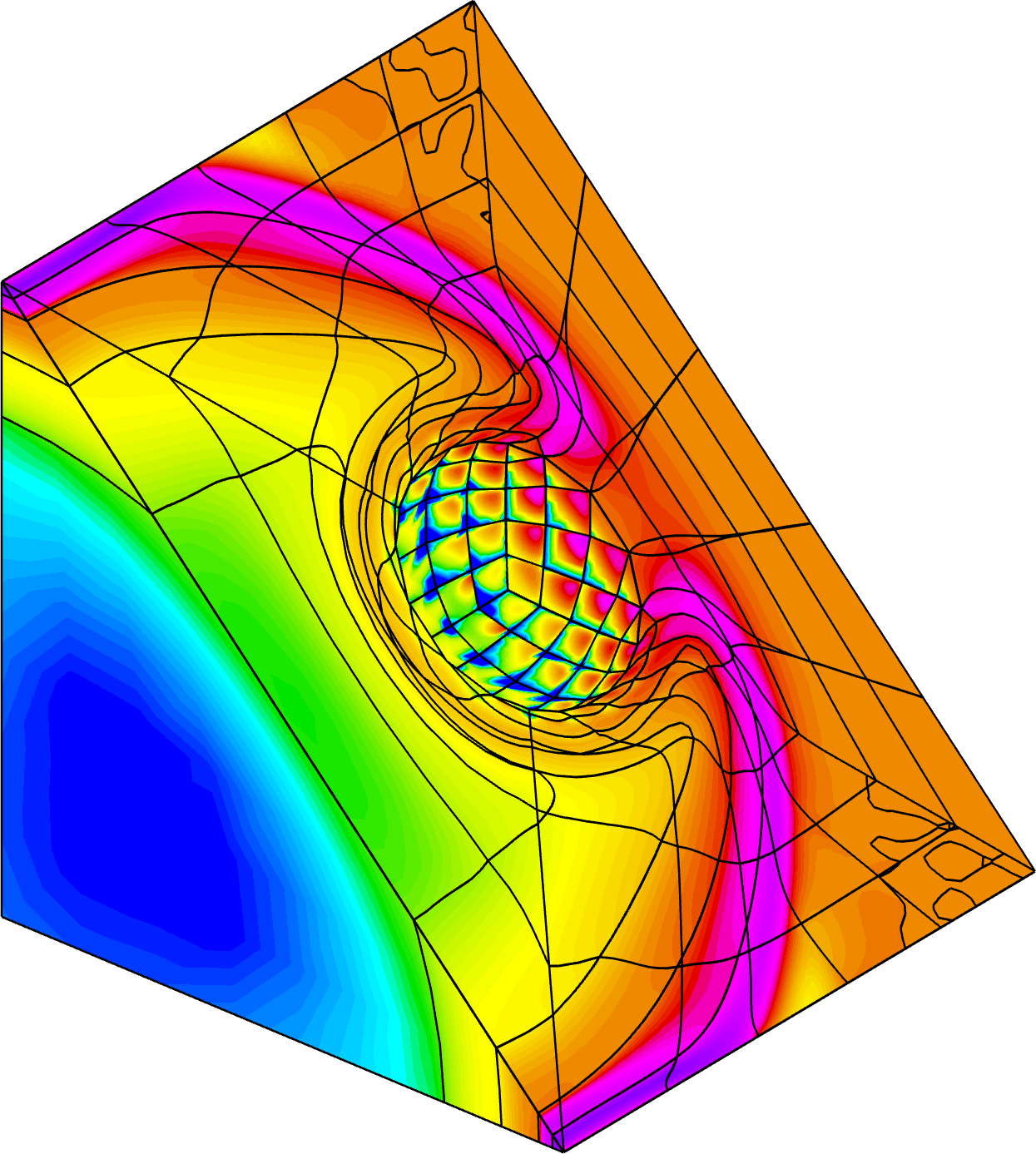}  \\[.4cm]
		\includegraphics[width=0.230\linewidth]{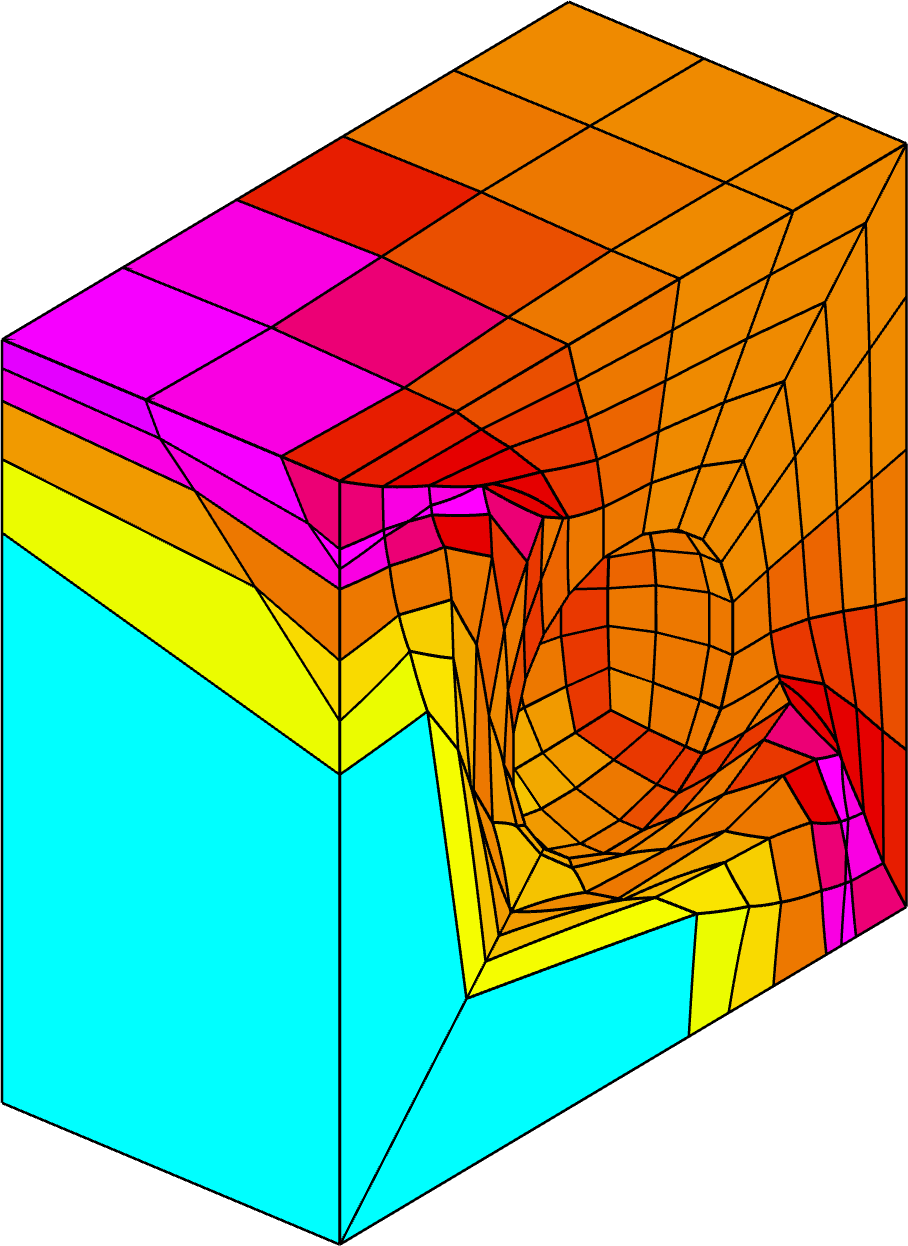} 
		\qquad & \qquad  
		\includegraphics[width=0.230\linewidth]{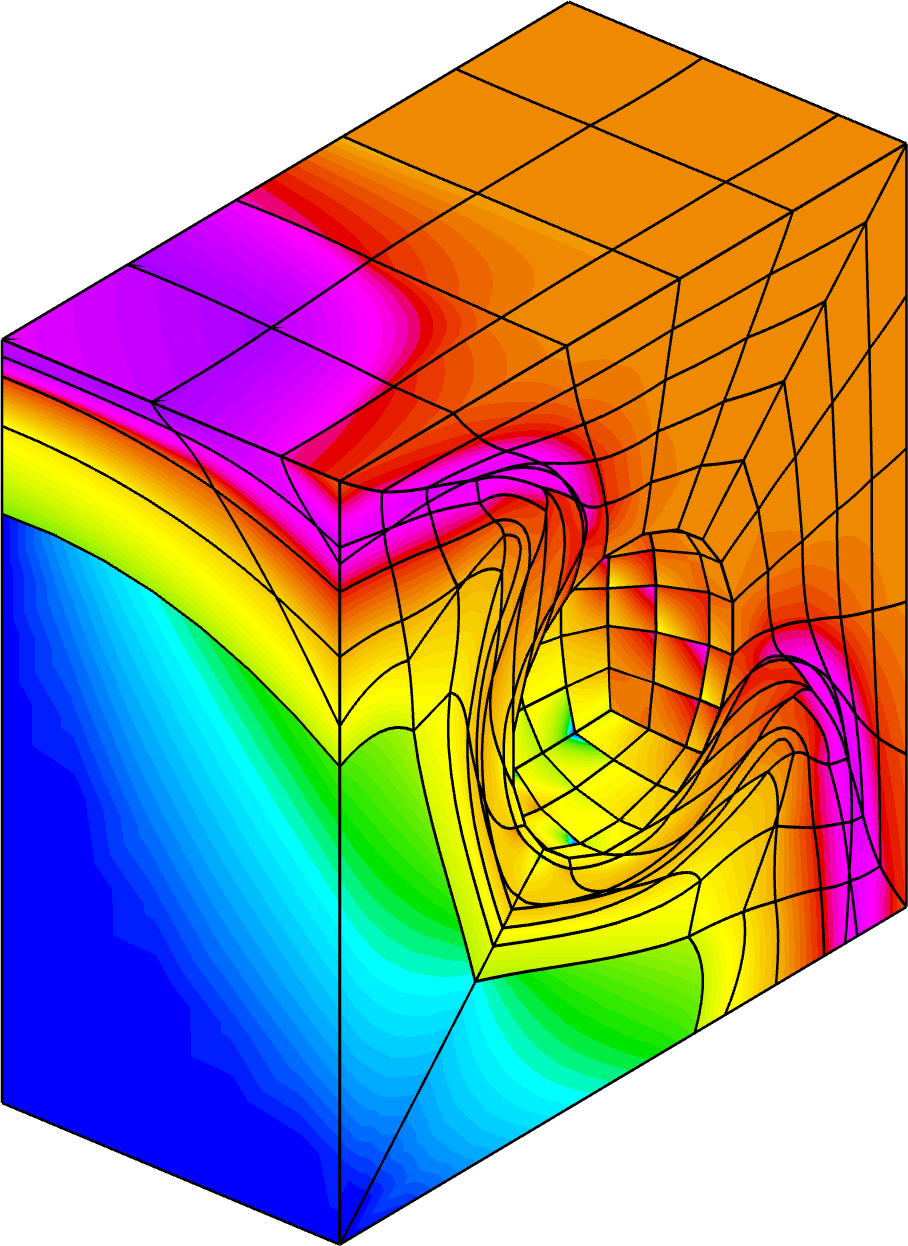}
		\qquad & \qquad  
		\includegraphics[width=0.230\linewidth]{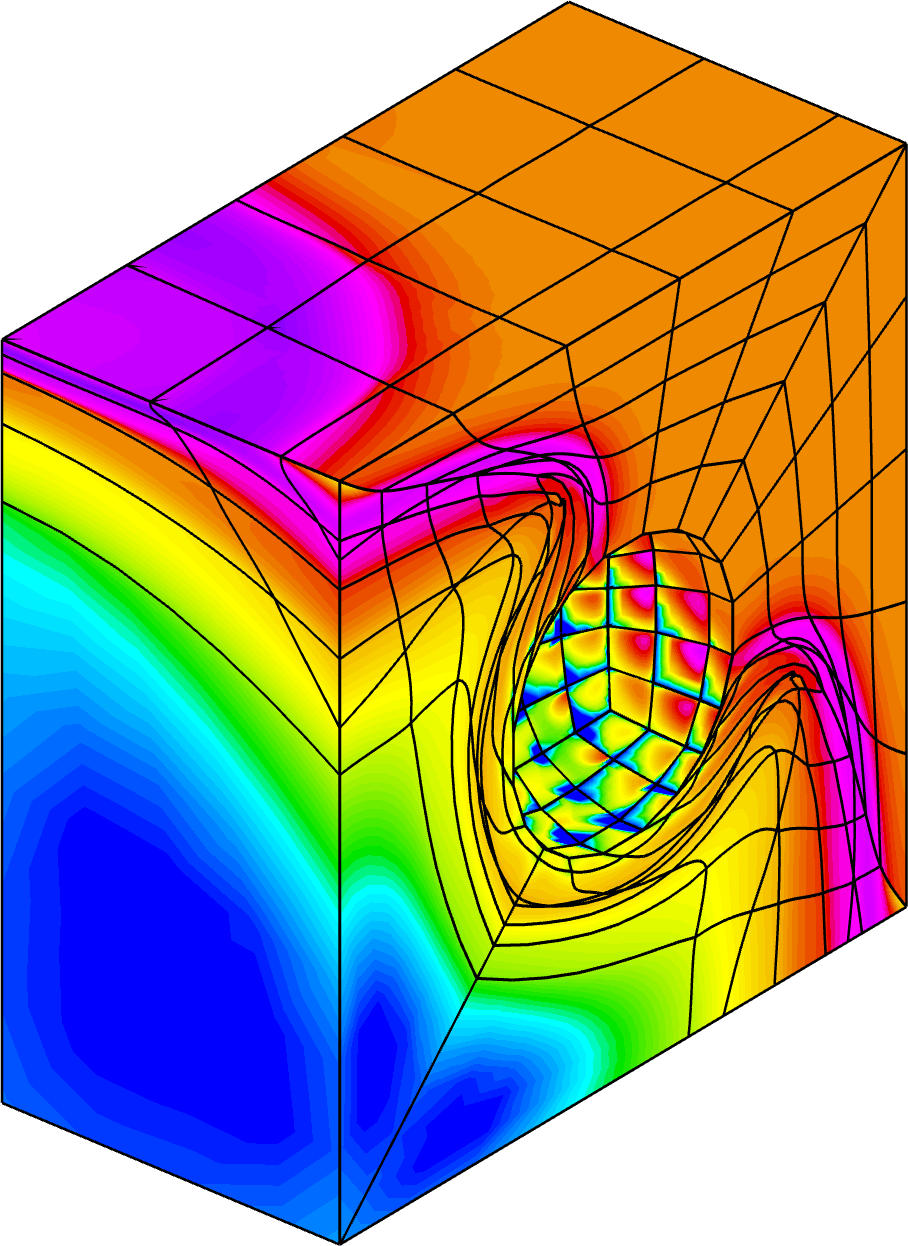} \\[0cm]
		$Q_{1}-Q_{0}$ 
		&
		$Q_{2}-Q_{1}$
		&
		$Q_{3}-Q_{2}$
	\end{tabular}
	\caption{Plots of the velocity and density fields in addition to the mesh deformation for the planar Sedov test using $Q_{1}-Q_{0}$, $Q_{2}-Q_{1}$, $Q_{3}-Q_{2}$ velocity-energy pairs. Various viewpoints and cuts are presented.}
	\label{fig:CubeHoleSedovresults}
\end{figure}